# Determination of Preferred Fiber Orientation State based on Newton – Raphson Method using on Exact Jacobian

Aigbe Awenlimobor [*], Douglas E. Smith

Department of Mechanical Engineering, School of Engineering and Computer Science, Baylor University, Waco, TX 76798, USA

*Corresponding author. Email: aigbe_awenlimobor1@baylor.edu

## Abstract

Fiber orientation is an important descriptor of the microstructure for short fiber polymer composite materials where accurate and efficient prediction of the orientation state is crucial when evaluating the bulk thermo-mechanical response of the material. Recent macroscopic fiber orientation models have employed the moment-tensor form in representing the fiber orientation state which all require a closure approximation for the higher order orientation tensors. In addition, various models have been developed to account for rotary diffusion due to fiber-fiber and fiber-matrix interactions which can now more accurately simulate the experimentally observed slow fiber kinematics in polymer composite processing. Traditionally explicit numerical IVP-ODE transient solvers like the 4th order Runge-Kutta method have been used to predict the steady-state fiber orientation state. Here we propose a computationally efficient method based on the Newton-Raphson iterative technique for determining steady state orientation tensor values by evaluating the exact derivatives of the moment-tensor evolution equation with respect to the independent components of the orientation tensor. We consider various existing macroscopic fiber orientation models and several closure ap-proximations to ensure the robustness and reliability of the method. The performance and stability of the approach for obtaining physical solutions in various homogeneous flow fields is demonstrated through several examples. Validation of the obtained exact derivatives of the orientation tensor is performed by benchmarking with results of finite difference techniques.

## Introduction

Characterization and evaluation of fiber suspensions has received considerable attention over the past four decades, particularly in the area of short fiber polymer composites produced in flow processes such as extrusion and molding. Understanding the flow of fiber suspensions is critical for predicting the thermal, mechanical, and electrical performance of products made of these materials. Fiber orientation within a fiber suspension is often simulated using fiber orientation tensors which quantify the degree of alignment. The goal of this paper is to compute the steady state orientation tensor of a fiber suspension in a homogeneous flow field using a Newton-Raphson iteration method with exact Jacobians.

Most models used to compute the fiber orientation, $p_i$ of a fiber suspension in polymer composite flows originated in the pioneering work by Jeffery [1] which described the evolution of a single rigid ellipsoid in a purely viscous flow of a Newtonian fluid. Figure 1 shows a typical orientation state of an ellipsoid and the coordinate systems and Euler angles definition.

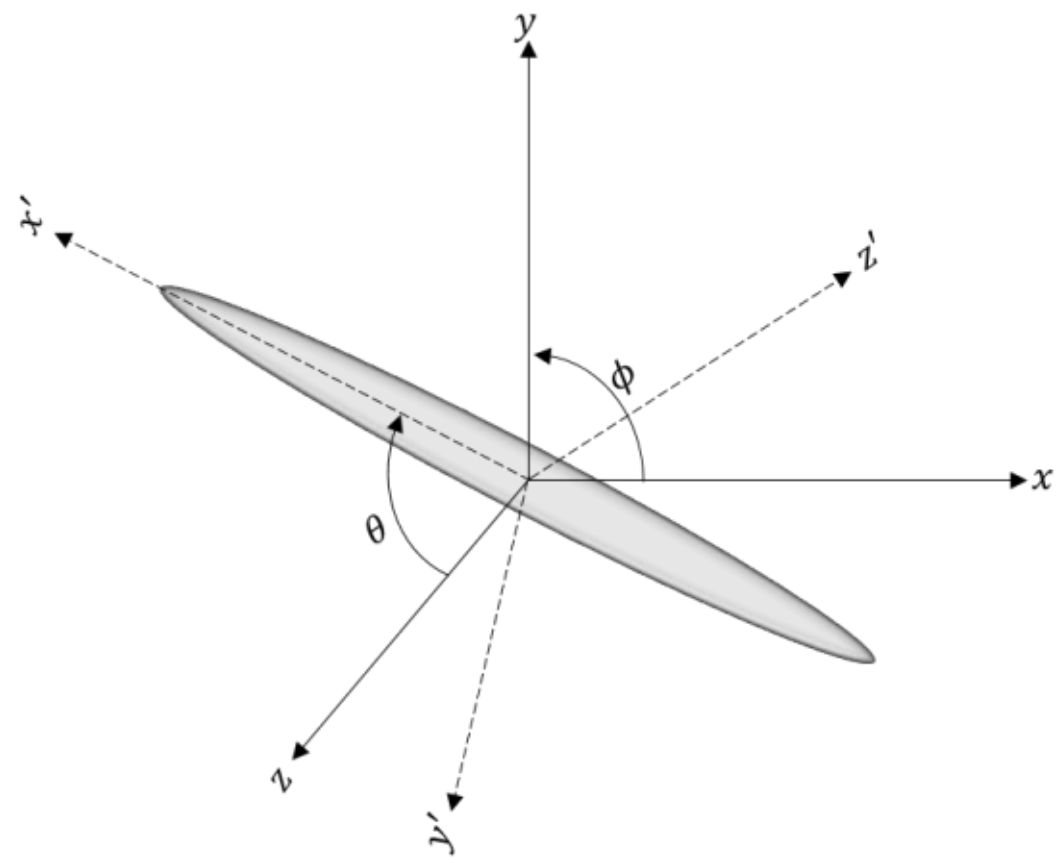

**Figure 1.** Single ‘rigid’ ellipsoidal fiber orientation

Jeffery’s hydrodynamic (HD) single ellipsoid model is limited to dilute suspension and ignores the effect of momentum diffusion due to fiber-fiber interactions. Moreover, Jeffery’s assumption ignores the flexural and fracture behavior of the fiber and assumes no slip contact between the fiber and surrounding fluid. Studies have shown that fibers suspended in industrial polymer melt flows tend to orient relative to the flow field [2] which is not captured in Jeffery’s work. As a result, enhancements to the Jeffery’s single fiber model have been made to better capture the bulk behavior of fibers in semi-dilute and concentrated suspension which include fiber-fiber interaction.

Although theoretically feasible, it is computationally expensive and nearly impractical to simulate the behavior of every individual particle in an industrially relevant short fiber polymer suspension flow. Folgar and Tucker [3], [4] introduced a phenomenological isotropic rotary diffusion (IRD) term with a linear dependence on the scalar magnitude of the rate of deformation tensor and derived from an orientation probability distribution function (ODF) and the hydrodynamic contribution from Jeffery’s model. The time evolution of the ODF is defined by the Fokker-Plank equation for probability distribution function (PDF) of fiber orientation. Conventionally a numerical method such as finite volume (cf. Bay [5]) and more recently a more computationally efficient exact spherical harmonics method (cf. Montgomery-Smith et al. [6]) has been used to solve the Folgar-Tuckers (FT) equation of change for fiber orientation, but these methods have yet to see significant application in molding or extrusion processes. Advani and Tucker [7] proposed a tensorial representation of fiber orientation evolution based on the series expansion of even order moments of the ODF. Their approach simplifies computations which has led to the widespread use of orientation tensors as the preferred method of evaluating fiber orientation for short fiber polymer composites. The Advani-Tucker nth order orientation evolution model commonly used to compute the 2nd order orientation tensor, $a_{ij}$ thus requires a closure approximation. Due to experimentally observed disparity in the fiber orientation kinetics compared to those computed from traditional orientation models, modifications to the Advani-Tucker model have been proposed in an effort to reduce the rate of alignment in polymer melt flow. Huynh [8] applied a strain reduction factor (SRF) to the Advani-Tucker’s model to slow the transient response of the orientation tensor. Unfortunately, Huynh’s SRF model lacked material objectivity. To better address material objectivity, Wang et al. [9] developed a phenomenological reduced strain closure (RSC) model that applies the reduction factor solely to the spectral decomposed principal rates of the orientation tensor, without modifying the evolution of the rotation tensor. Similarly, Tseng et al. [10], [11] proposed a retarding principal rate (RPR) model that involved a coaxial correction to the FT model by assuming the intrinsic orientation kinetics (IOK) describing the behavior of the fiber suspension involved a nonlinear modification to the principal directions of the material derivative.

Prediction of fiber orientation with IRD-based models have been validated experimentally for short-fiber/thermoplastic composites (SFT) with fiber length typically in the range of 0.2mm to 0.4mm [12]. For long-fiber/thermoplastic composites (LFT) with typical size above 10mm, IRD models possess unidirectional prediction effectiveness. Modifications to the rotary diffusion term have been proposed for improving the accuracy when predicting the components of the orientation tensor. Ranganathan et al. [13] assumed an isotropic rotary diffusivity that inversely varies with the degree of alignment of the orientation tensor which was implemented using a phenomenological interaction parameter that depends on the reciprocal of the inter fiber spacing. The applicability of their model is limited to the transient orientation state while being well suited for long range fiber-fiber interaction. However, their model was shown to be unsuitable for LFTs steady state orientation prediction as with other IRD models since its diffusivity was isotropic.

Fan et al. [14] and Phan-Thien [15], were the first to propose an anisotropic rotary diffusion (ARD) moment-tensor model by replacing the scalar phenomenological interaction parameter with a second order anisotropic rotary diffusion tensor. In a similar manner, Koch [16] developed an ARD model suited for semi-dilute suspensions with an anisotropic spatial diffusion tensor that depends on the orientation state and the rate of deformation tensor. However, their model was based on the more complicated PDF form for the orientation tensor representation rather than the moment-tensor form and proved ineffective in LFT modelling. Phelps et al. [12] built on the work of Fan [14] and Phan-Thien et al. [15] by developing a more general moment-tensor anisotropic diffusion model that depends on the spatial diffusion tensor and orientation tensor state. The derivation of the spatial diffusion tensor was written as a function of the orientation state and rate of deformation tensor in a similar manner to that proposed by Hand [17]. Phelps's model showed remarkable improvements in predicting orientation states of LFTs. Tseng et al. [18] proposed an improved anisotropic rotary diffusion model (iARD) which defines a two-parameter spatial diffusion tensor that couples the effect of fiber-matrix interaction and fiber-fiber interaction. Unfortunately, the iARD model lacked material objectivity, limiting its applicability. More recently, Tseng et al. [19] proposed a principal anisotropic rotary diffusion (pARD) model assuming a principal spatial diffusion tensor that corotates with the orientation tensor. Bakharev [20] proposed a moldflow rotational diffusion (MRD) model based on a reduction of the terms from the generic ARD model by Phelps to linear terms only with a spatial diffusion tensor like Tseng's model. Latz et al. [21] developed a fully coupled flow-orientation model for concentrated suspensions by replacing the diffusion term in the FT model with an effective collision tensor that incorporates both isotropic diffusion interaction term and a topological exclusion interaction term based on a nematic (NEM) 'Onsager' potential of non-Brownian Maier-Saupe form. They found the influence of the topological interaction on the fiber orientation state to be flow dependent having a significant effect on channel and contraction flows and a relatively lesser influence for flow around cylinders. Kugler et al. [22], Favaloro et al. [23], Agboola et al. [24] and Park and Park [25] presents detailed review and comparison of existing fiber orientation models. The foregoing ARD models find useful application in polymer composite industry and have been incorporated in mold-filling flow computations in injection molding process simulations [26], [27], [28], [29], [30], [31], [32].

Due to the absence of exact solutions for orientation state for inhomogeneous flows involving momentum diffusion, various closure approximations with different degrees of accuracy have been proposed to calculate higher order fiber orientation moment-tensors. These closure approximations are developed form lower order orientation tensors and their identity tensors and can be subdivided into non-fitted and fitted closures. The non-fitted closure approximations include the those constructed from combinations of lower order orientation tensors and their identity matrices such as the linear (LIN), quadratic (QDR), and the hybrid (HYB) closures. The general class of Hinch and Leal's [33] composite closure approximations precontracted with the deformation rate tensor also falls under this category. The class of fitted closures usually involve deriving analytical expressions for the independent components of the higher order tensor based on a polynomial fitting procedure using test data obtained from solutions to the probability distribution function (PDF) for different standard flow conditions. The fitted closures can be sub-divided into the Eigenvalue-Based Fitted (EBF) and Invariant-Based Fitted (IBF) closures. Higher order optimal fitted closures such as the Eigenvalue-Based Optimal Fitted (EBOF) and Invariant-Based Optimal Fitted (IBOF) closures which are extensions to the EBF and IBF closures respectively have also been developed for improved model accuracy. More details on different closure approximations and their associated advantages and disadvantages is discussed in later sections and also in [22], [34], [35], [36]. Other closure approximations include the neural network based fitted closures by Jack et al. [37] and the 6th order Invariant based orthotropic fitted closure by Jack and Smith [38], [39].

The steady state orientation tensor values have traditionally been computed with time evolution numerical IVP-ODE techniques such as the 4th order Runge-Kutta method or predictor-corrector methods. Here we present a computationally efficient and faster method based on Newton-Raphson iteration for determining the steady state or preferred orientation using explicit derivatives of the 2nd order moment-tensor equation of change with respect to its orientation tensor components. Here we consider various fiber orientation models and closure approximations and compare their performance in complex homogeneous flow fields. We benchmark the results of the explicit derivatives with those obtained using finite differences to ensure accuracy. The explicit derivatives are comparatively faster compared to the finite difference derivatives.

## Methodology

### *Determination of Steady State Orientations*

The numerical approach developed here for determining the steady state orientation vector $p_i$ and orientation tensor components $\mathrm{a}_{ij}$ computes a zero rate of change of the orientation state using the Newton-Raphson iteration by setting the rate of change equation (often referred to as the residual, $R_i$) to zero, respectively, as

$$R_i = \frac{Dp_i}{Dt} = 0 \qquad and \qquad R_{ij} = \frac{D\mathrm{a}_{ij}}{Dt} = 0 \tag{1}$$

based on Newton's algorithm, the successive iterative improvement to the approximation of a given root (in our case the orientation state) is computed as [40]

$$p_i^+ = p_i^- - J_{ij}^{-1} R_j \tag{2}$$

for the orientation vector and

$$\mathrm{a}_{ij}{}^{+} = \mathrm{a}_{ij}{}^{-} - J_{mnij}{}^{-1} R_{mn} \tag{3}$$

for the 2nd order orientation tensor where the $(+)$ and $(-)$ superscripts denote the current and previous iterations. The implication of Equations (2) and (3) is the need to compute the Jacobians $J_{ij}$ and $J_{mnij}$ which are the derivatives of the orientation vector and tensor, respectively, with respect to its components that define them i.e.

$$J_{ij} = \frac{\partial R_i}{\partial p_j} \qquad and \qquad J_{mnij} = \frac{\partial R_{mn}}{\partial \mathrm{a}_{ij}} \tag{4}$$

In the following sections, we present existing models for rate of change equations of the orientation vector and tensor in Equation (1) based on a review by Kugler et al. [22]. Equating each rate of change equation to zero yields the Newton-Raphson residual which is then differentiated to obtain expressions for the associated Jacobian for each model.

Since there are only 5 independent components of $a_{ij}$ (see discussion below), the residual $R_{mn}$ may be represented in contracted form as the vector $R_r$ (and similarly, the tensor $a_{ij}$ may be represented as the vector $a_r$), and the Jacobian $J_{mnij}$ as a matrix $J_{rs}$. To relate form the contracted notation, we employ the index mapping

$$r(m,n) = n - \frac{1}{2}(m-1)(m-6), \qquad for\; n = m \dots 3, \qquad for\; m = 1, 2 \tag{5}$$

It follows that Equation 3 in contracted notation becomes

$$a_r{}^{+} = a_r{}^{-} - J_{rs}{}^{-1} R_r \tag{6}$$

*Fiber Orientation*

The orientation state of a fiber as shown Figure 1 in can be described by the unit vector $\underline{p}$ associated as [7]

$$\underline{p} = [\sin\theta\cos\phi \quad \sin\theta\sin\phi \quad \cos\theta]^T \tag{7}$$

The orientation of the unit vector $\underline{p}$ is commonly described in terms of the probability distribution function (PDF) $\psi\left(\underline{p}\right)$ over all possible directions of $\underline{p}$. It can be shown that $\psi\left(\underline{p}\right) = \psi\left(-\underline{p}\right)$ satisfies the normalization condition

$$\oint \psi\left(\underline{p}\right) d\underline{p} = \int_{\theta=0}^{\pi} \int_{\phi=0}^{2\pi} \psi(\theta,\phi) \sin\theta \, d\theta d\phi = 1 \tag{8}$$

The PDF $\psi\left(\underline{p}\right)$ also satisfies the continuity condition [7]

$$\frac{D\psi}{Dt} = -\frac{\partial}{\partial \underline{p}}\left(\psi \underline{\dot{p}}\right) \tag{9}$$

*Fiber Orientation Modelling*

Macroscopic fiber Orientation modelling is usually required in polymers processing to predict the bulk response of chopped fibers in composites parts and ultimately determine part performance. The choice of macroscopic model depends on various processing parameters such as the concentration of fiber suspension, flow type and strength, fiber geometry and volume fraction, material rheology, etc. Fiber suspension concentration is classified into 3 regimes depending on the fiber volume fraction $\phi_f = nV_f$, $n$ is the number of fibers per unit volume and $V_f$ is the fiber volume. Depending on the degree of fiber alignment and fiber's length, the various suspension regimes are [39], [41], [42]

$$\begin{cases} \phi_f < \frac{1}{r_e^2}, & dilute \\ \frac{1}{r_e^2} \le \phi_f < \frac{1}{r_e}, & semi-concentrated \\ \phi_f \ge \frac{1}{r_e}, & concentrated \end{cases} \tag{10}$$

*Fiber Orientation Modelling in the Dilute Regime*

Jeffery's hydrodynamic model for the motion of a single rigid ellipsoidal particle in an incompressible Newtonian viscous fluid flow field forms the basis for evaluating fiber orientation for dilute suspensions[1]. Jeffery assumed that an ellipsoidal particle is convected with the bulk motion of the undisturbed surrounding fluid where the components of $\underline{p}$ are given as a function of time by [1], [22], [23].

$$\dot{p}_i^{JF} = \omega_{ij}p_j + \xi\left(\dot{\gamma}_{ij}p_j - \dot{\gamma}_{kl}p_k p_l p_i\right) \tag{11}$$

where the superscript $JF$ represents 'Jeffery' type particle, and $\omega_{ij}$ and $\dot{\gamma}_{ij}$ are the anti-symmetric and symmetric decomposition of the rate of deformation tensor $L_{ij} = {\partial v_i}/{\partial x_j}$ given respectively as

$$\omega_{ij} = \frac{1}{2}\left(\frac{\partial v_i}{\partial x_j} - \frac{\partial v_j}{\partial x_i}\right), \qquad \dot{\gamma}_{ij} = \frac{1}{2}\left(\frac{\partial v_i}{\partial x_j} + \frac{\partial v_j}{\partial x_i}\right) \tag{12}$$

such that $L_{ij} = \partial v_i/\partial x_j = \dot{\gamma}_{ij} + \omega_{ij}$, and $\xi$ is a particle shape parameter given as $\xi = (r_e^2 - 1)/(r_e^2 + 1)$ with particle aspect geometric ratio $r_e$. We define components $R_i$ of the Newton-Raphson residual in Equation (2) for computing a steady-state orientation vector $\underline{p}$ in Equation (9) as

$$R_i^{JF} = \dot{p}_i^{JF} \tag{13}$$

where the Jacobian $J_{mn}^{JF} = {\partial R_m}/{\partial p_n}$ is obtained by taking derivatives of the components $R_m$ with respect to $p_n$, i.e.

$$J_{mn}^{JF} = \frac{\partial \dot{p}_m^{JF}}{\partial p_n} = \omega_{mj}\delta_{jn} + \xi\left[\dot{\gamma}_{mj}\delta_{jn} - \dot{\gamma}_{kl}(\delta_{kn}p_l p_m + p_k\delta_{ln}p_m + p_k p_l \delta_{mn})\right] \tag{14}$$

where we note that derivatives of $\underline{p}$ with respect to itself form the identity matrix., i.e., $\partial p_i/\partial p_j = \delta_{ij}$. Equation (14) contains only 2 independent components of the orientation vector $R_i^{JF}$ making $J_{mn}^{JF}$ a $2 \times 2$ matrix. Jeffery's model for dilute suspensions has limited direct application since fiber-fiber interaction is ignored, requiring that more advanced models be employed to capture these effects.

*Fiber Orientation Modelling in Semi-dilute and Concentrated Regime*

In the concentrated regime, the average interparticle spacing is very small (orders of magnitude less than the smallest particle dimension) such that the fiber motion is affected by hydrodynamic forces and possible direct mechanical contact with other particles and physical boundaries. These semi-dilute and concentrated fiber suspension categories are defined by the second and third inequalities of Equation (10) respectively. The effect of the interparticle interaction on the single particle's motion is typically modeled by the incorporation of a momentum diffusion term to the equation defining the free particle's motion in dilute suspension. Different diffusion models with varying degree of accuracy are presented in subsections following. Most commercial SFRP composites suspensions fall within the concentrated class of fiber suspension.

1. The Folgar-Tucker Diffusion Model

A model that accounts for the effect of momentum diffusion due to short- and long-range fiber-fiber interaction in non-dilute fiber suspensions was first proposed by Folgar-Tucker model [3], [4] by incorporating a rotary diffusion term into Jeffery's single fiber model as

$$\dot{p}_i^{FT} = \dot{p}_i^{JF} - D_r \frac{1}{\psi} \frac{\partial \psi}{\partial p_i} \tag{15}$$

where $D_r$ is the rotary diffusivity constant that introduces a Brownian diffusion effect among contacting particles. For long slender particles (i.e., $\xi \approx 1$) Folgar and Tucker set $D_r = C_I \dot{\gamma}$ where $C_I$ is a phenomenological interaction coefficient and $\dot{\gamma} = \sqrt{2\dot{\gamma}_{ij}\dot{\gamma}_{ji}}$ is the scalar magnitude of the strain rate tensor $\dot{\gamma}_{ij}$. The orientation PDF $\psi$ in Equation (15) defines the probability of a given fiber has a particular orientation and the rate of change of $\psi$ is given by the Fokker-Planck's continuity equation describing its time evolution as

$$\frac{D\psi}{Dt} = -\frac{\partial}{\partial p_i}(\psi \dot{p}_i) \tag{16}$$

---

[1] Dilute suspension is a heterogenous suspension where the average interparticle spacing is relatively large such that there is no restriction on the fibers motion due to hydrodynamic forces or mechanical contact, and typically defined as one that satisfies the first inequality of Equation (10)

The above PDF $\psi$ form of Folger-Tuckers model leads to a computationally intensive solution which limits its usefulness in real application. Advani and Tucker [7] derived the moment-tensor form of the Folger-Tuckers model by defining a set of even order orientation tensors as

$$\mathrm{a}_{ij} = \oint p_i p_j \psi\left(\underline{p}\right) d\underline{p}, \qquad \mathrm{a}_{ijkl} = \oint p_i p_j p_k p_l \psi\left(\underline{p}\right) d\underline{p} \tag{17}$$

which are, respectively, the 2nd and 4th order orientation tensors that quantify the orientation state of a fiber suspension. Orientation tensors defined in this form are completely symmetric i.e.

$$\begin{gathered} \mathrm{a}_{ij} = \mathrm{a}_{ji} \\ \mathrm{a}_{ijkl} = \mathrm{a}_{jikl} = \mathrm{a}_{kijl} = \mathrm{a}_{lijk} = \mathrm{a}_{ikjl} = \mathrm{a}_{iljk} = \cdots, \qquad 24\, permutations \end{gathered} \tag{18}$$

where a normalization condition requires that

$$\mathrm{a}_{ii} = 1, \qquad \mathrm{a}_{ijkk} = \mathrm{a}_{ij} \tag{19}$$

Consequently, there are only 5 independent components in $\mathrm{a}_{ij}$ and 9 independent components in $\mathrm{a}_{ijkl}$. Advani and Tucker developed an equation of change for the 2nd order orientation tensors, identified here and below as the orientation material derivative tensor $\dot{\mathrm{a}}_{ij}$, in terms of the 2nd and 4th order tensor as

$$\frac{D\mathrm{a}_{ij}}{Dt} = \dot{\mathrm{a}}_{ij}^{FT} = \left\{\dot{\mathrm{a}}_{ij}^{HD} + \dot{\mathrm{a}}_{ij}^{IRD}\right\} \tag{20}$$

where $\dot{\mathrm{a}}_{ij}^{HD}$ is the hydrodynamic tensor component of the Folger-Tuckers that represents Jeffery's equation and given as

$$\dot{\mathrm{a}}_{ij}^{HD} = -\left(\omega_{ik}\mathrm{a}_{kj} - \mathrm{a}_{ik}\omega_{kj}\right) + \xi\left(\dot{\gamma}_{ik}\mathrm{a}_{kj} - \mathrm{a}_{ik}\dot{\gamma}_{kj} - 2\dot{\gamma}_{kl}\mathrm{a}_{ijkl}\right) \tag{21}$$

and $\dot{\mathrm{a}}_{mn}^{IRD}$ is the isotropic rotary diffusion term modelling fiber interaction written as

$$\dot{\mathrm{a}}_{ij}^{IRD} = 2D_r\left(\delta_{ij} - \alpha \mathrm{a}_{ij}\right) \tag{22}$$

In the above, $\alpha$ is a dimension factor (i.e., $\alpha = 3$ for 3D orientation and $\alpha = 2$ for 2D planar orientation). It follows that the residual in Equation (2) for the Folger-Tuckers model is

$$R_{mn}^{FT} = \dot{\mathrm{a}}_{ij}^{FT} \tag{23}$$

The associated Jacobian $J_{mnij}^{FT}$ for the Folgar-Tucker model is obtained by differentiating $R_{mn}^{FT}$with-respect-to the components of $\mathrm{a}_{ij}$ yielding

$$J_{mnij}^{FT} = \frac{\partial R_{mn}^{FT}}{\partial \mathrm{a}_{ij}} = \frac{\partial \dot{\mathrm{a}}_{mn}^{HD}}{\partial \mathrm{a}_{ij}} + \frac{\partial \dot{\mathrm{a}}_{mn}^{IRD}}{\partial \mathrm{a}_{ij}} \tag{24}$$

Given that the derivative of $\mathrm{a}_{ij}$with respect itself is

$$\frac{\partial \mathrm{a}_{mn}}{\partial \mathrm{a}_{ij}} = \delta_{mi}\delta_{nj} \tag{25}$$

the derivative terms on the right-hand side of Equation (24) may be written as

$$\frac{\partial \dot{a}_{mn}^{HD}}{\partial a_{ij}} = (-\omega_{mk} + \xi\dot{\gamma}_{mk})\delta_{ki}\delta_{nj} + (\omega_{kn} + \xi\dot{\gamma}_{kn})\delta_{mi}\delta_{kj} - 2\xi\dot{\gamma}_{kl}\frac{\partial a_{mnkl}}{\partial a_{ij}} \tag{26}$$

$$\frac{\partial \dot{a}_{mn}^{IRD}}{\partial a_{ij}} = -2D_r\alpha\delta_{mi}\delta_{nj} \tag{27}$$

We obtain expressions for derivatives of $a_{mnkl}$ (i.e., $\partial a_{mnkl}/\partial a_{ij}$) in Equation (26) and others using well known closure approximations in sections to follow.

2. Strain Reduction Factor (SRF) Model

The SRF model was developed by Huynh [8] to reduce the rate of alignment as compared to the FT model by introducing a strain reduction factor $\kappa$ $(0 < \kappa < 1)$ into Equation 16 (i.e., multiplying the right-hand-side of Equation 16 by $\kappa$) to slow down the orientation kinetics as observed experimentally. He based his premise on a reduced bulk strain of fiber clusters in a concentrated suspension flow. It was shown that predictions of steady state orientation in simple shear flow using a suitable value of $\kappa$ matched experimental results [9]. Unfortunately, the reduced strain model was shown to introduce an initial overshoot at small strain not observed in the test data.

The residual and Jacobian in this case is just a multiplication of $\kappa$ with that previously obtained for the FT model.

$$R_{mn}^{SRF} = \kappa R_{mn}^{FT}, \qquad J_{mnij}^{FT} = \kappa J_{mnij}^{FT}, \qquad 0 < \kappa < 1 \tag{28}$$

The SRF model does not satisfy the rheological test of material objectivity and results are dependent on the coordinate system and cannot be applied to complex flows.

3. Reduced Strain Closure (RSC) Model

To address the material objectivity drawback of the SRF model, Wang et al. [9] developed a reduced strain closure (RSC) model which applied a reduction factor only to the evolution rate of the spectral decomposed principal directions of the orientation tensor $\underline{\dot{\lambda}}$, without modifying the rate of rotation $\underline{\underline{\dot{\Phi}}}$ as

$$\dot{\lambda}_i^{RSC} = \kappa\dot{\lambda}_i^{FT}, \qquad \dot{\Phi}_{ij}^{RSC} = \dot{\Phi}_{ij}^{FT}, \qquad a_{mn}|a_{mn} = \lambda_i\Phi_{mi}\Phi_{ni} \tag{29}$$

Based on this model, the modified material derivative is [9]

$$\dot{a}_{mn}^{RSC} = \dot{a}_{mn}^{FT} - (1-\kappa)\dot{a}_{mn}^{\Delta FT} \tag{30}$$

where the kernel of the last term, $\dot{a}_{mn}^{\Delta FT}$ defining the part of that standard FT model affected by the modification factor, $\kappa$ to yield the RSC model is given as

$$\dot{a}_{mn}^{\Delta FT} = \dot{\lambda}_i^{FT}\Phi_{mi}\Phi_{ni} = 2\xi\dot{\gamma}_{kl}(L_{mnkl} - M_{mnrs}a_{rskl}) + \dot{a}_{mn}^{IRD} \tag{31}$$

and

$$L_{mnkl} = \dot{\lambda}_i\Phi_{mi}\Phi_{ni}\Phi_{ki}\Phi_{li}, \qquad M_{mnkl} = \Phi_{mi}\Phi_{ni}\Phi_{ki}\Phi_{li} \tag{32}$$

The Newton-Raphson residual for the RSC model is $R_{mn}^{RSC} = \dot{a}_{mn}^{RSC}$ and the Jacobian is obtained by taking partial derivatives as

$$J_{mnij}^{RSC} = J_{mnij}^{FT} - (1-\kappa)\frac{\partial \dot{a}_{mn}^{\Delta FT}}{\partial a_{ij}} \tag{33}$$

where,

$$\frac{\partial \dot{a}_{mn}^{\Delta FT}}{\partial a_{ij}} = 2\xi\dot{\gamma}_{kl}\frac{\partial}{\partial a_{ij}}\{L_{mnkl} - M_{mnrs}a_{rskl}\} + \frac{\partial \dot{a}_{mn}^{IRD}}{\partial a_{ij}} \tag{34}$$

Expanding Equation (34) and applying the distributive property of differentiation we obtain

$$\frac{\partial \dot{a}_{mn}^{\Delta FT}}{\partial a_{ij}} = \xi \dot{\gamma}_{kl} \left[ \frac{\partial L_{mnkl}}{\partial a_{ij}} - a_{rskl} \frac{\partial M_{mnrs}}{\partial a_{ij}} - M_{mnrs} \frac{\partial a_{rskl}}{\partial a_{ij}} \right] + \frac{\partial \dot{a}_{mn}^{IRD}}{\partial a_{ij}} \tag{35}$$

Finally, applying the product rule of differentiation, we obtain the derivatives of 4th order tensors $M_{mnkl}$ and $L_{mnkl}$, respectively, as

$$\begin{aligned} \frac{\partial M_{mnkl}}{\partial a_{rs}} &= \frac{\partial}{\partial a_{rs}} \{\Phi_{mi}\Phi_{ni}\Phi_{ki}\Phi_{li}\} \\ &= \Phi_{ni}\Phi_{ki}\Phi_{li} \frac{\partial \Phi_{mi}}{\partial a_{rs}} + \Phi_{mi}\Phi_{ki}\Phi_{li} \frac{\partial \Phi_{ni}}{\partial a_{rs}} + \Phi_{mi}\Phi_{ni}\Phi_{li} \frac{\partial \Phi_{ki}}{\partial a_{rs}} \\ &\quad + \Phi_{mi}\Phi_{ni}\Phi_{ki} \frac{\partial \Phi_{li}}{\partial a_{rs}} \end{aligned} \tag{36}$$

and

$$\frac{\partial L_{mnkl}}{\partial a_{rs}} = \Phi_{mi}\Phi_{ni}\Phi_{ki}\Phi_{li} \frac{\partial \dot{\lambda}_i}{\partial a_{rs}} + \dot{\lambda}_i \frac{\partial}{\partial a_{rs}} \{\Phi_{mi}\Phi_{ni}\Phi_{ki}\Phi_{li}\} \tag{37}$$

where the procedure for obtaining the derivatives of the eigenvalues and eigenvectors appear in Appendix A.

## 4. Retarding Principal Rate (RPR) Model

Tseng et al. [10], [11]developed a Retarding Principal Rate (RPR) model, which, like the RSC model, reduce fiber orientation kinetics. Their approach is based on a coaxial modification to the principal directions of the orientation tensor evolution rate via a nonlinear correlation. The material derivative tensor of any standard model X, (i.e., $\dot{a}_{mn}^{X}$) is linearly combined with the RPR correction to slow down the response rate. i.e.

$$\dot{a}_{mn}^{X-RPR} = \dot{a}_{mn}^{X} + \dot{a}_{mn}^{RPR} \tag{38}$$

where the RPR correction $\dot{a}_{mn}^{RPR}$ is written in terms of the eigenvalue matrix $\dot{\Lambda}_{kl}^{IOK}$ as

$$\dot{a}_{mn}^{RPR} = -\Phi_{mk}\dot{\Lambda}_{kl}^{IOK}\Phi_{nl}, \qquad \dot{\Lambda}_{kl}^{IOK} = \dot{\Lambda}_{kl}^{IOK}\left(\underline{\dot{\lambda}}^{X}\right) \tag{39}$$

The superscript on $\dot{\Lambda}_{kl}^{IOK}$ indicates the intrinsic orientation kinetics (IOK) assumption. Given that the correction is coaxial, the rotation tensor growth rate is unaffected and is obtained from the spectral decomposition of $a_{mn}^{X}$. i.e.,

$$\underline{\underline{\Phi}} \mid \Lambda_{mn}^{X} = \Phi_{km} a_{kl}^{X} \Phi_{ln} \tag{40}$$

where | denotes the mathematical abbreviation for 'such that'. The columns of the eigenmatrix obtained in this way are reordered in descending order with respect to the magnitude of the eigenvalues i.e.,

$$\Phi_{ij} = \left\{ \Phi_{ij} \mid \lambda_j^{X} = \Lambda_{jj}^{X}, \qquad \lambda_1^{X} \geq \lambda_2^{X} \geq \lambda_3^{X} \right\} \tag{41}$$

The growth rate of the principal eigenvalues of the standard model – X, i.e. $\dot{\Lambda}_{kl}^{X}$ is obtained from

$$\dot{\Lambda}_{kl}^{X} = \Phi_{km}\dot{a}_{kl}^{X}\Phi_{ln}, \qquad \dot{\lambda}_k^{X} = \dot{\Lambda}_{kk}^{X} \tag{42}$$

The correction to the principal values of the orientation tensor material derivative based on the IOK assumption $\dot{\Lambda}_{kk}^{IOK}$ is defined by a 2-parameter non-linear correlation to the principal values of the standard model $\dot{\Lambda}_{kl}^{X}$ such that

$$\dot{\Lambda}_{kk}^{IOK} = \dot{\lambda}_k^{IOK} = \alpha \left[ \dot{\lambda}_k^{X} - \beta \left( \left\{ \dot{\lambda}_k^{X} \right\}^2 + 2\dot{\lambda}_l^{X}\dot{\lambda}_m^{X} \right) \right], \qquad \dot{\Lambda}_{kl}^{IOK}\Big|_{k\neq l} = 0 \tag{43}$$

For an RPR corrected model, the NR residual $R_{mn}^{X-RPR}$ is simply the material derivative

$$R_{mn}^{X-RPR} = \dot{a}_{mn}^{X-RPR} \tag{44}$$

and the Jacobian $J_{mnij}^{X-RPR}$ is given as

$$J_{mnij}^{X-RPR} = \frac{\partial \dot{a}_{mn}^{X}}{\partial a_{ij}} + \frac{\partial \dot{a}_{mn}^{RPR}}{\partial a_{ij}} \tag{45}$$

The partial derivative of the RPR correction term $\dot{a}_{mn}^{RPR}$ is

$$\frac{\partial \dot{\mathrm{a}}_{mn}^{RPR}}{\partial \mathrm{a}_{ij}} = -\left\{\frac{\partial \Phi_{mk}}{\partial \mathrm{a}_{ij}} \dot{\Lambda}_{kl}^{IOK} \Phi_{nl} + \Phi_{mk} \frac{\partial \dot{\Lambda}_{kl}^{IOK}}{\partial \mathrm{a}_{ij}} \Phi_{nl} + \Phi_{mk} \dot{\Lambda}_{kl}^{IOK} \frac{\partial \Phi_{nl}}{\partial \mathrm{a}_{ij}}\right\} \tag{46}$$

and the partial derivative of the modified growth rate of eigenvalues tensor $\dot{\Lambda}_{kl}^{IOK}$ is obtained by taking the derivative of Equation (43) as

$$\frac{\partial \dot{\Lambda}_{kk}^{IOK}}{\partial \mathrm{a}_{ij}} = \frac{\partial \dot{\lambda}_{k}^{IOK}}{\partial \mathrm{a}_{ij}} = \alpha \left[\frac{\partial \dot{\lambda}_{k}^{X}}{\partial \mathrm{a}_{ij}} - 2\beta \left(\dot{\lambda}_{k}^{X} \frac{\partial \dot{\lambda}_{k}^{X}}{\partial \mathrm{a}_{ij}} + \frac{\partial \dot{\lambda}_{l}^{X}}{\partial \mathrm{a}_{ij}} \dot{\lambda}_{m}^{X} + \dot{\lambda}_{l}^{X} \frac{\partial \dot{\lambda}_{m}^{X}}{\partial \mathrm{a}_{ij}}\right)\right], \qquad \left.\frac{\partial \dot{\Lambda}_{kl}^{IOK}}{\partial \mathrm{a}_{ij}}\right|_{k \neq l} = 0 \tag{47}$$

where

$$\frac{\partial \dot{\Lambda}_{kl}^{X}}{\partial \mathrm{a}_{ij}} = \left\{\frac{\partial \Phi_{km}}{\partial \mathrm{a}_{ij}} \dot{\mathrm{a}}_{kl}^{X} \Phi_{ln} + \Phi_{km} \frac{\partial \dot{\mathrm{a}}_{kl}^{X}}{\partial \mathrm{a}_{ij}} \Phi_{ln} + \Phi_{km} \dot{\mathrm{a}}_{kl}^{X} \frac{\partial \Phi_{ln}}{\partial \mathrm{a}_{ij}}\right\}, \qquad \frac{\partial \dot{\lambda}_{k}^{X}}{\partial \mathrm{a}_{ij}} = \frac{\partial \dot{\Lambda}_{kk}^{X}}{\partial \mathrm{a}_{ij}} \tag{48}$$

In the above, $\partial \dot{\mathrm{a}}_{mn}^{X} / \partial \mathrm{a}_{ij}$ is the partial derivative of material derivative $\dot{\mathrm{a}}_{mn}^{X}$ with respect to the components of $\mathrm{a}_{ij}$ obtained a priori and the partial derivatives of the eigenmatrix with respect to the same $\left(i.e., \partial \Phi_{mn} / \partial \mathrm{a}_{ij}\right)$ can be obtained through any method in Appendix A (cf. [43]).

## 5. Anisotropic Rotary Diffusion (ARD) Models

While the IRD models were experimental observed to be accurate in predicting the orientation state of short-fiber/thermoplastic composites (SFT) they were ineffective in accurately predicting the complete set of orientation tensor components for the long-fiber/thermoplastic composites (LFT) which was the motivation for ARD model development. Various ARD models with different modifications have been developed based on the definition of the spatial diffusion tensor. Most models utilize the moment-tensor form for the ARD representation developed by Phelps and Tucker [12]. The general expression for the 2nd order orientation tensor evolution rate is a linear combination of the Jeffery's model and the and the rotary diffusion term given as

$$\dot{\mathrm{a}}_{mn}^{PT} = \dot{\mathrm{a}}_{mn}^{HD} + \dot{\mathrm{a}}_{mn}^{ARD} \tag{49}$$

where the rotary diffusion term $\dot{\mathrm{a}}_{mn}^{ARD}$ is defined in terms of the spatial diffusion coefficient and the orientation state and is given as

$$\dot{\mathrm{a}}_{mn}^{ARD} = \dot{\gamma}[2C_{mn} - 2C_{rs}\delta_{rs}\mathrm{a}_{mn} - 5(C_{mk}\mathrm{a}_{kn} + \mathrm{a}_{mk}C_{kn}) + 10\mathrm{a}_{mnkl}C_{kl}] \tag{50}$$

and $\underline{\underline{C}}$ is the spatial diffusion tensor. Based on this model, the NR residual and Jacobian are respectively given as

$$R_{mn}^{PT} = \dot{\mathrm{a}}_{mn}^{PT} \tag{51}$$

$$J_{mnij}^{PT} = \frac{\partial \dot{\mathrm{a}}_{mn}^{PT}}{\partial \mathrm{a}_{ij}} = \frac{\partial \dot{\mathrm{a}}_{mn}^{HD}}{\partial \mathrm{a}_{ij}} + \frac{\partial \dot{\mathrm{a}}_{mn}^{ARD}}{\partial \mathrm{a}_{ij}} \tag{52}$$

where the derivative of the rotary diffusion ($ARD$) term is obtained using product rule as

$$\begin{aligned} \frac{\partial \dot{\mathrm{a}}_{mn}^{ARD}}{\partial \mathrm{a}_{ij}} = \dot{\gamma} \Bigg[ & 2\frac{\partial C_{mn}}{\partial \mathrm{a}_{ij}} - 2\left(\frac{\partial C_{rs}}{\partial \mathrm{a}_{ij}} \delta_{rs}\mathrm{a}_{mn} + C_{rs}\delta_{rs}\delta_{mi}\delta_{nj}\right) + \cdots \\ & - 5\left(\frac{\partial C_{mk}}{\partial \mathrm{a}_{ij}} \mathrm{a}_{kn} + C_{mk}\delta_{ki}\delta_{nj} + \delta_{mi}\delta_{kj}C_{kn} + \mathrm{a}_{mk}\frac{\partial C_{kn}}{\partial \mathrm{a}_{ij}}\right) + \cdots \\ & + 10\left(\frac{\partial \mathrm{a}_{mnkl}}{\partial \mathrm{a}_{ij}} C_{kl} + \mathrm{a}_{mnkl}\frac{\partial C_{kl}}{\partial \mathrm{a}_{ij}}\right)\Bigg] \end{aligned} \tag{53}$$

In the mold-flow ARD (mARD) model developed by Bakharev [20], the Phelps & Tucker's rotary diffusion ($ARD$) expression is truncated to include just the linear terms. i.e.

$$\dot{\mathrm{a}}_{mn}^{mARD} = \dot{\gamma}[2C_{mn} - 2C_{kl}\delta_{kl}\mathrm{a}_{mn}] \tag{54}$$

$$\frac{\partial \dot{a}_{mn}^{mARD}}{\partial a_{ij}} = \dot{\gamma}\left[2\frac{\partial C_{mn}}{\partial a_{ij}} - 2\left(\frac{\partial C_{kl}}{\partial a_{ij}}\delta_{kl}a_{mn} + C_{kl}\delta_{kl}\delta_{mi}\delta_{nj}\right)\right] \tag{55}$$

The corresponding evolution rate equation for the 2nd order orientation tensor based on mARD model is given as

$$\dot{a}_{mn}^{mPT} = \dot{a}_{mn}^{HD} + \dot{a}_{mn}^{mARD} \tag{56}$$

Various models for the spatial diffusion coefficient $C_{mn}$ used in the ARD model have been developed by various researchers. The basic representation of $C_{mn}$ by Phelps and Tucker [280] based on a modification of Hand's anisotropic tensor [17] is given as a function of the rate of deformation tensor and orientation state as

$$C_{mn}^{PT} = b_1\delta_{mn} + b_2 a_{mn} + b_3 a_{mk}a_{nk} + \frac{b_4}{\dot{\gamma}}\dot{\gamma}_{mn} + \frac{b_5}{\dot{\gamma}^2}\dot{\gamma}_{mk}\dot{\gamma}_{nk} \tag{57}$$

where $b_i$ are dimensionless constants obtained from regression analysis of experimental data. For this model, the derivative of the $C_{mn}^{PT}$ with respect to $a_{ij}$ is given as

$$\frac{\partial C_{mn}^{PT}}{\partial a_{ij}} = b_2\delta_{mi}\delta_{nj} + b_3\left(\delta_{mi}\delta_{kj}a_{nk} + a_{mk}\delta_{ni}\delta_{kj}\right) \tag{58}$$

The sensitivity of the PT model parameters $b_i$ to ensure numerical stability of the model response coupled with the complicated process involved in their determination were the major limitations to this model application. Tseng et al. [18] developed an improved anisotropic rotary diffusion model (iARD) based on a definition of a two-parameter spatial diffusion tensor model in terms of the rate of deformation tensor that couples the effect of fiber-matrix interaction and fiber-fiber interaction given as

$$C_{mn}^{iARD} = C_I\left(\delta_{mn} - 4C_M\frac{\dot{\gamma}_{mk}\dot{\gamma}_{nk}}{\dot{\gamma}^2}\right) \tag{59}$$

where $C_I$ & $C_M$ are the fiber-fiber and fiber-matrix interaction parameters respectively. An alternate definition is given as

$$C_{mn}^{iARD} = C_I\left(\delta_{mn} - C_M\tilde{L}_{mn}\right), \qquad \tilde{L}_{mn} = (L_{mk}L_{nk})/(L_{rs}L_{rs}) \tag{60}$$

The derivative of the spatial diffusion tensor with respect to the 2nd order orientation tensor is simply zero. i.e.,

$$\frac{\partial C_{mn}^{iARD}}{\partial a_{ij}} = 0 \tag{61}$$

Because of the material non-objectivity of the rate of deformation tensor $L_{mn}$ used in the definition of the spatial diffusion tensor $C_{mn}$ in the iARD model, Tseng et al. [19] developed an improved objective principal spatial tensor ARD model (pARD) that is coaxial with the orientation tensor given as

$$C_{mn}^{pARD} = \left\{C_I\Phi_{mk}D_{kl}\Phi_{nl}, \qquad \underline{\underline{\Phi}} \mid a_{mn} = \Phi_{mk}\bar{a}_{kl}\Phi_{nl}\right\} \tag{62}$$

where the tensor $D_{kl}$ contains only diagonal terms and its trace is unity. i.e.

$$D_{kl}\delta_{kl} = D_{kk} = 1, \qquad D_{kl}|_{k\neq l} = 0 \tag{63}$$

The derivative of $C_{mn}^{pARD}$ with respect to the 2nd order orientation tensor is given as

$$\frac{\partial C_{mn}^{pARD}}{\partial a_{ij}} = C_I\left\{\frac{\partial \Phi_{mk}}{\partial a_{ij}}D_{kl}\Phi_{nl} + \Phi_{mk}D_{kl}\frac{\partial \Phi_{nl}}{\partial a_{ij}}\right\} \tag{64}$$

Another ARD model reduction suggested by Wang [44] called the WPT model involved truncating the terms of the PT model to just the first and third term such that,

$$C_{mn}^{WPT} = b_1\delta_{mn} + b_3 a_{mk}a_{nk} \tag{65}$$

Falvoro et al. [23] provided an alternative form of the spatial diffusion tensor where he replaced the coefficients with a weighted superposition of the interaction coefficient, i.e.

$$C_{mn}^{WPT} = C_I\big((1-w)\delta_{mn} + w\text{a}_{mk}\text{a}_{nk}\big) \tag{66}$$

where $w$ is the weighting factor. The derivative of $C_{mn}^{pARD}$ with respect to the 2nd order orientation tensor is given as

$$\frac{\partial C_{mn}^{WPT}}{\partial \text{a}_{ij}} = wC_I(\delta_{mk}\text{a}_{nk} + \text{a}_{mk}\delta_{nk}) \tag{67}$$

Lastly, we consider the $D_z$ ARD model development (cf. Falvoro et al. [23]) by Moldflow for simulating 2.5D flow processes. Their model is defined in terms of the interaction coefficient, $C_I$, a moment of interaction thickness parameter $D_z$, and the unit normal to the mold surface $\hat{n}$. The expression for $C_{mn}$ here is given as

$$C_{mn}^{Dz} = C_I^{Dz}(\delta_{mn} - (1-D_z)\hat{n}_m\hat{n}_n) \tag{68}$$

and

$$\frac{\partial C_{mn}^{Dz}}{\partial \text{a}_{ij}} = 0 \tag{69}$$

6. Nematic Potential (NEM) Model

Latz et al. [21] proposed a 2-parameter nematic potential ARD (NEM) model for the diffusion term that couples the phenomenological effect of the momentum diffusion due to fiber-fiber interaction and a topological interaction effect of diffusion due to an exclusion volume mechanism. i.e.

$$\dot{\text{a}}_{mn}^{IRD-MS} = \dot{\gamma}[C_I(\delta_{mn} - \alpha\text{a}_{mn}) + U_0(\text{a}_{mk}\text{a}_{kn} - \text{a}_{kl}\text{a}_{mnkl})] \tag{70}$$

where $U_0$ is the 'Onsager' nematic topological interaction coefficient of the Maier-Saupe potential. Typically, for stability, $U_0 \leq 4C_I$ for 2D analysis and $U_0 > 8C_I$ for 3D analysis. The material derivative of the 2nd order orientation tensor based on the nematic diffusion model is thus given as

$$\dot{\text{a}}_{mn}^{nem} = \dot{\text{a}}_{mn}^{HD} + \dot{\text{a}}_{mn}^{IRD-MS} \tag{71}$$

The NR residual and Jacobian are respectively given as

$$R_{mn}^{nem} = \dot{\text{a}}_{mn}^{nem} \tag{72}$$

$$J_{mnij}^{nem} = \frac{\partial \dot{\text{a}}_{mn}^{HD}}{\partial \text{a}_{ij}} + \frac{\partial \dot{\text{a}}_{mn}^{IRD-MS}}{\partial \text{a}_{ij}} \tag{73}$$

where the derivative of the nematic diffusion term is given as

$$\frac{\partial \dot{\text{a}}_{mn}^{IRD-MS}}{\partial \text{a}_{ij}} = \dot{\gamma}\left[-C_I\alpha\delta_{mi}\delta_{nj} + U_0\left(\delta_{mi}\delta_{kj}\text{a}_{kn} + \text{a}_{mk}\delta_{ki}\delta_{nj} - \delta_{ki}\delta_{lj}\text{a}_{mnkl} - \text{a}_{kl}\frac{\partial \text{a}_{mnkl}}{\partial \text{a}_{ij}}\right)\right] \tag{74}$$

Most commercial software used in simulation of the injection molding process such as Autodesk Moldflow and Moldex3D usually combines multiple models in predicting the orientation state for improved accuracy. One such combination is the ARD-RSC models whose material derivative is expressed as

$$\dot{\text{a}}_{mn}^{pARD-RSC} = \dot{\text{a}}_{mn}^{RSC} - \kappa\dot{\text{a}}_{mn}^{IRD} + \dot{\text{a}}_{mn}^{ARD} + \dot{\text{a}}_{mn}^{\Delta RSC} \tag{75}$$

where,

$$\dot{\text{a}}_{mn}^{\Delta RSC} = -2\dot{\gamma}(1-\kappa)[M_{mnkl} - \delta_{kl}\text{a}_{mn} - 5(L_{mnkl} - M_{mnrs}\text{a}_{rskl})]C_{kl} \tag{76}$$

and $\dot{\text{a}}_{mn}^{RSC}$, $\dot{\text{a}}_{mn}^{IRD}$, $\dot{\text{a}}_{mn}^{ARD}$ have been defined in preceding sections. In this case the NR residual $R_{mn}^{pARD-RSC}$ is the material derivative $\dot{\text{a}}_{mn}^{pARD-RSC}$, i.e.,

$$R_{mn}^{pARD-RSC} = \dot{\text{a}}_{mn}^{pARD-RSC} \tag{77}$$

and the Jacobian is obtained by taking partial derivatives with respect to the 2nd order tensor as usual and can be expressed as

$$J_{mnij}^{pARD-RSC} = \frac{\partial \dot{\text{a}}_{mn}^{RSC}}{\partial \text{a}_{ij}} - \kappa\frac{\partial \dot{\text{a}}_{mn}^{IRD}}{\partial \text{a}_{ij}} + \frac{\partial \dot{\text{a}}_{mn}^{ARD}}{\partial \text{a}_{ij}} + \frac{\partial \dot{\text{a}}_{mn}^{\Delta RSC}}{\partial \text{a}_{ij}} \tag{78}$$

where

$$\frac{\partial \dot{a}_{mn}^{\Delta RSC}}{\partial a_{ij}} = -2\dot{\gamma}(1-\kappa)\left\{\begin{matrix}\left[\frac{\partial M_{mnkl}}{\partial a_{ij}} - \delta_{kl}\delta_{mi}\delta_{nj} - 5\frac{\partial}{\partial a_{ij}}\{L_{mnkl} - M_{mnrs}a_{rskl}\}\right]C_{kl} + \cdots \\ +[M_{mnkl} - \delta_{kl}a_{mn} - 5(L_{mnkl} - M_{mnrs}a_{rskl})]\frac{\partial C_{kl}}{\partial a_{ij}}\end{matrix}\right\} \quad (79)$$

The terms of the partial derivatives have been previously derived in the preceding section.

*Closure Approximations and Their Explicit Derivatives*

Derivatives of the orientation tensor closure approximation are used in the Newton-Raphson iteration method to compute the steady-state fiber orientation tensor state. These derivatives for the various closure approximation used in this study appear below.

1. Quadratic Closure Approximation

The quadratic (QDR) closure was introduced by Doi [297] and Lipscomb [298] and defined as dyadic product of the 2nd order orientation tensor, $a_{ij}$. We denote the quadratic closure approximate $\tilde{a}_{ijkl}$ and is mathematically given as

$$\tilde{a}_{ijkl} = a_{ij}a_{kl} \quad (80)$$

The derivative of $\tilde{a}_{ijkl}$ above with respect to the 2nd order tensor $a_{mn}$ is simply.

$$\frac{\partial \tilde{a}_{ijkl}}{\partial a_{mn}} = \frac{\partial a_{ij}}{\partial a_{mn}}a_{kl} + a_{ij}\frac{\partial a_{kl}}{\partial a_{mn}} = \delta_{im}\delta_{jn}a_{kl} + a_{ij}\delta_{km}\delta_{ln} \quad (81)$$

The QDR closure inherently lacks symmetry but preserves the required symmetry of the computed lower order tensor.

2. Linear Closure Approximation

The linear closure approximation of $\hat{a}_{ijkl}$ first proposed by Hand [284] using all the products of the 2nd order orientation tensor $a_{ij}$ and the identity matrix $\delta_{ij}$ is given as

$$\begin{aligned}\hat{a}_{ijkl} = &-h_1(\delta_{ij}\delta_{kl} + \delta_{ik}\delta_{jl} + \delta_{il}\delta_{jk}) \\ &+ h_2(a_{ij}\delta_{kl} + a_{ik}\delta_{jl} + a_{il}\delta_{jk} + \delta_{ij}a_{kl} + \delta_{ik}a_{jl} \\ &+ \delta_{il}a_{jk})\end{aligned} \quad (82)$$

The derivative of $\hat{a}_{ijkl}$ above with respect to components of the 2nd order tensor $a_{mn}$ is given as

$$\begin{aligned}\frac{\partial \hat{a}_{ijkl}}{\partial a_{mn}} = h_2(&\delta_{im}\delta_{jn}\delta_{kl} + \delta_{im}\delta_{kn}\delta_{jl} + \delta_{im}\delta_{ln}\delta_{jk} + \delta_{im}\delta_{jn}\delta_{kl} \\ &+ \delta_{im}\delta_{kn}\delta_{jl} + \delta_{im}\delta_{ln}\delta_{jk})\end{aligned} \quad (83)$$

where $h_1$ and $h_2$ are numerical factors which vary based on spatial dimensionality and given in **Table 1** below

**Table 1.** Numerical factors of the linear closure

| | $\boldsymbol{Solid\ (3D)}$ | $\boldsymbol{Planar\ (2D)}$ |
|---|---|---|
| $h_1$ | $1/35$ | $1/24$ |
| $h_2$ | $1/7$ | $1/6$ |

The LIN closure is exact for random orientation distribution while the QDR closures are exact for uniaxially aligned fiber orientation.

3. Hybrid Closure Approximation

The hybrid closure approximation, $a_{ijkl}$ is simply a weighted combination of both linear (LIN) $\hat{a}_{ijkl}$ and quadratic (QDR) $\tilde{a}_{ijkl}$ closure approximation above by some scalar measure of orientation tensor given as [7]

$$\mathrm{a}_{ijkl} = f\tilde{\mathrm{a}}_{ijkl} + (1-f)\hat{\mathrm{a}}_{ijkl} \tag{84}$$

where $f$ is a generalization of Herman's Orientation factor. Advani & Tucker [7] proposed an appropriate approximation of the weighting factor as an invariant of the orientation state given as $f = a_f \mathrm{a}_{ij}\mathrm{a}_{ji} - b_f$, where $a_f$ and $b_f$ are constants that depends on the spatial dimension given in Table 2 below

**Table 2.** Constants of the hybrid closure

| | ***Solid*** $(\mathbf{3D})$ | ***Planar*** $(\mathbf{2D})$ |
|---|---|---|
| $a_f$ | $3/2$ | $2$ |
| $b_f$ | $1/2$ | $1$ |

The derivative of the hybrid closure approximation $\mathrm{a}_4$ with respect to components of the 2nd order tensor $\mathrm{a}_{mn}$ is given as

$$\frac{\partial \mathrm{a}_{ijkl}}{\partial \mathrm{a}_{mn}} = f\frac{\partial \tilde{\mathrm{a}}_{ijkl}}{\partial \mathrm{a}_{mn}} + (1-f)\frac{\partial \hat{\mathrm{a}}_{ijkl}}{\partial \mathrm{a}_{mn}} + \frac{\partial f}{\partial \mathrm{a}_{mn}}\left(\tilde{\mathrm{a}}_{ijkl} - \hat{\mathrm{a}}_{ijkl}\right) \tag{85}$$

where,

$$\frac{\partial f}{\partial \mathrm{a}_{mn}} = a_f\left(\delta_{im}\delta_{jn}\mathrm{a}_{ji} + \mathrm{a}_{ij}\delta_{jm}\delta_{in}\right) \tag{86}$$

An alternative estimation of the factor $f$ by Advani & Tucker [7] is given as

$$f = 1 - \alpha^{\alpha} e_{ijk}\mathrm{a}_{i1}\mathrm{a}_{j2}\mathrm{a}_{k3} \tag{87}$$

$$\frac{\partial f}{\partial \mathrm{a}_{mn}} = -\alpha^{\alpha} e_{ijk}\left\{\delta_{im}\delta_{1n}\mathrm{a}_{j2}\mathrm{a}_{k3} + \mathrm{a}_{i1}\delta_{jm}\delta_{2n}\mathrm{a}_{k3} + \mathrm{a}_{i1}\mathrm{a}_{j2}\delta_{km}\delta_{3n}\right\}$$

The hybrid model is observed to perform better for transient state orientation prediction; however, the hybrid closure was shown to over-predict the steady state fiber alignment compared with the more accurate orientation distribution function closure (ODFC) predictions. ODFCs are, however, computationally expensive since they require a finite difference grid in space and time.

4. Hinch and Leal Closure Approximation

Hinch and Leal [33] developed numerous composite closure approximations for the 4th order tensor in precontracted forms with the deformation rate tensor. The accuracy of their predictions was dependent on the flow field considered and dependent on flow magnitude. The Hinch and Leal closure approximations were not explicit expressions of the 4th order orientation tensor $\mathrm{a}_{ijkl}$ but were in contracted form with the deformation rate tensor i.e., $\gamma_{kl}\mathrm{a}_{ijkl}$. Advani and Tucker developed a general explicit expression of $\mathrm{a}_{ijkl}$ (cf. Equation (88)) summarizing all the Hinch and Leal closures forms given as

$$\begin{aligned}\mathrm{a}_{ijkl} = {} & \beta_1\left(\delta_{ij}\delta_{kl}\right) + \beta_2\left(\delta_{ik}\delta_{jl} + \delta_{il}\delta_{jk}\right) + \beta_3\left(\delta_{ij}\mathrm{a}_{kl} + \mathrm{a}_{ij}\delta_{kl}\right) \\ & + \beta_4\left(\mathrm{a}_{ik}\delta_{jl} + \cdots . + \mathrm{a}_{jl}\delta_{ik} + \mathrm{a}_{il}\delta_{jk} + \mathrm{a}_{jk}\delta_{il}\right) \\ & + \beta_5\left(\mathrm{a}_{ij}\mathrm{a}_{kl}\right) + \beta_6\left(\mathrm{a}_{ik}\mathrm{a}_{jl} + \mathrm{a}_{il}\mathrm{a}_{jk}\right) \\ & + \cdots . + \beta_7\left(\delta_{ij}\mathrm{a}_{km}\mathrm{a}_{ml} + \mathrm{a}_{im}\mathrm{a}_{mj}\delta_{kl}\right) \\ & + \beta_8\left(\mathrm{a}_{im}\mathrm{a}_{mj}\mathrm{a}_{kn}\mathrm{a}_{nl}\right)\end{aligned} \tag{88}$$

and the partial derivative of the above expression with respect to components of the 2nd order orientation tensor $\mathrm{a}_{rs}$ based on product rule is given as

$$
\begin{aligned}
\frac{\partial \mathrm{a}_{ijkl}}{\partial \mathrm{a}_{rs}} = & \Big[\frac{\partial \beta_1}{\partial \mathrm{a}_{rs}}\left(\delta_{ij}\delta_{kl}\right) + \frac{\partial \beta_2}{\partial \mathrm{a}_{rs}}\left(\delta_{ik}\delta_{jl} + \delta_{il}\delta_{jk}\right) + \frac{\partial \beta_3}{\partial \mathrm{a}_{rs}}\left(\delta_{ij}\mathrm{a}_{kl} + \mathrm{a}_{ij}\delta_{kl}\right) \\
& + \frac{\partial \beta_4}{\partial \mathrm{a}_{rs}}\left(\mathrm{a}_{ik}\delta_{jl} + \mathrm{a}_{jl}\delta_{ik} + \mathrm{a}_{il}\delta_{jk} + \mathrm{a}_{jk}\delta_{il}\right) + \frac{\partial \beta_5}{\partial \mathrm{a}_{rs}}\left(\mathrm{a}_{ij}\mathrm{a}_{kl}\right) + \frac{\partial \beta_6}{\partial \mathrm{a}_{rs}}\left(\mathrm{a}_{ik}\mathrm{a}_{jl} + \mathrm{a}_{il}\mathrm{a}_{jk}\right) \\
& + \frac{\partial \beta_7}{\partial \mathrm{a}_{rs}}\left(\delta_{ij}\mathrm{a}_{km}\mathrm{a}_{ml} + \mathrm{a}_{im}\mathrm{a}_{mj}\delta_{kl}\right) + \frac{\partial \beta_8}{\partial \mathrm{a}_{rs}}\left(\mathrm{a}_{im}\mathrm{a}_{mj}\mathrm{a}_{kn}\mathrm{a}_{nl}\right)\Big] \\
& + \Big[\beta_3\left(\delta_{ij}\delta_{kr}\delta_{ls} + \delta_{ir}\delta_{js}\delta_{kl}\right) + \beta_4\left(\delta_{ir}\delta_{ks}\delta_{jl} + \delta_{jr}\delta_{ls}\delta_{ik} + \delta_{ir}\delta_{ls}\delta_{jk} + \delta_{jr}\delta_{ks}\delta_{il}\right) \\
& + \beta_5\left(\delta_{ir}\delta_{js}\mathrm{a}_{kl} + \mathrm{a}_{ij}\delta_{kr}\delta_{ls}\right) + \beta_6\left(\delta_{ir}\delta_{ks}\mathrm{a}_{jl} + \mathrm{a}_{ik}\delta_{jr}\delta_{ls} + \delta_{ir}\delta_{ls}\mathrm{a}_{jk} + \mathrm{a}_{il}\delta_{jr}\delta_{ks}\right) \\
& + \beta_7\left(\delta_{ij}\delta_{kr}\delta_{ms}\mathrm{a}_{ml} + \delta_{ij}\mathrm{a}_{km}\delta_{mr}\delta_{ls} + \delta_{ir}\delta_{ms}\mathrm{a}_{mj}\delta_{kl} + \mathrm{a}_{im}\delta_{mr}\delta_{js}\delta_{kl}\right) \\
& + \beta_8\left(\delta_{ir}\delta_{ms}\mathrm{a}_{mj}\mathrm{a}_{kn}\mathrm{a}_{nl} + \mathrm{a}_{im}\delta_{mr}\delta_{js}\mathrm{a}_{kn}\mathrm{a}_{nl} + \mathrm{a}_{im}\mathrm{a}_{mj}\delta_{kr}\delta_{ns}\mathrm{a}_{nl} + \mathrm{a}_{im}\mathrm{a}_{mj}\mathrm{a}_{kn}\delta_{nr}\delta_{ls}\right)\Big]
\end{aligned}
\tag{89}
$$

Mullens [36] provided a summary Table (cf. Table 3) for the $\beta_i$ factors of the Hinch and Leal closures subdivided into weak flow (WF - Isotropic, Linear and Quadratic), strong flow (SF), and Hinch and Leal composite flows (HL – HL1 & HL2) closure forms.

**Table 3.** Summary of the Hinch and Leal closure $\beta_i$ factors for the different flow classifications

| | | $\boldsymbol{\beta_1}$ | $\boldsymbol{\beta_2}$ | $\boldsymbol{\beta_3}$ | $\boldsymbol{\beta_4}$ | $\boldsymbol{\beta_5}$ | $\boldsymbol{\beta_6}$ | $\boldsymbol{\beta_7}$ | $\boldsymbol{\beta_8}$ |
|---|---|---|---|---|---|---|---|---|---|
| WF | ISO | $\frac{1}{15}$ | $\frac{1}{15}$ | ⋯ | ⋯ | ⋯ | ⋯ | ⋯ | ⋯ |
| | LIN | $-\frac{1}{35}$ | $-\frac{1}{35}$ | $\frac{1}{7}$ | $\frac{1}{7}$ | ⋯ | ⋯ | ⋯ | ⋯ |
| | QDR | ⋯ | ⋯ | ⋯ | ⋯ | 1 | ⋯ | ⋯ | ⋯ |
| SF | SF2 | ⋯ | ⋯ | ⋯ | ⋯ | 1 | 1 | ⋯ | $-\frac{2}{\langle \mathrm{a}^2\rangle}$ |
| HL | HL1 | ⋯ | ⋯ | $\frac{2}{5}$ | ⋯ | $-\frac{1}{5}$ | $\frac{3}{5}$ | $-\frac{2}{5}$ | ⋯ |
| | HL2 | $\frac{26}{315}\alpha$ | $\frac{26}{315}\alpha$ | $\frac{16}{63}\alpha$ | $-\frac{4}{21}\alpha$ | 1 | 1 | ⋯ | $-\frac{2}{\langle \mathrm{a}^2\rangle}$ |

where the parameters $\langle \mathrm{a}^2\rangle$ and $\alpha$ are respectively

$$
\langle \mathrm{a}^2\rangle = \mathrm{a}_{ij}\mathrm{a}_{ji}, \qquad \alpha = \exp\left[2\frac{1-3\langle \mathrm{a}^2\rangle}{1-\langle \mathrm{a}^2\rangle}\right] \tag{90}
$$

and the partial derivatives are respectively given as

$$
\begin{gathered}
\frac{\partial \langle \mathrm{a}^2\rangle}{\partial \mathrm{a}_{rs}} = \delta_{ir}\delta_{js}\mathrm{a}_{ji} + \mathrm{a}_{ij}\delta_{jr}\delta_{is} \\
\frac{\partial \alpha}{\partial \mathrm{a}_{rs}} = -\frac{4\alpha}{(1-\langle \mathrm{a}^2\rangle)^2}\frac{\partial \langle \mathrm{a}^2\rangle}{\partial \mathrm{a}_{rs}}, \quad \frac{\partial}{\partial \mathrm{a}_{rs}}\left\{\frac{k}{\langle \mathrm{a}^2\rangle}\right\} = -\frac{k}{\langle \mathrm{a}^2\rangle^2}\frac{\partial \langle \mathrm{a}^2\rangle}{\partial \mathrm{a}_{rs}}
\end{gathered}
\tag{91}
$$

## 5. Eigenvalue based Fitted (EBF) Closure Approximations

Recently, more accurately fitted closure approximations have been developed including the Eigenvalue-Based Fitted (EBF) closures that require a principal axis transformation and the Invariant-Based Fitted (IBF) closures. The idea of orthotropic closure approximations for the 4th order tensor was to impose objectivity such that the approximation is independent of the coordinate frame selection. The Orthotropic Smooth (ORS), Orthotropic Fitted (ORF) and ORF low fiber-fiber interaction coefficient (ORL) closures are included in the EBF class of closures originally developed by Cintra and Tucker [35] based on the assumption of coincident principal axis of the 2nd and 4th order tensor. The (9x9) term 4th order tensor can be represented in (6 x 6) contracted notation like in structural analysis of composite material based on symmetry property. i.e.

$$A_{rs} = \mathrm{a}_{ijkl} \tag{92}$$

where the index of the contracted notation is related to the index notation according to

$$r = \begin{cases} i = j & \delta_{ij} = 1 \\ (9 - i - j) & \delta_{ij} = 0 \end{cases} \;\&\; s = \begin{cases} k = l & \delta_{kl} = 1 \\ (9 - k - l) & \delta_{kl} = 0 \end{cases} \tag{93}$$

The derivative of the 4th order tensor with respect to the 2nd order tensor is such that

$$\frac{\partial A_{rs}}{\partial \mathrm{a}_{mn}} = \frac{\partial \mathrm{a}_{ijkl}}{\partial \mathrm{a}_{mn}} \tag{94}$$

Symmetry property of the 4th order tensor requires $\mathrm{a}_{ijkl} = \mathrm{a}_{klij}$ which implies that $A_{rs} = A_{sr}$ . The contracted tensor $A_{rs}$ transformed to the principal axes has the orthotropic form $\bar{A}_{rs}$ given thus.

$$\underline{\underline{\bar{A}}} = \begin{bmatrix} \bar{A}_{11} & \bar{A}_{12} & \bar{A}_{13} & & & \\ \bar{A}_{21} & \bar{A}_{22} & \bar{A}_{23} & & & \\ \bar{A}_{31} & \bar{A}_{32} & \bar{A}_{33} & & & \\ & & & \bar{A}_{44} & & \\ & & & & \bar{A}_{55} & \\ & & & & & \bar{A}_{66} \end{bmatrix} \tag{95}$$

The contracted tensor transforms from its principal reference frame to the original coordinate axes according to

$$A_{rs} = \mathcal{M}_{ri}\mathcal{M}_{sj}\bar{A}_{ij} \tag{96}$$

The 6x6 transformation matrix $\mathcal{M}_{ij}$ is given as $\mathcal{M}_{ij} = F_{im}\, Q_{mn}\, F_{nj}^{-1}$, where $F_{ij} = k\delta_{ij},\ k = \begin{cases} 1 & i \leq 3 \\ 2 & i > 3 \end{cases}$ and $Q_{rs} = \Phi_{ik}\Phi_{jl} + (1 - \delta_{kl})\Phi_{jk}\Phi_{il}$. The modal matrix $\Phi_{ij}$ whose kth column are the corresponding eigenvectors $\underline{x}^k$ of eigenvalues $\lambda_k = \Lambda_{kk}$ is obtained from the spectral decomposition of $\mathrm{a}_{ij}$ is such that:

$$\underline{\underline{\Phi}} \mid \mathrm{a}_{mn} = \Phi_{mk}\Lambda_{kl}\Phi_{nl} \tag{97}$$

The indices of the contracted 4th order modal tensor $Q_{rs}$ relates to the those of the 2nd order modal matrix $\underline{\underline{\Phi}}$ according to the above equation. A more direct way is to reconstruct the 4th order orientation tensor $\bar{\mathrm{a}}_{mnpq}$ from the contracted form $A_{rs}$ and transform from the principal reference frame to the original axes according to Equation (98) below.

$$\mathrm{a}_{ijkl} = \Phi_{im}\Phi_{jn}\Phi_{kp}\Phi_{lq}\bar{\mathrm{a}}_{mnpq} \tag{98}$$

and using the product rule

$$\begin{aligned} \frac{\partial \mathrm{a}_{ijkl}}{\partial \mathrm{a}_{rs}} = q_{im}q_{jn}q_{kp}q_{lq}\frac{\partial \bar{a}_{mnpq}}{\partial \mathrm{a}_{rs}} + \cdots \\ + \left(\frac{\partial q_{im}}{\partial \mathrm{a}_{rs}} q_{jn}q_{kp}q_{lq} + q_{im}\frac{\partial q_{jn}}{\partial \mathrm{a}_{rs}} q_{kp}q_{lq} + q_{im}q_{jn}\frac{\partial q_{kp}}{\partial \mathrm{a}_{rs}} q_{lq} \right. \\ \left. + q_{im}q_{jn}q_{kp}\frac{\partial q_{lq}}{\partial \mathrm{a}_{rs}}\right)\bar{a}_{mnpq} \end{aligned} \tag{99}$$

Derivative of the eigentensor $\underline{\underline{\Phi}}$ can be found in [43] (cf. Appendix A). Symmetry requirements of the transformed orthotropic tensor reduces the total number of independent non-zero components to 9, and additional special symmetry properties of the exact 4th order tensor requires that $\bar{\mathrm{a}}_{ijkl} = \bar{\mathrm{a}}_{kjil} = \bar{\mathrm{a}}_{ljki} = \bar{\mathrm{a}}_{ikjl} = \bar{\mathrm{a}}_{ilkj}$ reduces the non-zero independent components to the 6 diagonal terms. i.e.

$$\bar{A}_{ij} = \bar{A}_{kk} \quad \{k:\ k = 9 - i - j\,, \qquad i \neq j \tag{100}$$

The normalization property $\mathrm{a}_{ijkk} = \mathrm{a}_{ij}$ of the exact 4th order tensor further requires that:

$$\begin{bmatrix} \bar{A}_{44} \\ \bar{A}_{55} \\ \bar{A}_{66} \end{bmatrix} = \underline{\underline{ß}}^{-1} \left\{ \begin{bmatrix} \lambda_1 \\ \lambda_2 \\ \lambda_3 \end{bmatrix} - \begin{bmatrix} \bar{A}_{11} \\ \bar{A}_{22} \\ \bar{A}_{33} \end{bmatrix} \right\} \tag{101}$$

where $\lambda_i$ are the eigenvalues of the 2nd order orientation tensor $\mathrm{a}_{ij}$, $\sum_i \lambda_i = 1$ and $ß_{ij} = 1 - \delta_{ij}$. Based on the foregoing conditions, the only three surviving non-zero independent terms are $\bar{A}_{11}$, $\bar{A}_{22}$ & $\bar{A}_{33}$. The general form for orthotropic closure is to express the three surviving non-zero independent components ($\bar{A}_{11}$, $\bar{A}_{22}$, $\bar{A}_{33}$) of the contracted 4th order tensor in the principal reference frame after imposing all symmetric and normalization conditions of the exact 4th order tensor, as a scalar function $F_k(\lambda_1, \lambda_2)$ of the two largest eigenvalues $(\lambda_1, \lambda_2)$ of the 2nd order tensor. Most fitted closures take the form of an nth-order binomial function in $\lambda_1$ & $\lambda_2$ to represents the scalar function i.e.,

$$\bar{A}_{kk} = F_k(\lambda_1, \lambda_2) = f_k^{(n)}(\lambda_1, \lambda_2), \quad \lambda_1 \geq \lambda_2 \geq \lambda_3, \qquad k = 1,2,3 \tag{102}$$

Polynomial order exceeding $n \geq 4$ fall under the class of eigenvalue based optimal fitting (EBOF) closures. Generally, we can represent the function $f_k^{(n)}$ as a tensor product of a constant coefficient matrix $\mathfrak{C}_{ij}^{(n)}$ and a nth order permuted bivariate polynomial vector $\underline{\mathbb{A}}^{(n)} = \underline{\mathbb{A}}^{(n)}\ (\lambda_1, \lambda_2)$, i.e.

$$f^{(n)}(\lambda_1, \lambda_2) = \mathfrak{C}_{kj}^{(n)}\ \mathbb{A}_j^{(n)} \tag{103}$$

Different representation of $\underline{\underline{\mathfrak{C}}}^{(n)}$ and $\underline{\mathbb{A}}^{(n)}$ depending on the polynomial order fit $(n)$ can be found in Appendix B. The derivative of the components of the orthotropic closure with respect to the 2nd order tensor are thus:

$$\frac{\partial \bar{A}_{kk}}{\partial \mathrm{a}_{rs}} = \frac{\partial \bar{A}_k}{\partial \mathrm{a}_{rs}} = \mathfrak{C}_{kj}^{(n)} \frac{\partial\ \mathbb{A}_j^{(n)}}{\partial \mathrm{a}_{rs}} = \mathfrak{C}_{kj}^{(n)}\ \mathbb{A}'^{(n)}_{jrs} = \mathfrak{C}_{kj}^{(n)}\ \tilde{\mathbb{A}}_{jl}^{(n)}\ \lambda'_{lrs} \qquad k = 1,2,3, \quad l = 1,2 \tag{104}$$

The nth order binomial permutation vector $\mathbb{A}_k^{(n)}$ and its derivative coefficient matrix $\tilde{\mathbb{A}}_{kl}^{(n)}$ for the quadratic closure are given from terms of binomial expansion respectively as

$$\begin{cases} \mathbb{A}_k^{(n)}\ (\lambda_1, \lambda_2) = \lambda_1^{i-j} \lambda_2^j \\ \tilde{\mathbb{A}}_{kl}^{(n)} = \dfrac{\partial\ \mathbb{A}_k^{(n)}}{\partial \lambda_l} = \begin{cases} (i-j) \cdot \lambda_1^{i-j-1} \lambda_2^j & l = 1\,, \quad k | k = j + \dfrac{1}{2} i(i+1) \\ j \cdot \lambda_1^{i-j} \lambda_2^{j-1} & l = 2 \end{cases} \end{cases} \qquad j = 0 \cdots i, \quad i = 0 \cdots n \tag{105}$$

For a special case of orthotropic fitted closure called rational elliptical (RE) closure by Wetzel and Tucker [45], the scalar function for the 3 independent tensor component is given as

$$F\ (\lambda_1, \lambda_2) = \frac{f^{(n)}\ (\lambda_1, \lambda_2)}{f^{(m)}(\lambda_1, \lambda_2)} \tag{106}$$

The derivative of the components of the above with respect to the 2nd order tensor based on the quotient rule is thus:

$$\frac{\partial \bar{A}_{kk}}{\partial \mathrm{a}_{rs}} = \frac{1}{[f^{(m)}]^2} \left[ f^{(m)}\ \frac{\partial f^{(n)}}{\partial \mathrm{a}_{rs}} - f^{(n)}\ \frac{\partial f^{(m)}}{\partial \mathrm{a}_{rs}} \right] \tag{107}$$

From normalization condition of the 4th order tensor, we obtain for the derivative of $\bar{A}_{kk}$, $(k = 4,5,6)$

$$\frac{\partial \bar{A}_{kk}}{\partial \mathrm{a}_{rs}} = ß_{ki}^{-1}\left\{\frac{\partial \lambda_i}{\partial \mathrm{a}_{rs}} - \frac{\partial \bar{A}_{ii}}{\partial \mathrm{a}_{rs}}\right\}, \quad \frac{\partial \lambda_i}{\partial \mathrm{a}_{rs}} = \mathcal{E}_{il}\lambda'_{lrs}, \quad \underline{\underline{\mathcal{E}}} = \begin{bmatrix} 1 & 0 & -1 \\ 0 & 1 & -1 \end{bmatrix}^T \qquad (108)$$

$$i = 1,2,3, \quad l = 1,2$$

For the partial derivative of the eigenvalues with respect to the components of the 2nd order orientation tensor, kindly refer to Appendix A. EBF closures are computationally more involved in numerical calculations of actual flows because of the principal axis transformation.

## 6. Invariant Based Fitted (IBF) Closure Approximations

Of the class of IBF closures, the natural (NAT) closure approximation of Verleye et al. [46] was built on the work of Lipscomb et al. [47] and formed the basis for other IBF developments. They developed a general expression for the full symmetric 4th order tensor in terms of the 2nd order tensor, the identity matrix and fitted coefficients as functions of the tensor invariant which were derived from analytical calculations based on a least square fitting process. The NAT closure assumed the absence of fiber-fiber interaction and infinitely long fiber geometry. The closure is exact based on the foregoing assumptions however it has been reported to possess singularities for axisymmetric orientation states. The Invariant based optimal fitting (IBOF) closure developed by Chung and Kwon [48] was an extension to the NAT closure development however the independent coefficients are derived from regression analysis based on actual data obtained from DFC considering different flow types like EBF closures. In contracted form the 4th order tensor based on symmetry properties is given as

$$\underline{\underline{A}} = \begin{bmatrix} A_{11} & A_{12} & A_{13} & A_{14} & A_{15} & A_{16} \\ & A_{22} & A_{23} & A_{24} & A_{24} & A_{26} \\ & & A_{33} & A_{34} & A_{35} & A_{36} \\ & & & A_{44} & A_{45} & A_{46} \\ & & & & A_{55} & A_{56} \\ \ldots Sym & & & & & A_{66} \end{bmatrix} \qquad (109)$$

based on special symmetry requirement

$$\begin{matrix} A_{44} = A_{23}, & A_{45} = A_{36}, & A_{46} = A_{25} \\ A_{55} = A_{13}, & A_{56} = A_{14}, & A_{66} = A_{12} \end{matrix} \qquad (110)$$

and from the normalization condition

$$\sum_{n=1}^{3} A_{nm} = \mathrm{a}_m, \qquad \mathrm{a}_m = \mathrm{a}_{ij}, \qquad m = \begin{cases} i = j & i = j \\ 9 - i - j & i \neq j \end{cases} \qquad (111)$$

Or more explicitly we derive the sets of equations in Equation (112) below.

$$\begin{matrix} A_{11} + A_{12} + A_{13} = \mathrm{a}_{11}, & A_{12} + A_{22} + A_{23} = \mathrm{a}_{22} \\ A_{13} + A_{23} + A_{33} = \mathrm{a}_{33}, & A_{14} + A_{24} + A_{34} = \mathrm{a}_{23} \\ A_{15} + A_{25} + A_{35} = \mathrm{a}_{13}, & A_{16} + A_{26} + A_{36} = \mathrm{a}_{12} \end{matrix} \qquad (112)$$

Taking partial derivatives of Equations (110) and (111) we obtain the following

$$\frac{\partial A_{mn}}{\partial \mathrm{a}_{rs}} = \frac{\partial A_{ij}}{\partial \mathrm{a}_{rs}}, \qquad \&, \qquad \sum_{n=1}^{3} \frac{\partial A_{nm}}{\partial \mathrm{a}_{rs}} = \frac{\partial \mathrm{a}_m}{\partial \mathrm{a}_{rs}} \qquad (113)$$

There are thus only 9 independent components for the 4th order tensor. The IBOF is developed in terms of the full symmetric 4th order expansion of $\mathrm{a}_{ijkl}$ as a combination of the 2nd order tensor $\mathrm{a}_{ij}$ and identity matrix $\delta_{kl}$ based on Cayley-Hamilton theory is given as

$$\mathrm{a}_{ijkl} = \beta_1 \mathbb{S}(\delta_{ij}\delta_{kl}) + \beta_2 \mathbb{S}(\delta_{ij}\mathrm{a}_{kl}) + \beta_3 \mathbb{S}(\mathrm{a}_{ij}\mathrm{a}_{kl}) + \beta_4 \mathbb{S}(\delta_{ij}\mathrm{a}_{km}\mathrm{a}_{ml}) + \beta_5 \mathbb{S}(\mathrm{a}_{ij}\mathrm{a}_{km}\mathrm{a}_{ml}) + \beta_6 \mathbb{S}(\mathrm{a}_{im}\mathrm{a}_{mj}\mathrm{a}_{kn}\mathrm{a}_{nl}) \qquad (114)$$

where the S operator represents the symmetric permutation expansion of its argument, for example,

$$\mathbb{S}(\mathbb{T}_{ijkl}) = \frac{1}{24}\big[\mathbb{T}_{ijkl} + \mathbb{T}_{ijlk} + \mathbb{T}_{ikjl} + \mathbb{T}_{iklj} + \mathbb{T}_{iljk} + \mathbb{T}_{ilkj} + \mathbb{T}_{jikl} + \mathbb{T}_{jilk} + \mathbb{T}_{jkil} + \mathbb{T}_{jkli} + \mathbb{T}_{jlik} + \mathbb{T}_{jlki} + \mathbb{T}_{kijl} + \mathbb{T}_{kilj} + \mathbb{T}_{kjil} + \mathbb{T}_{kjli} + \mathbb{T}_{klij} + +\mathbb{T}_{klji} + \mathbb{T}_{lijk} + \mathbb{T}_{likj} + \mathbb{T}_{ljik} + \mathbb{T}_{ljki} + \mathbb{T}_{lkij} + \mathbb{T}_{lkji}\big] \tag{115}$$

We obtain the derivative of the 4th order tensor with respect to components of 2nd order tensor by product rule thus.

$$\begin{aligned}\frac{\partial}{\partial \mathrm{a}_{rs}}\{\mathrm{a}_{ijkl}\} = \Big[&\frac{\partial\beta_1}{\partial \mathrm{a}_{rs}}\mathbb{S}(\delta_{ij}\delta_{kl}) + \frac{\partial\beta_2}{\partial \mathrm{a}_{rs}}\mathbb{S}(\delta_{ij}\mathrm{a}_{kl}) + \frac{\partial\beta_3}{\partial \mathrm{a}_{rs}}\mathbb{S}(\mathrm{a}_{ij}\mathrm{a}_{kl}) + \frac{\partial\beta_4}{\partial \mathrm{a}_{rs}}\mathbb{S}(\delta_{ij}\mathrm{a}_{km}\mathrm{a}_{ml}) \\ &+ \frac{\partial\beta_5}{\partial \mathrm{a}_{rs}}\mathbb{S}(\mathrm{a}_{ij}\mathrm{a}_{km}\mathrm{a}_{ml}) + \frac{\partial\beta_6}{\partial \mathrm{a}_{rs}}\mathbb{S}(\mathrm{a}_{im}\mathrm{a}_{mj}\mathrm{a}_{kn}\mathrm{a}_{nl})\Big] + \cdots \\ &+ \Big[\beta_2\mathbb{S}(\delta_{ij}\delta_{kr}\delta_{ls}) + \beta_3\{\mathbb{S}(\delta_{ir}\delta_{js}\mathrm{a}_{kl}) + \mathbb{S}(\mathrm{a}_{ij}\delta_{kr}\delta_{ls})\} \\ &+ \beta_4\{\mathbb{S}(\delta_{ij}\delta_{kr}\delta_{ms}\mathrm{a}_{ml}) + \mathbb{S}(\delta_{ij}\mathrm{a}_{km}\delta_{mr}\delta_{ls})\} \\ &+ \beta_5\{\mathbb{S}(\delta_{ir}\delta_{js}\mathrm{a}_{km}\mathrm{a}_{ml}) + \mathbb{S}(\mathrm{a}_{ij}\delta_{kr}\delta_{ms}\mathrm{a}_{ml}) + \mathbb{S}(\mathrm{a}_{ij}\mathrm{a}_{km}\delta_{mr}\delta_{ls})\} \\ &+ \beta_6\{\mathbb{S}(\delta_{ir}\delta_{ms}\mathrm{a}_{mj}\mathrm{a}_{kn}\mathrm{a}_{nl}) + \mathbb{S}(\mathrm{a}_{im}\delta_{mr}\delta_{js}\mathrm{a}_{kn}\mathrm{a}_{nl}) + \mathbb{S}(\mathrm{a}_{im}\mathrm{a}_{mj}\delta_{kr}\delta_{ns}\mathrm{a}_{nl}) \\ &+ \mathbb{S}(\mathrm{a}_{im}\mathrm{a}_{mj}\mathrm{a}_{kn}\delta_{nr}\delta_{ls})\}\Big]\end{aligned} \tag{116}$$

The $\beta_i$ coefficients are expressed as functions of the second and third invariants (II & III) of the 2nd order tensor $\mathrm{a}_{ij}$. Based on normalization condition and full symmetry requirement coupled with the Cayley-Hamilton theorem, there remains only 3 independent coefficients to determine. The expressions for the IBOF dependent coefficients $(\beta_1, \beta_2, \beta_5)$ are given as

$$\begin{aligned}\boldsymbol{\beta_1} &= \frac{3}{5}\Big[-\frac{1}{7} + \frac{1}{5}\boldsymbol{\beta_3}\Big(\frac{1}{7} + \frac{4}{7}\mathrm{II} + \frac{8}{3}\mathrm{III}\Big) - \boldsymbol{\beta_4}\Big(\frac{1}{5} - \frac{8}{15}\mathrm{II} - \frac{14}{15}\mathrm{III}\Big) \\ &\quad - \boldsymbol{\beta_6}\Big(\frac{1}{35} - \frac{4}{35}\mathrm{II} - \frac{24}{105}\mathrm{III} + \frac{16}{15}\mathrm{II\,III} + \frac{8}{35}\mathrm{II}^2\Big)\Big] \\ \boldsymbol{\beta_2} &= \frac{6}{7}\Big[1 - \frac{1}{5}\boldsymbol{\beta_3}(1 + 4\mathrm{II}) + \frac{7}{5}\boldsymbol{\beta_4}\Big(\frac{1}{6} - \mathrm{II}\Big) - \boldsymbol{\beta_6}\Big(-\frac{1}{5} + \frac{4}{5}\mathrm{II} + \frac{2}{3}\mathrm{III} - \frac{8}{5}\mathrm{II}^2\Big)\Big] \\ \boldsymbol{\beta_5} &= -\frac{4}{5}\boldsymbol{\beta_3} - \frac{7}{5}\boldsymbol{\beta_4} - \frac{6}{5}\beta_6\Big(1 - \frac{4}{3}\mathrm{II}\Big)\end{aligned} \tag{117}$$

We obtain the explicit derivatives of the dependent coefficients via the product rule thus

$$\begin{aligned}\frac{\partial\beta_1}{\partial \mathrm{a}_{rs}} = \frac{3}{5}\Big[&\frac{1}{5}\frac{\partial\beta_3}{\partial \mathrm{a}_{rs}}\Big(\frac{1}{7} + \frac{4}{7}\mathrm{II} + \frac{8}{3}\mathrm{III}\Big) - \frac{\partial\beta_4}{\partial \mathrm{a}_{rs}}\Big(\frac{1}{5} - \frac{8}{15}\mathrm{II} - \frac{14}{15}\mathrm{III}\Big) + \cdots \\ &- \frac{\partial\beta_6}{\partial \mathrm{a}_{rs}}\Big(\frac{1}{35} - \frac{4}{35}\mathrm{II} - \frac{24}{105}\mathrm{III} + \frac{16}{15}\mathrm{II\,III} + \frac{8}{35}\mathrm{II}^2\Big)\Big] + \cdots \\ &+ \frac{3}{5}\Bigg[\Big[\frac{4}{35}\beta_3 - \beta_4 + \beta_6\Big(\frac{4}{35} - \frac{16}{35}\mathrm{II} - \frac{16}{15}\mathrm{III}\Big)\Big]\frac{\partial\mathrm{II}}{\partial \mathrm{a}_{rs}} \\ &+ \Big[\frac{1}{5}\frac{8}{3}\beta_3 + \frac{14}{15}\beta_4 + \beta_6\Big(\frac{24}{105} - \frac{16}{15}\mathrm{II}\Big)\Big]\frac{\partial\mathrm{III}}{\partial \mathrm{a}_{rs}}\Bigg]\end{aligned} \tag{118}$$

$$\begin{aligned}\frac{\partial\beta_2}{\partial \mathrm{a}_{rs}} = \frac{6}{7}\Big[&-\frac{1}{5}\frac{\partial\beta_3}{\partial \mathrm{a}_{rs}}(1 + 4\mathrm{II}) + \frac{7}{5}\frac{\partial\beta_4}{\partial \mathrm{a}_{rs}}\Big(\frac{1}{6} - \mathrm{II}\Big) - \frac{\partial\beta_6}{\partial \mathrm{a}_{rs}}\Big(-\frac{1}{5} + \frac{4}{5}\mathrm{II} + \frac{2}{3}\mathrm{III} - \frac{8}{5}\mathrm{II}^2\Big)\Big] \\ &+ \frac{6}{7}\Big[-\frac{1}{5}[4\boldsymbol{\beta_3} + 7\boldsymbol{\beta_4} + \boldsymbol{\beta_6}(4 - 16\mathrm{II})]\frac{\partial\mathrm{II}}{\partial \mathrm{a}_{rs}} - \frac{2}{3}\boldsymbol{\beta_6}\frac{\partial\mathrm{III}}{\partial \mathrm{a}_{rs}}\Big]\end{aligned}$$

$$\frac{\partial \beta_5}{\partial \mathrm{a}_{rs}} = -\frac{4}{5}\frac{\partial \beta_3}{\partial \mathrm{a}_{rs}} - \frac{7}{5}\frac{\partial \beta_4}{\partial \mathrm{a}_{rs}} - \frac{6}{5}\frac{\partial \beta_6}{\partial \mathrm{a}_{rs}}\left(1 - \frac{4}{3}\mathrm{II}\right) + \frac{8}{5}\boldsymbol{\beta_6}\frac{\partial \mathrm{II}}{\partial \mathrm{a}_{rs}}$$

The independent coefficients $(\beta_3, \beta_4, \beta_6)$ by Chung et al. [48] were obtained from a 5th order binomial fitted function in terms of II & III thus:

$$\beta_{\boldsymbol{m}} = \sum_{i=0}^{5}\sum_{j=0}^{i} a_k^m \cdot \mathrm{II}^{i-j}\mathrm{III}^{j}, \qquad k = j + \frac{1}{2}i(i+1) \tag{119}$$

where the coefficients of the binomial terms can be found in Table 18. The non-unity invariants of $\mathrm{a}_2$ are respectively given as

$$\mathrm{II} = \lambda_1\lambda_2 + \lambda_2\lambda_3 + \lambda_3\lambda_1, \qquad \mathrm{III} = \lambda_1\lambda_2\lambda_3 \tag{120}$$

The derivative of the independent coefficient with respect to the components of the 2nd order tensor is

$$\frac{\partial \beta_{\boldsymbol{m}}}{\partial \mathrm{a}_{rs}} = \sum_{i=0}^{5}\sum_{j=0}^{i} a_k^m \left\{(i-j)\cdot \mathrm{II}^{i-j-1}\mathrm{III}^{j}\frac{\partial \mathrm{II}}{\partial \mathrm{a}_{rs}} + j \cdot \mathrm{II}^{i-j-1}\mathrm{III}^{j-1}\frac{\partial \mathrm{III}}{\partial \mathrm{a}_{rs}}\right\} \tag{121}$$

where,

$$\frac{\partial \mathrm{II}}{\partial \mathrm{a}_{rs}} = (\lambda_2 + \lambda_3)\frac{\partial \lambda_1}{\partial \mathrm{a}_{rs}} + (\lambda_1 + \lambda_3)\frac{\partial \lambda_2}{\partial \mathrm{a}_{rs}} + (\lambda_1 + \lambda_2)\frac{\partial \lambda_3}{\partial \mathrm{a}_{rs}} \tag{122}$$

$$\frac{\partial \mathrm{III}}{\partial \mathrm{a}_{rs}} = (\lambda_2\lambda_3)\frac{\partial \lambda_1}{\partial \mathrm{a}_{rs}} + (\lambda_1\lambda_3)\frac{\partial \lambda_2}{\partial \mathrm{a}_{rs}} + (\lambda_1\lambda_2)\frac{\partial \lambda_3}{\partial \mathrm{a}_{rs}}$$

*Error Estimate*

The performance of the Newton-Raphson (NR) method in accurately predicting the steady-state values of the 2nd order orientation tensor component, is accessed based on the relative absolute error between results of the focus NR method and a reference method, in this case the explicit 4th order Runge-Kutta (RK4) numerical method. We define the error percentage as

$$err = \frac{\mathrm{a}_{mn}^{NR} - \mathrm{a}_{mn}^{RK4}}{\mathrm{a}_{mn}^{RK4}} \times 100\% \tag{123}$$

## Results

We present results of validation carried out for the derived partial derivatives of material derivative for the 2nd order tensor with respect to its components for each model and closure approximations discussed in preceding sections using finite differences. We also present result of the validation for the steady state orientation obtained using the Newton Raphson method by comparing with those obtained using the explicit 4th order Runge-Kutta ODE method. Validation exercise is carried out for different flow conditions.

*Validation of Derivatives based on Finite Difference Approximation*

The results of the validation based on comparison of the Jacobian obtained with the exact derivative to the finite difference approximation is presented below. We present the error defined as the Euclidean norm of the difference between the results obtained from both methods. i.e.

$$err = \left\| \underline{\underline{J}}^{exact} - \underline{\underline{J}}^{FD} \right\|_2 \tag{124}$$

The central difference finite difference approximation is used according to

$$J_{mnij}^{FD} = \frac{\Sigma_{mn}\left(\mathrm{a}_{ij} + \delta\mathrm{a}_{ij}\right) - \Sigma_{mn}\left(\mathrm{a}_{ij} - \delta\mathrm{a}_{ij}\right)}{2\delta\mathrm{a}_{ij}} + O(\delta^2) \tag{125}$$

The model parameters used here can be found in Table 7. The results of the error are shown for different models and closure approximations below. We assume for this validation exercise a 'randomly' generated orientation state $\underline{\underline{a}}^0$ given below:

$$\underline{\underline{a}}^0 = \begin{bmatrix} 0.0622 & 0.0765 & 0.0398 \\ 0.0765 & 0.5521 & 0.0186 \\ 0.0398 & 0.0186 & 0.3857 \end{bmatrix}$$

**Table 4.** Result of error ($\times 10^{-8}$) obtained for different evolution models and different permutation closure approximations.

| | **$HYB_1$** | **$HYB_2$** | **ISO** | **LIN** | **QDR** | **SF2** | **$HL_1$** | **$HL_2$** |
|---|---|---|---|---|---|---|---|---|
| FT | 0.6436 | 0.9385 | 0.2220 | 0.4188 | 0.2691 | 2.0949 | 0.9618 | 4.3940 |
| PT | 0.8088 | 0.7549 | 0.5837 | 0.5003 | 0.4244 | 1.6776 | 0.8241 | 3.4809 |
| iARD | 0.5737 | 1.2712 | 0.3444 | 0.6336 | 0.5728 | 0.5100 | 0.8774 | 1.6148 |
| pARD | 0.7169 | 0.5475 | 0.2722 | 0.2438 | 0.4155 | 1.4805 | 0.9818 | 3.6185 |
| WPT | 0.8563 | 1.0386 | 0.3773 | 0.2525 | 0.2926 | 1.4543 | 0.9632 | 3.5284 |
| Dz | 0.5899 | 0.8248 | 0.2373 | 0.5484 | 0.3233 | 0.7137 | 1.0594 | 2.7732 |
| NEM | 0.6490 | 0.9306 | 0.4012 | 0.4314 | 0.1612 | 2.1012 | 0.9846 | 4.4062 |
| pARD-RSC | 1.0030 | 1.3343 | 1.3062 | 1.0506 | 1.3699 | 0.5378 | 1.3441 | 1.6482 |
| iARD-RPR | 0.5645 | 0.6900 | 0.3478 | 0.3687 | 0.5222 | 1.0731 | 1.0324 | 2.0512 |

**Table 5.** Result of error ($\times 10^{-7}$) obtained for different evolution models and different orthotropic fitted and IBOF closure approximations.

| | **IBOF** | **ORS** | **ORT** | **$NAT_1$** | **ORW** | **$NAT_2$** |
|---|---|---|---|---|---|---|
| FT | 6.1748 | 0.4573 | 0.3351 | 0.3557 | 0.5994 | 0.5812 |
| PT | 5.6437 | 0.3327 | 0.2517 | 0.3076 | 0.5693 | 0.3946 |
| iARD | 4.4942 | 0.2129 | 0.1934 | 0.2512 | 0.4762 | 0.2663 |
| pARD | 5.5458 | 0.3216 | 0.2479 | 0.2975 | 0.5354 | 0.3956 |
| WPT | 5.6412 | 0.3430 | 0.2622 | 0.3040 | 0.5720 | 0.4068 |
| Dz | 6.4226 | 0.2805 | 0.2920 | 0.3532 | 0.6704 | 0.3276 |
| NEM | 6.1978 | 0.4615 | 0.3388 | 0.3649 | 0.6062 | 0.5869 |
| pARD-RSC | 4.1821 | 0.2074 | 0.1802 | 0.2529 | 0.4772 | 0.2800 |
| iARD-RPR | 3.4882 | 0.1601 | 0.1404 | 0.1964 | 0.3687 | 0.1709 |

**Table 6.** Result of error ($\times 10^{-7}$) obtained for different evolution models and different EBOF closure approximations.

| | **WTZ** | **LAR32** | **ORW3** | **VST** | **FFLAR4** | **LAR4** |
|---|---|---|---|---|---|---|
| FT | 4.3147 | 5.0800 | 0.6496 | 3.1567 | 4.3188 | 4.3101 |
| PT | 4.0286 | 4.7162 | 0.5284 | 2.9443 | 4.0435 | 3.9967 |
| iARD | 3.0776 | 3.6782 | 0.4213 | 2.2662 | 3.0115 | 3.0665 |
| pARD | 3.8741 | 4.5415 | 0.5229 | 2.8548 | 3.8500 | 3.8818 |
| WPT | 4.0391 | 4.7234 | 0.5612 | 2.9350 | 4.0486 | 3.9851 |
| Dz | 4.8010 | 5.5454 | 0.6010 | 3.3924 | 4.8016 | 4.6691 |
| NEM | 4.3254 | 5.0936 | 0.6520 | 3.1653 | 4.3271 | 4.3182 |
| pARD-RSC | 2.9401 | 3.4878 | 0.4155 | 2.1358 | 2.9236 | 2.9070 |
| iARD-RPR | 2.4509 | 2.9082 | 0.2820 | 1.7359 | 2.3590 | 2.4110 |

*Validation of NR Method using explicit 4$^{th}$-order Runge-Kutta (RK4) Method*

In this section, results for the steady state values of the preferred orientation states obtained for various cases using the Newton Raphson algorithm are compared to those obtained based on the 4th order explicit Runge-Kutta method. Three (3) sample cases were studied here, the first set of models are based on study by Falvoro et al. [23] and the two (2) other model set were based on study by Tseng et al. [10]. The EBOF closure approximation of Verweyst [32] has been utilized for all analysis. The following data have been used for the different models considered in the first case study [23].

**Table 7.** Case Study 1 parameters for the FT, Dz, iARD, pARD, WPT, MRD and PT models [23]

| | $\boldsymbol{C_I}$ | **ARD Parameters** |
|---|---|---|
| $FT$ | 0.0311 | - |
| $Dz$ | 0.0258 | $D_z = 0.051, \hat{n} = [0 \quad 0 \quad 1]$ |
| $iARD$ | 0.0562 | $C_M = 09977$ |
| $pARD$ | 0.0169 | $\Omega = 0.9868$ |
| $WPT$ | 0.0504 | $w = 0.9950$ |
| $MRD$ | 0.0198 | $\begin{bmatrix} D_1 & D_2 & D_3 \\ 1.000 & 0.7946 & 0.0120 \end{bmatrix}$ |
| $PT$ | - | $\begin{bmatrix} b_1 & b_2 & b_3 & b_4 & b_5 \\ 1.924 & 58.39 & 400 & 0.1168 & 0 \end{bmatrix} \times 10^{-4}$ |

A random orientation state was considered for the initial tensor in the RK4 while for the NR method we consider an initial guess value $\underline{\underline{a}}^0$ for the 2$^{nd}$ order orientation tensor below.

$$\underline{\underline{a}}^0 = \begin{bmatrix} 0.30 & 0.00 & 0.00 \\ 0.00 & 0.60 & 0.10 \\ 0.00 & 0.10 & 0.10 \end{bmatrix}$$

The transient profiles for the component of the 2$^{nd}$ order orientation tensor based on RK4 method for the models presented in Table 7 are shown in Figure 2 below.

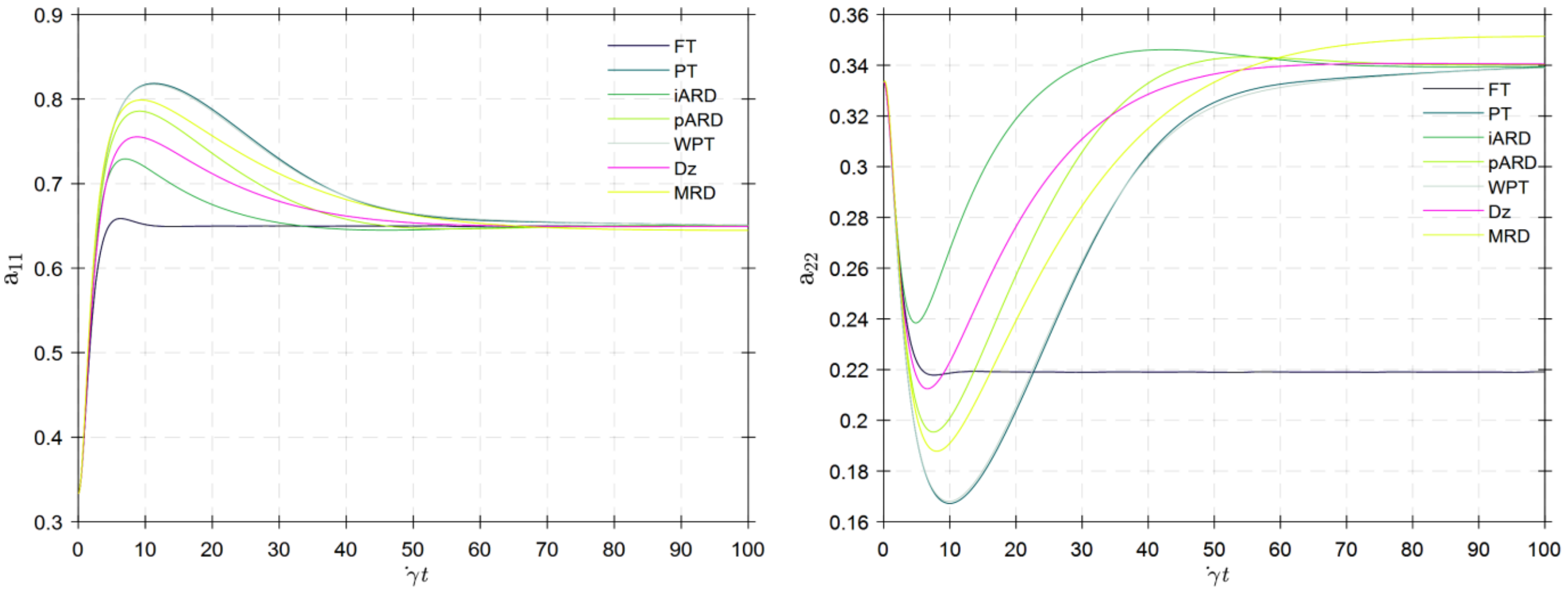

(**a**) (**b**)

**Figure 2.** Time evolution of the 2nd order orientation tensor for calibrated FT, PT, iARD, pARD, WPT, Dz and MRD models for (a) $a_{11}$ component (b) $a_{22}$ component.

Table 8 below shows the result of the error estimate of the steady state values of the orientation tensor components obtained by the NR method for the various models considered using the RK4 values as reference. From the result we see the NR predictions possess good accuracy.

**Table 8.** Error estimates of the $a_{11}, a_{22}$& $a_{13}$ steady-state orientation tensor component values for FT, Dz, iARD, pARD, WPT, MRD and PT models

| | $\mathbf{a_{11}}$ | $\mathbf{a_{22}}$ | $\mathbf{a_{13}}$ |
|---|---|---|---|
| FT | 0.0014 | 0.0009 | 0.0053 |
| PT | 0.0038 | 0.0022 | 0.0067 |
| iARD | 0.0032 | 0.0015 | 0.0033 |
| pARD | 0.0073 | 0.0035 | 0.0534 |
| WPT | 0.0026 | 0.0015 | 0.0099 |
| Dz | 0.0297 | 0.0155 | 0.0086 |

The second case study is based on work by Tseng et al. [10], the calibrated data based on the different model improvements for slow orientation kinetics which they utilized are presented in Table 9. below.

**Table 9.** Case Study 2 parameters for the FT, SRF, RSC and RPR models [10].

| | **FT** | **SRF** | **RSC** | **RPR** |
|---|---|---|---|---|
| $C_I$ | 0.01 | 0.01 | 0.01 | 0.01 |
| ꞔ | – | 0.1 | 0.1 | – |
| $\alpha$ | – | – | – | 0.9 |
| $\beta$ | – | – | – | 0 |

A random orientation state was used as the starting orientation for the RK4 analysis while the initial guess $\underline{\underline{a}}^0$ given below was used for the Newton Raphson method.

$$\underline{\underline{a}}^0 = \begin{bmatrix} 0.35 & 0.00 & 0.00 \\ 0.00 & 0.55 & 0.10 \\ 0.00 & 0.10 & 0.10 \end{bmatrix}$$

Two flow cases were considered:

- Simple shear flow in the 1-2 plane, $L_{12} = \dot{\gamma}$ ($L_1$).
- Balanced shear/planar-elongation flow, simple shear in 1-2 plane superimposed on planar elongation in 1-2 plane. $L_{11} = -\dot{\varepsilon}, L_{22} = \dot{\varepsilon}, L_{12} = \dot{\gamma}$ given $\dot{\gamma}/\dot{\varepsilon} = 10$ ($L_2$).

The time evolution of the components of the 2nd order orientation tensor based on the RK4 method are shown in Figure 3. below.

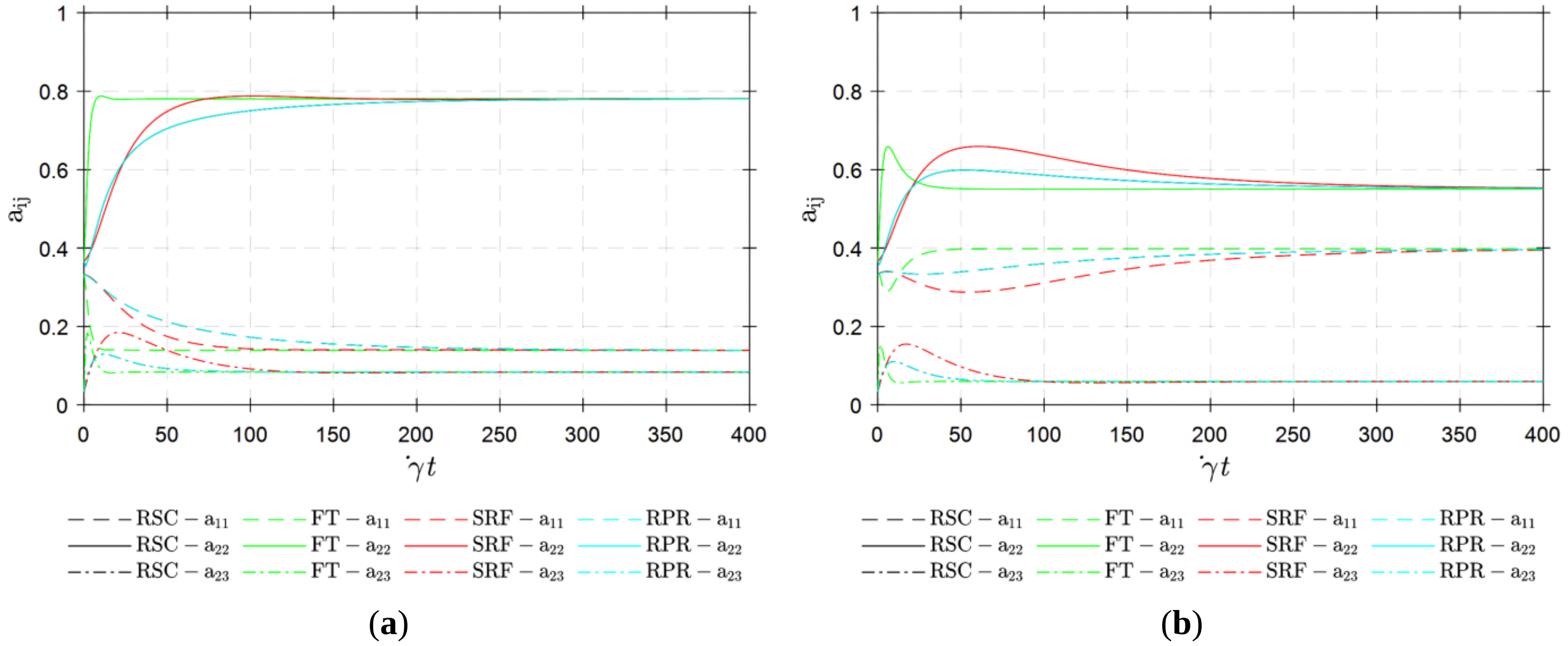

**Figure 3.** Time evolution of the $a_{11}, a_{22}$& $a_{13}$ components of the 2nd order orientation tensor for calibrated FT, SRF, RSC and FT-RPR models for (a) simple shear flow, and (b) shearing/stretching combination flow.

The percentage error estimate between the NR steady state values and the reference RK4 values are presented in Table 10. below. Results show a high level in accuracy in prediction based on the NR method.

**Table 10.** Error estimates of the $a_{11}, a_{22}$& $a_{13}$ steady-state orientation tensor component values for RSC, FT, SRF, and RPR models and for the 2 different flow fields (L1 & L2)

| | **L1** | | | **L2** | | |
|---|---|---|---|---|---|---|
| | $\mathbf{a_{11}}$ | $\mathbf{a_{22}}$ | $\mathbf{a_{13}}$ | $\mathbf{a_{11}}$ | $\mathbf{a_{22}}$ | $\mathbf{a_{13}}$ |
| RSC | 0.0000 | 0.0000 | 0.0000 | 0.0010 | 0.0029 | 0.0634 |
| FT | 0.0079 | 0.0026 | 0.0203 | 0.0000 | 0.0005 | 0.0050 |
| SRF | 0.0022 | 0.0015 | 0.0119 | 0.0010 | 0.0002 | 0.0150 |
| RPR | 0.0000 | 0.0000 | 0.0000 | 0.0000 | 0.0002 | 0.0017 |

In the third case, we consider more complex model development usually involving the combination of two models typically found in injection molding simulation packages such as Moldex3D. The different cases are based on [10] and the model parameter used for the analysis are given in Table 11 & Table 12 below. We assume a random initial orientation state for the reference RK4 method and the same initial guess as with case study 2 for the NR method. The result of the steady state values based on the RK4 method for the different methods are shown in Figure 4. . The percentage error estimate of the NR steady state values with respect to the RK4 reference values are given in Table 13. and the results show negligible discrepancy in values obtained. The results shown in Table 13. reveals good performance in terms of accuracy for the NR method based on the calculated error estimates of the steady state orientation values for the 3-tensor components and for the various models.

**Table 11.** ARD-RSC Parameters [10] for (a) 40 wt. % glass-fiber/PP (b) 31 wt. % carbon-fiber/PP (c) 40 wt. % glass-fiber/nylon

| | **(a)** | **(b)** | **(c)** |
|---|---|---|---|
| ķ | 1/30 | 1/30 | 1/20 |
| $b_1$ | $3.842 \times 10^{-4}$ | $3.728 \times 10^{-3}$ | $4.643 \times 10^{-4}$ |
| $b_2$ | $-1.786 \times 10^{-3}$ | $-1.695 \times 10^{-2}$ | $-6.169 \times 10^{-4}$ |
| $b_3$ | $5.250 \times 10^{-2}$ | $1.750 \times 10^{-1}$ | $1.900 \times 10^{-2}$ |
| $b_4$ | $1.168 \times 10^{-5}$ | $-3.367 \times 10^{-3}$ | $9.650 \times 10^{-4}$ |
| $b_5$ | $-5.000 \times 10^{-4}$ | $-1.000 \times 10^{-2}$ | $7.000 \times 10^{-4}$ |

**Table 12.** *iARD-RPR & pARD-RPR* Parameters [10] for (a) 40 wt. % glass-fiber/PP (b) 31 wt. % carbon-fiber/PP (c) 40 wt. % glass-fiber/nylon

| | **(a)** | **(b)** | **(c)** |
|---|---|---|---|
| $C_I$ | 0.0165 | 0.0630 | 0.0060 |
| $C_M$ | 0.9990 | 1.0100 | 0.9000 |
| $\Omega$ | 0.9880 | 0.9650 | 0.9000 |
| $\alpha$ | 0.9650 | 0.9650 | 0.9500 |
| $\beta$ | 0.0000 | 0.0000 | 0.0000 |

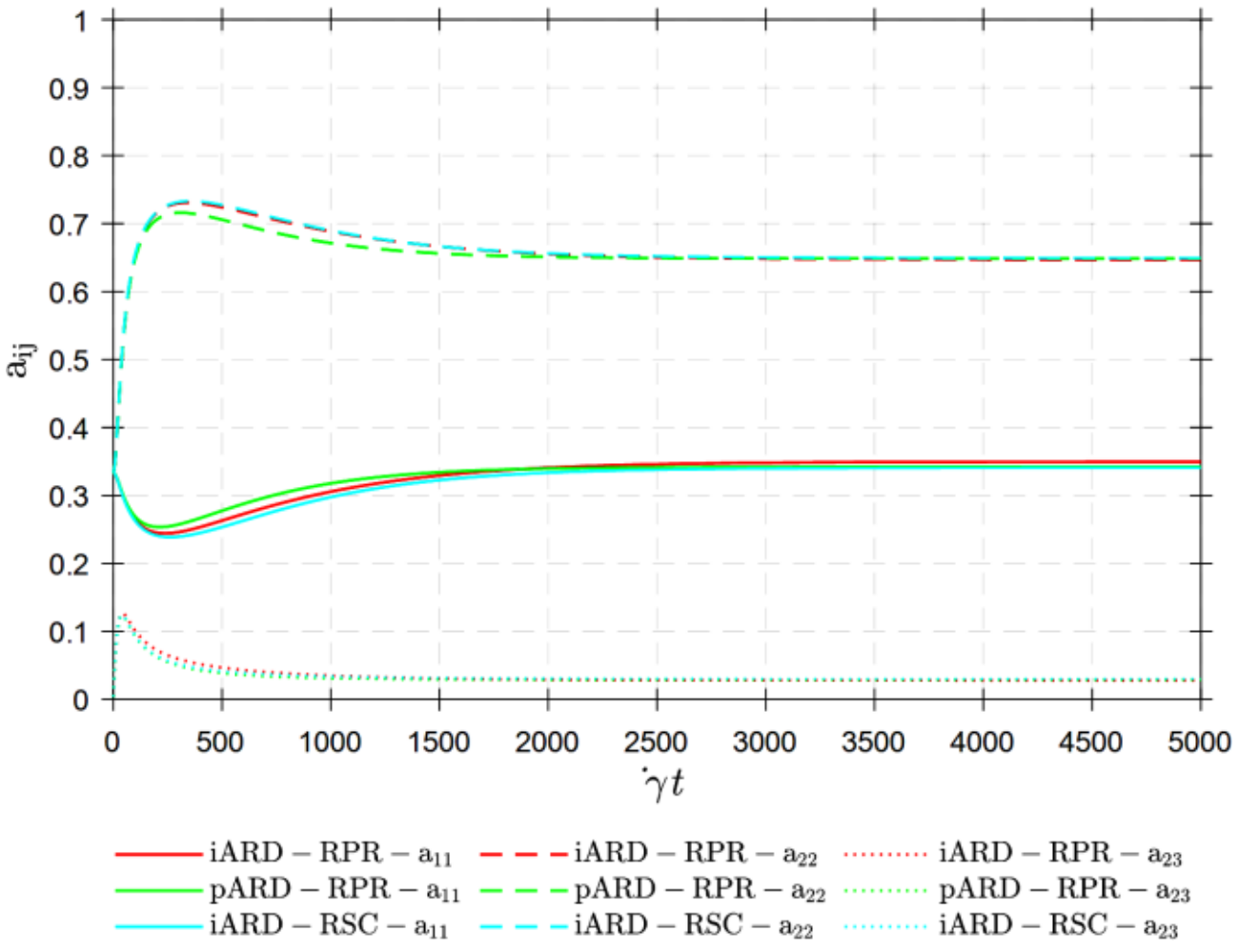

(a)

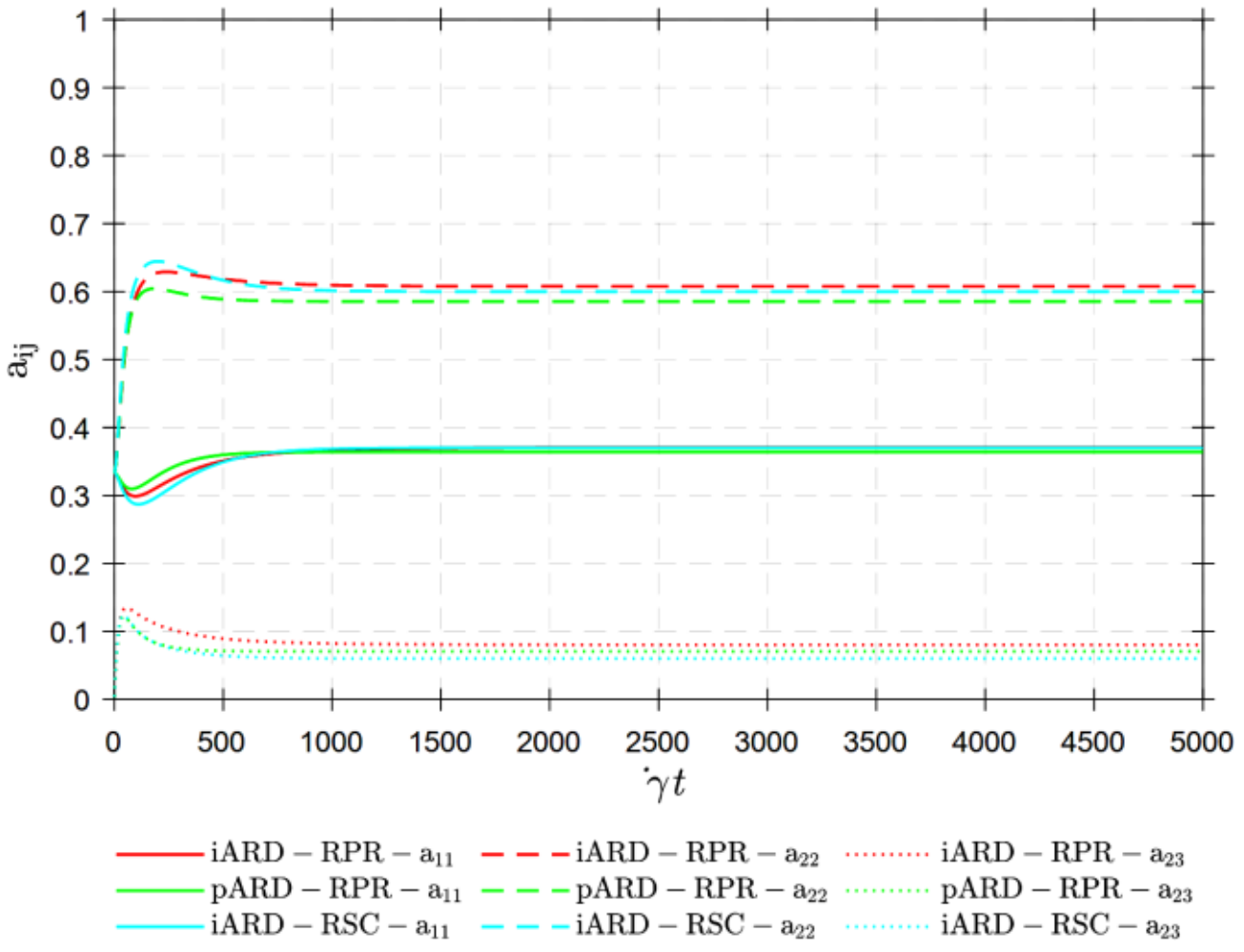

(b)

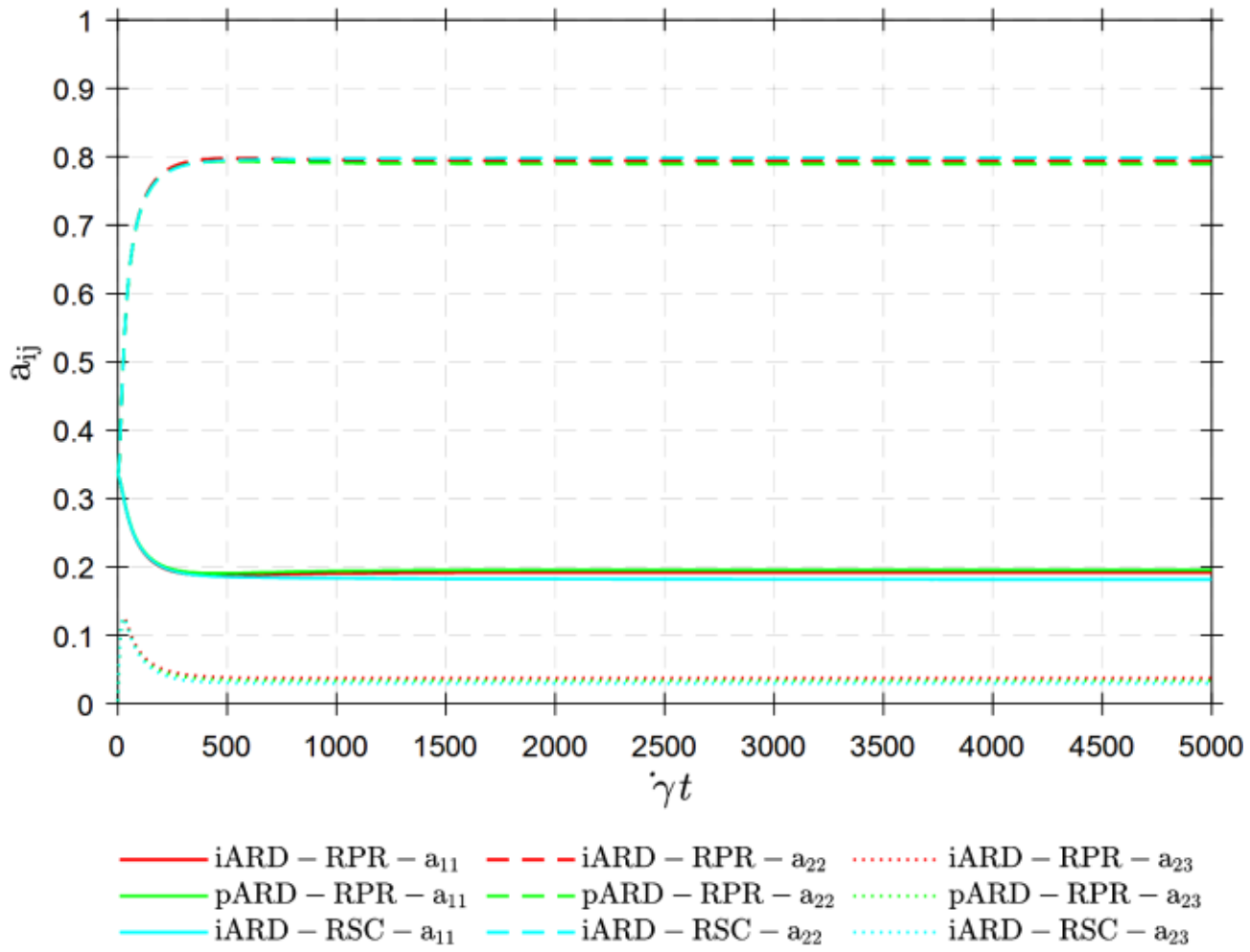

(c)

**Figure 4.** Time evolution of the $a_{11}, a_{22}$& $a_{13}$ components of the 2nd order orientation tensor for calibrated iARD-RPR, iARD-RSC, pARD-RPR models for (a) 40% wt. glass-fiber/PP (b) 31% wt. carbon-fiber/PP (c) 40% wt. glass-fiber/nylon

**Table 13.** Error estimates of the $a_{11}, a_{22}$& $a_{13}$ steady-state orientation tensor component values for iARD-RPR, iARD-RSC, pARD-RPR models for (a) 40% wt. glass-fiber/PP, (b) 31% wt. carbon-fiber/PP, (c) 40% wt. glass-fiber/nylon.

| | | $\mathbf{a_{11}}$ | $\mathbf{a_{22}}$ | $\mathbf{a_{13}}$ |
|---|---|---|---|---|
| | iARD-RPR | 0.0060 | 0.0032 | 0.0036 |
| (a) | pARD-RPR | 0.0009 | 0.0006 | 0.0311 |
| | iARD-RSC | 0.0070 | 0.0037 | 0.0134 |
| | iARD-RPR | 0.0000 | 0.0001 | 0.0053 |
| (b) | pARD-RPR | 0.0000 | 0.0001 | 0.0059 |
| | iARD-RSC | 0.0066 | 0.0014 | 0.0067 |
| | iARD-RPR | 0.0003 | 0.0002 | 0.0012 |
| (c) | pARD-RPR | 0.0003 | 0.0007 | 0.0156 |
| | iARD-RSC | 0.0008 | 0.0017 | 0.0551 |

*Performance of the different Closure Approximations*

The performance of the NR method in obtaining the steady state values for different closure approximations of the 4th order orientation tensor in terms of accuracy and stability has also been assessed. We consider for this assessment the FT model with a $C_I = 0.01$. The initial orientation state for the RK4 reference method is assumed to be random and we assume the same initial guess for the NR method as that of the preceding section. From Table 14, except for the HL2 closure approximation, all other Hinch and Leal closures behaved well. Because of the inherent nature of the transient behavior of the orientation tensor based on the HL2 closure approximation which shows a sudden transition in steady state values at a time fraction of about 100 (cf. Figure 5), we observe a discrepancy in the result for this closure since the NR method has no memory of the history of the orientation state and the accuracy of its prediction is based on the initial guess. The NR method predicts the initial steady state values of $a_{11} = 0.6103, a_{12} = 0.0206$ while the RK4 method transitions to a final steady state orientation of $a_{11} = 0.5759, a_{12} = 0.0467$. The higher order fitted closure approximations behave well in the NR methods and show good accuracy in predictions (cf. Table 15. ).

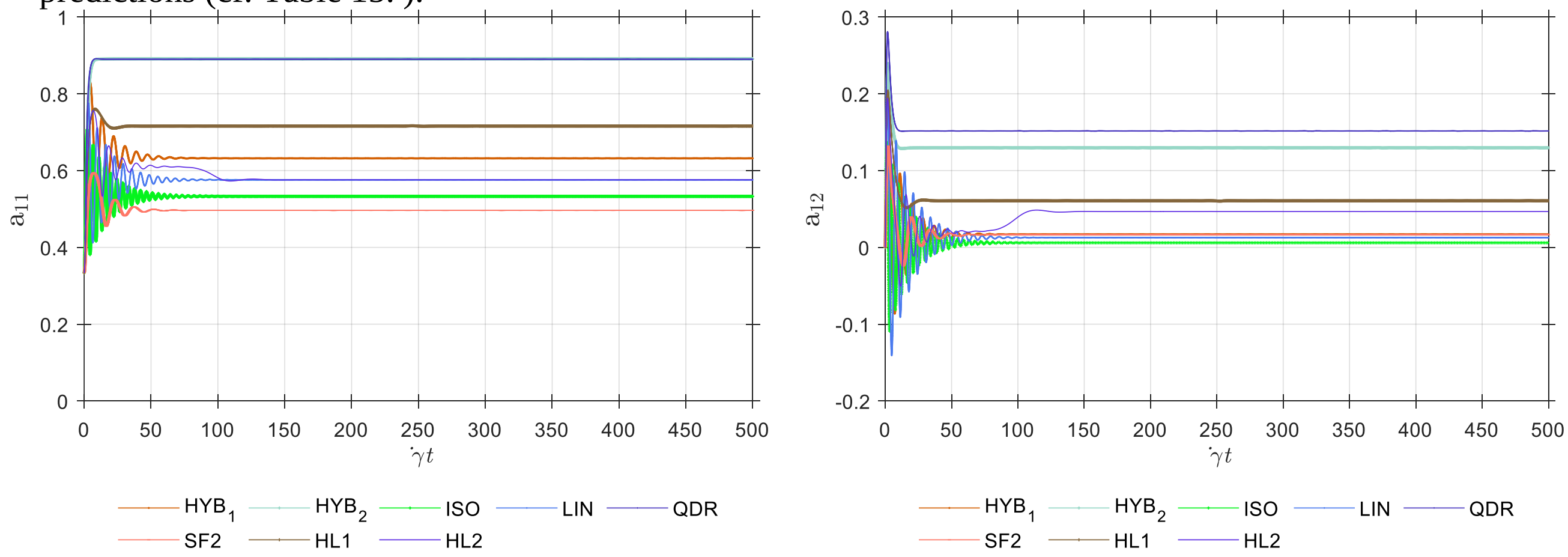

(**a**) (**b**)

**Figure 5.** Transient profiles of 2nd order orientation tensor evolution for (a) component $a_{11}$ and (b) component $a_{12}$ for the various Hinch and Leal closure approximations of the 4th order orientation tensor.

**Table 14.** Error estimates of the $a_{11}, a_{22}$& $a_{12}$ steady-state orientation tensor component values based on the various Hinch and Leal closure approximations of the 4th order orientation tensor.

| | $\mathbf{a_{11}}$ | $\mathbf{a_{22}}$ | $\mathbf{a_{12}}$ |
|---|---|---|---|
| $HYB_1$ | 0.0000 | 0.0005 | 0.2306 |
| $HYB_2$ | 0.0151 | 0.0017 | 0.0223 |
| ISO | 0.0000 | 0.0049 | 0.6305 |
| LIN | 0.0004 | 0.0003 | 0.4655 |
| QDR | 0.0036 | 0.0006 | 0.0053 |
| SF2 | 0.0027 | 0.0101 | 0.0419 |
| HL1 | 0.0042 | 0.0059 | 0.0313 |
| HL2 | 25.9705 | 5.9662 | 55.8952 |

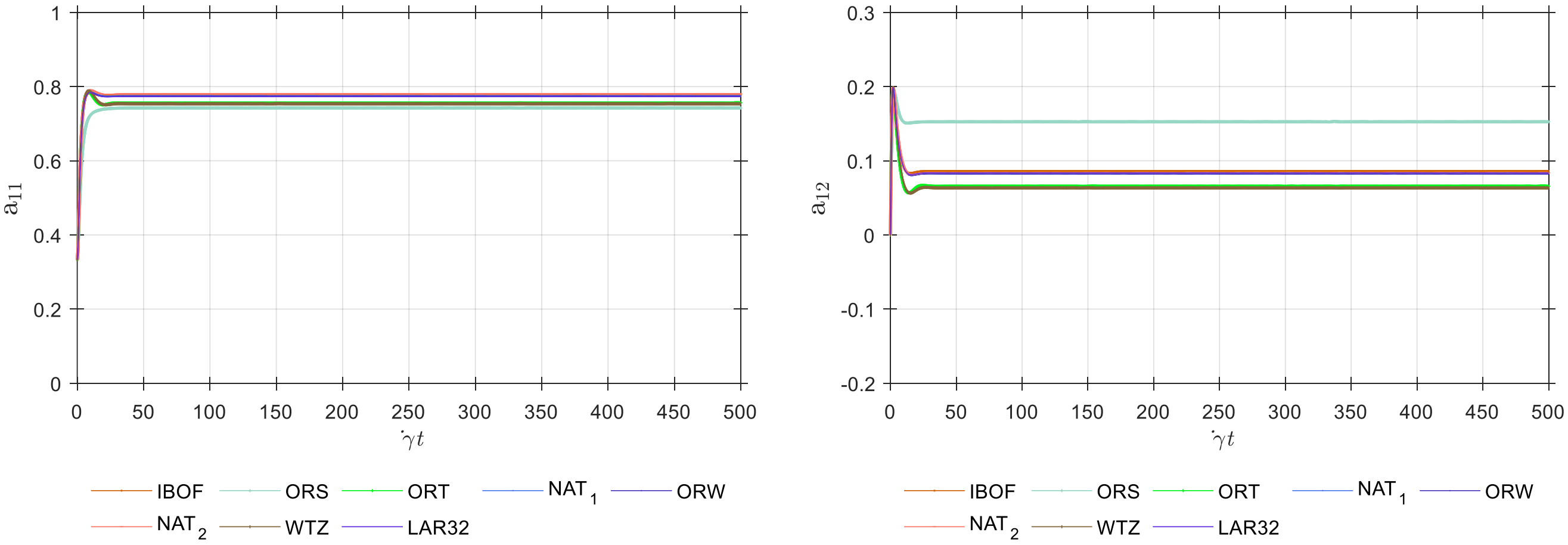

(**a**) (**b**)

**Figure 6.** Transient profiles of 2nd order orientation tensor evolution for (a) component $a_{11}$ and (b) component $a_{12}$ for the higher order orthotropic fitted, IBOF and EBOF closure approximations of the 4th order orientation tensor.

**Table 15.** Error estimates of the $a_{11}, a_{22}$& $a_{12}$ steady-state orientation tensor component values based on the higher order fitted closure approximations of the 4th order orientation tensor.

| | $\mathbf{a_{11}}$ | $\mathbf{a_{22}}$ | $\mathbf{a_{12}}$ |
|---|---|---|---|
| IBOF | 0.0007 | 0.0015 | 0.0232 |
| ORS | 0.0050 | 0.0055 | 0.0151 |
| ORT | 0.0781 | 0.0415 | 0.1600 |
| $NAT_1$ | 0.0000 | 0.0004 | 0.0096 |
| ORW | 0.0007 | 0.0001 | 0.0061 |
| $NAT_2$ | 0.0022 | 0.0017 | 0.0131 |
| WTZ | 0.0000 | 0.0003 | 0.0079 |
| LAR32 | 0.0046 | 0.0013 | 0.0203 |
| ORW3 | 0.0013 | 0.0003 | 0.0036 |
| VST | 0.0257 | 0.0123 | 0.0240 |
| FFLAR4 | 0.0068 | 0.0029 | 0.0327 |
| LAR4 | 0.0020 | 0.0008 | 0.0036 |

*Homogenous Flow Considerations*

We consider different homogenous flows to ensure the stability of the Newtons method in finding stable roots. The following flows were considered:

- Simple Shear (SS), $L_{12} = \dot{\gamma}$
- Two Stretching/Shearing flow (SUA), simple shear in 1-2 plane superimposed with uniaxial elongation in 3-direction. $L_{11} = -\dot{\varepsilon},\ L_{22} = \dot{\varepsilon},\ L_{33} = 2\dot{\varepsilon}\ , L_{12} = \dot{\gamma}$. Two cases consider, balanced shear/stretch, $\dot{\gamma}/\dot{\varepsilon} = 10$, dominant stretch, $\dot{\gamma}/\dot{\varepsilon} = 1$
- Uniaxial Elongation (UA), $L_{11} = 2\dot{\varepsilon},\ L_{22} = L_{33} = -\dot{\varepsilon}$
- Biaxial Elongation, (BA), $L_{11} = L_{22} = \dot{\varepsilon},\ L_{33} = -2\dot{\varepsilon}$
- Two shear/planar-elongation flow (PST), simple shear in 1-3 plane superimposed on planar elongation in 1-2 plane. $L_{11} = -\dot{\varepsilon}, L_{22} = \dot{\varepsilon}, L_{12} = \dot{\gamma}$. Two cases are considered: balanced shear-planar elongation, $\dot{\gamma}/\dot{\varepsilon} = 10$, & dominant planar elongation, $\dot{\gamma}/\dot{\varepsilon} = 1$
- Balanced shear/bi-axial elongation flow, (SBA), simple shear in 1-3 plane superimposed on biaxial elongation. $L_{11} = \dot{\varepsilon}, L_{22} = \dot{\varepsilon}, L_{12} = \dot{\gamma}, L_{33} = -2\dot{\varepsilon}$. A range of $\dot{\gamma}$ is used such that, $2 \leq \dot{\gamma}/\dot{\varepsilon} \leq 5$
- Triaxial Elongation, (TA), $L_{11} = L_{22} = L_{33} = \dot{\varepsilon}$
- Balanced shear/tri-axial elongation flow, (STA), simple shear in 1-3 plane superimposed on biaxial elongation. $L_{11} = \dot{\varepsilon}, L_{22} = L_{33} = \dot{\varepsilon}, L_{12} = \dot{\gamma}$. A range of $\dot{\gamma}$ is used such that, $2 \leq \dot{\gamma}/\dot{\varepsilon} \leq 5$

The initial orientation state for the RK4 reference method is assumed to be random and the initial guess $\underline{\underline{a}}^0$ assumed for each flow consideration is presented in Table 16 below.

**Table 16.** NR initial guess values for different flow conditions.

| **(i)** | **(ii) & (v)** | **(iii)** | **(iv, vii & viii)** | **(vi)** |
|---|---|---|---|---|
| $\begin{bmatrix} 0.35 & 0.00 & 0.00 \\ 0.00 & 0.55 & 0.00 \\ 0.00 & 0.00 & 0.10 \end{bmatrix}$ | $\begin{bmatrix} 0.70 & 0.00 & 0.00 \\ 0.00 & 0.20 & 0.00 \\ 0.00 & 0.00 & 0.10 \end{bmatrix}$ | $\begin{bmatrix} 0.10 & 0.00 & 0.00 \\ 0.00 & 0.10 & 0.00 \\ 0.00 & 0.00 & 0.80 \end{bmatrix}$ | $\begin{bmatrix} 0.40 & 0.00 & 0.00 \\ 0.00 & 0.40 & 0.00 \\ 0.00 & 0.00 & 0.20 \end{bmatrix}$ | $\begin{bmatrix} 0.20 & 0.00 & 0.00 \\ 0.00 & 0.70 & 0.00 \\ 0.00 & 0.00 & 0.10 \end{bmatrix}$ |

Among all closure approximations, the natural closure approximations (exact midpoint fit and extended quadratic fit (cf. Kuzmin [34]), and the Wetzel rational ellipsoid closures behaved well in all flows while the other orthotropic closures had stability issues for one or more of the complex flows and gave non-physical roots. The ability of the NR method to predict accurate results depends on a reasonable initial guess based on the flow type and a suitable closure approximation.

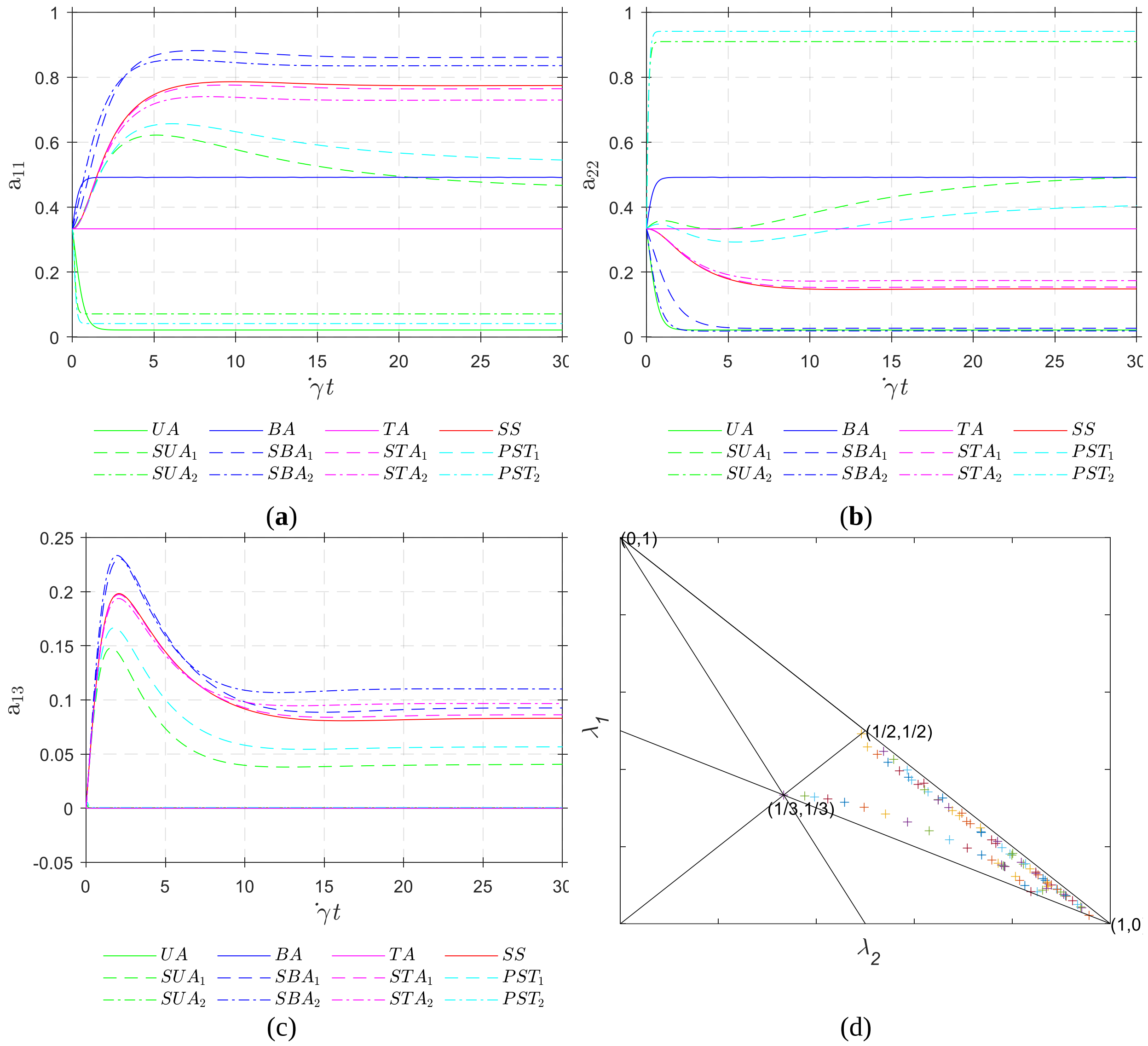

**Figure 7.** Transient profiles of $2^{nd}$ order orientation tensor evolution for (a) component $a_{11}$ and (b) component $a_{22}$ (c) component $a_{12}$ for the various flow considerations. (d) shows the eigenspace for the steady state values obtained for the different flow conditions based on NR method.

**Table 17.** Error estimates of the $a_{11}, a_{22}$& $a_{12}$ steady-state orientation tensor component values for the various flow considerations.

| | $\mathbf{a_{11}}$ | $\mathbf{a_{22}}$ | $\mathbf{a_{12}}$ |
|---|---|---|---|
| SS | 0.0027 | 0.0014 | 0.0012 |
| $SUA_1$ | 0.0004 | 0.0004 | 0.0073 |
| $SUA_2$ | 0.0001 | 0.0000 | 0.0000 |
| UA | 0.0046 | 0.0046 | 0.0000 |
| BA | 0.0175 | 0.0175 | 0.0000 |
| $PST_1$ | 0.0019 | 0.0013 | 0.0053 |
| $PST_2$ | 0.0004 | 0.0024 | 0.1631 |
| $SBA_1$ | 0.0148 | 0.0003 | 0.0000 |
| $SBA_2$ | 0.0053 | 0.0000 | 0.0000 |
| TA | 0.0000 | 0.0000 | 0.0000 |
| $STA_1$ | 0.0084 | 0.0039 | 0.0116 |
| $STA_2$ | 0.0029 | 0.0012 | 0.0124 |

## Conclusion

In conclusion, a Newton-Raphson (NR) method has been successfully implemented in determining the steady state $2^{nd}$ order fiber orientation tensor using exact 4th order Jacobian obtained from partial derivatives of the $2^{nd}$ order fiber orientation tensor material derivative with respect to the $2^{nd}$ order fiber orientation tensor itself. Different macroscopic fiber orientation moment-tensor models and closure approximations of the 4th order fiber orientation tensor are also considered and the performance of the NR method in different homogenous flows have been studied. Like with any typical application of the NR root finding method, a good initial guess of the steady state orientation is required to yield non-physical values. The numerical stability of the NR method depends on the complexity of the flow and the closure approximations. The Natural orthotropic and the IBOF closure approximations performed best for very complex flows. The NR method is comparatively faster compared to the RK4 method. Although obtaining exact derivatives of the $2^{nd}$ order moment-tensor equation of change can be very cumbersome, once they are modelled, they are computationally more efficient since they require less function evaluations compared to a higher order finite difference method of matching accuracy. Moreover, round off error and truncation error may become significant when dealing with relatively small quantities that may lead to instability of the numerical scheme.

**Author Contributions:** “Conceptualization, A.E.A. and D.E.S.; methodology, A.E.A. and D.E.S.; software, D.E.S..; validation, D.E.S.; formal analysis, A.E.A.; investigation, A.E.A.; resources, D.E.S..; data curation, A.E.A.; writing—original draft preparation, A.E.A.; writing—review and editing, D.E.S.; visualization, A.E.A.; supervision, D.E.S.; project administration, D.E.S.; funding acquisition, D.E.S. All authors have read and agreed to the published version of the manuscript.”

**Funding:** This research was funded by the National Science Foundation (NSF) grant number 2055628

**Data Availability Statement:** The data presented in this study are available on request from the corresponding author due to privacy restrictions.

**Conflicts of Interest:** The authors declare no conflicts of interest.

## Appendix A

*Eigenvalues and Eigenvector Derivatives*

The eigenvalue definition for any system is typically given in terms of the eigen-values $\lambda^k$ and corresponding eigen-vectors $\underline{\Phi}^k$ as [43].

$$\left[\underline{\underline{K}} - \lambda_k \underline{\underline{M}}\right] \underline{\Phi}^k = \underline{F}^k \quad \text{(A1)}$$

For most undamped systems, $\underline{F}^k = \mathbf{0}$, and since $\underline{\Phi}^k \neq 0$ to yield non-trivial solutions, then by setting $\left|\underline{\underline{K}} - \lambda_k \underline{\underline{M}}\right| = 0$, we can obtain solutions for $\lambda^k$. By reason of the nature of the system matrix $[K_{mn} - \lambda^k M_{mn}]$ being rank deficient with one order less than the matrix size one may adopt a scaling algorithm to obtain the corresponding eigen-vectors $\underline{\Phi}^k$ usually by defining a Mode $I$ normalization technique for scaling the eigen-vectors $\underline{\Phi}^k$ via a scalar functions $G^k\left(\underline{\Phi}^k\right)$ such that $G^k = 0$ , which may be non-linear in nature. The eigenvectors are thus obtained by replacing the $n^{th}$ row of the residual column vector $\Sigma_i^k = \left(K_{ij} - \lambda^k M_{ij}\right)\Phi_j^k - \underline{F_i^k}$ with $G^k\left(\underline{\Phi}^k\right)$ and solving for $\underline{\Phi}^k$ from the equation $\underline{\Sigma}^k\left(\underline{\Phi}^k\right) = 0$ through any iterative algorithm or explicit solvers. Here we employ Newton-Raphson's method to obtain $\underline{\Phi}^k$ such that:

$$\underline{\Phi}^{k^+} = \underline{\Phi}^{k^-} - \underline{\underline{J}}^{k^{-1}} \Sigma^k \tag{A2}$$

$$\Sigma_i^k = (1 - \delta_{in})\left[S_{ij}^k \Phi_j^k - \underline{F_i^k}\right] + \delta_{in} G^k\left(\underline{\Phi}^k\right) \tag{A3}$$

$$J_{ij}^k = \frac{\partial \Sigma_i^k}{\partial \Phi_j^k} = (1 - \delta_{in}) S_{ij}^k + \delta_{in} \frac{\partial G^k}{\partial \Phi_j^k} \tag{A4}$$

where

$$S_{ij}^k = K_{ij} - \lambda^k M_{ij} \tag{A5}$$

The derivative of the eigenvalues with respect to components $\mathrm{a}_{ij}$ can thus be obtained by differentiating Equation (A3) assuming symmetry of system matrix, i.e., $S_{ij}^k = S_{ji}^k$ such that:

$$\Phi_i^k \frac{\partial S_{ij}^k}{\partial \mathrm{a}_{rs}} \Phi_j^k + \underline{F_i^k} \frac{\partial \Phi_i^k}{\partial \mathrm{a}_{rs}} - \Phi_i^k \frac{\partial F_i^k}{\partial \mathrm{a}_{rs}} = 0 \tag{A6}$$

$$\Phi_i^k S_{ij}^k = S_{ji}^k \Phi_j^k = S_{ij}^k \Phi_j^k = \underline{F_i^k} \tag{A7}$$

where

$$\frac{\partial S_{ij}^k}{\partial \mathrm{a}_{rs}} = \frac{\partial K_{ij}}{\partial \mathrm{a}_{rs}} - \frac{\partial \lambda^k}{\partial \mathrm{a}_{rs}} M_{ij} - \lambda^k \frac{\partial M_{ij}}{\partial \mathrm{a}_{rs}} \tag{A8}$$

Since $\underline{F_i^k} = 0$ for undamped systems

$$\frac{\partial \lambda^k}{\partial \mathrm{a}_{rs}} = \frac{1}{\left(\Phi_i^k M_{ij} \Phi_j^k\right)} \left\{ \Phi_i^k \left[ \frac{\partial K_{ij}}{\partial \mathrm{a}_{rs}} - \lambda^k \frac{\partial M_{ij}}{\partial \mathrm{a}_{rs}} \right] \Phi_j^k \right\} \tag{A9}$$

Consequently, we can obtain the corresponding derivatives $\partial \Phi_i^k / \partial \mathrm{a}_{rs}$ given $\partial \lambda^k / \partial \mathrm{a}_{rs}$ from

$$S_{ij}^k \frac{\partial \Phi_i^k}{\partial \mathrm{a}_{rs}} = -\frac{\partial S_{ij}^k}{\partial \mathrm{a}_{rs}} \Phi_j^k \tag{A10}$$

Recalling the system matrix $\left[\underline{\underline{K}} - \lambda^k \underline{\underline{M}}\right]$ is inherently singular, we can substitute an *n^th^* row of the above equation adopting one the normalization techniques with the equation [43].

$$\frac{\partial G^k}{\partial \Phi_j^k} \frac{\partial \Phi_j^k}{\partial \mathrm{a}_{rs}} = \frac{\partial G^k}{\partial \mathrm{a}_{rs}} \tag{A11}$$

The modified differential equation can be recast as given in equation xx below allowing the inversion of the modified system matrix.

$$J_{ij}^k \frac{\partial \Phi_i^k}{\partial \mathrm{a}_{rs}} = \frac{\partial Q_i^k}{\partial \mathrm{a}_{rs}} \tag{A12}$$

where

$$\frac{\partial Q_i^k}{\partial \mathrm{a}_{rs}} = -(1 - \delta_{in}) \frac{\partial S_{ij}^k}{\partial \mathrm{a}_{rs}} \Phi_j^k - \delta_{in} \frac{\partial G^k}{\partial \mathrm{a}_{rs}} \tag{A13}$$

Smith et. al [43] presents 3 common Mode *I* normalization techniques thus:

1. Mass Normalization

$$G^k\left(\underline{\Phi}^k\right) = \Phi_i^k M_{ij} \Phi_j^k - 1, \qquad \frac{\partial G^k}{\partial \Phi_j^k} = 2\Phi_i^k M_{ij}, \qquad \frac{\partial G^k}{\partial \mathrm{a}_{rs}} = \Phi_i^k \frac{\partial M_{ij}}{\partial \mathrm{a}_{rs}} \Phi_j^k \tag{A14}$$

2. Preassigning an $m^{th}$ Component of $\underline{\Phi}^k$

$$G^k(\underline{\Phi}^k) = \delta_{mj}\Phi_j^k - \alpha, \qquad \frac{\partial G^k}{\partial \Phi_j^k} = \delta_{mj}, \qquad \frac{\partial G^k}{\partial \mathrm{a}_{rs}} = 0 \tag{A15}$$

3. Predefining the Euclidean norm of $\underline{\Phi}^k$

$$G^k(\underline{\Phi}^k) = \sqrt{\Phi_j^k\Phi_j^k} - \beta, \qquad \frac{\partial G^k}{\partial \Phi_j^k} = \Phi_j^k, \qquad \frac{\partial G^k}{\partial \mathrm{a}_{rs}} = 0 \tag{A16}$$

A more direct and efficient approach by Nelson [49] which utilizes mass normalization technique is given below:

$$\begin{aligned} \frac{\partial \Phi_i^k}{\partial \mathrm{a}_{rs}} &= V_i^k + c^k\Phi_i^k \\ c^k &= -\Phi_i^k M_{ij} V_j^k - \frac{1}{2}\Phi_i^k \frac{\partial M_{ij}}{\partial \mathrm{a}_{rs}} \Phi_j^k, \\ V_i^k &= \left\{\underline{\underline{SP}}^{k^{-1}}\right\}_{ij} \frac{\partial QP_j^k}{\partial \mathrm{a}_{rs}} \end{aligned} \tag{A17}$$

Given a pre-selected fixed index – m and noting in Equations (A18)-(A19) below that repeated indices do not imply a summation.

$$SP_{ij}^k = (K_{ij} - \lambda^k M_{ij})(1 - \delta_{mi})(1 - \delta_{mj}) + \delta_{mi}\delta_{mj} \tag{A18}$$

$$\frac{\partial QP_i^k}{\partial \mathrm{a}_{rs}} = -\frac{\partial S_{ij}^k}{\partial \mathrm{a}_{rs}} \Phi_j^k (1 - \delta_{mi}) \tag{A19}$$

## Appendix B

*Optimal Fitted Closure Approximation Constants/Coefficients*

**Eigenvalue Based Optimal Fitting Closure (EBOF) Approximation**

We consider four (4) fitted closures approximations of this form. Linear, quadratic, cubic and quartic binomial fitted closures with $(n+1)(n+2)/2$ number of parameters and their constants below.

1. Linear Orthotropic Fitted Closure $(n = 1)$

For the general linear orthotropic closure, the constant coefficient matrix $\underline{\underline{\mathfrak{C}}}'$ is given as

$$\underline{\underline{\mathfrak{C}}}^{(1)} = \frac{1}{7}\begin{bmatrix} -3/5 & 6 & 0 \\ -3/5 & 0 & 6 \\ 27/5 & -6 & -6 \end{bmatrix}$$

and for the smooth orthotropic closure by Cintra and Tucker (cf. [35]), the constant coefficient matrix $\underline{\underline{\mathfrak{C}}}'$ is given as

$$\underline{\underline{\mathfrak{C}}}^{(1)} = \begin{bmatrix} -0.15 & 1.15 & -0.10 \\ -0.15 & 0.15 & 0.90 \\ 0.60 & -0.60 & -0.60 \end{bmatrix}$$

2. Quadratic Orthotropic Fitted Closures $(n = 2)$

The simple general orthotropic quadratic closure has constant coefficient matrix $\underline{\underline{\mathfrak{C}}}'$ given as

$$\underline{\underline{\mathfrak{C}}}^{(2)} = \begin{bmatrix} 0 & 0 & 0 & 1 & 0 & 0 \\ 0 & 0 & 0 & 0 & 0 & 1 \\ 1 & -2 & -2 & 1 & 2 & 1 \end{bmatrix}$$

The orthotropic natural closure - exact midpoint fit [34] has constant coefficient matrix $\underline{\underline{\mathfrak{C}}}'$ given as

$$\underline{\underline{\mathfrak{C}}}^{(2)} = \begin{bmatrix} 0.0708 & 0.3236 & -0.3776 & 0.6056 & 0.4124 & 0.3068 \\ 0.0708 & -0.2792 & 0.2252 & 0.2084 & 0.4124 & 0.7040 \\ 1.1880 & -2.0136 & -2.1264 & 0.8256 & 1.7640 & 0.9384 \end{bmatrix}$$

The ORF independent coefficients are derived from a 2nd order polynomial fit of the principal axis data obtained from DFC via a minimization technique. For the orthotropic fitted closure by Cintra and Tucker (ORF) (cf. [35]), the constant coefficient matrix $\underline{\underline{\mathfrak{C}}}'$ is given as

$$\underline{\underline{\mathfrak{C}}}^{(2)} = \begin{bmatrix} 0.060964 & 0.371243 & -0.369160 & 0.555301 & 0.371218 & 0.318266 \\ 0.124711 & -0.389402 & 0.086169 & 0.258844 & 0.544992 & 0.796080 \\ 1.228982 & -2.054116 & -2.260574 & 0.821548 & 1.819756 & 1.053907 \end{bmatrix}$$

The ORF has been shown to have better performance compared to non-fitted closure approximations, however, the ORF behaved poorly for flows with very low interaction coefficients and sometimes gave non-physical oscillations similar to that seen with the Hinch and Leal closure (HL2) under the same condition. The ORL behaves better for low interaction coefficient in simple shear flow, but overpredicts flow direction alignment and is unstable for radial diverging flows. Chung and Kwon [50], improved the ORF and developed the 2nd order ORW closure for wide interaction coefficients that proved to be stable for all flow conditions considered. The improved orthotropic fitted closure (ORW2) by Chung and Kwon (cf. [50]), has constant coefficient matrix $\underline{\underline{\mathfrak{C}}}'$ given as

$$\underline{\underline{\mathfrak{C}}}^{(2)} = \begin{bmatrix} 0.070055 & 0.339376 & -0.396796 & 0.590331 & 0.411944 & 0.333693 \\ 0.115177 & -0.368267 & 0.094820 & 0.252880 & 0.535224 & 0.800181 \\ 1.249811 & -2.148297 & -2.290157 & 0.898521 & 1.934914 & 1.044147 \end{bmatrix}$$

Kuzmin [34] presents details on derivations of some orthotropic fitted closures via a numerical bottom up approach.

3. Cubic Orthotropic Fitted Closures ($n = 3$)

Recently higher order polynomial fitted closures were developed for improved accuracy. The orthotropic natural closure - extended quadratic fit though of cubic order is essentially quadratic.

- Non-rational Fitted Closure

The constant coefficient matrix $\underline{\underline{\mathfrak{C}}}^{(3)}$ for this closure approximation is given as

$$\underline{\underline{\mathfrak{C}}}^{(3)} = \begin{bmatrix} 0 & 0.5 & 0 & 0.5 & -0.6 & 0 & 0 & 0.6 & 0.6 & 0 \\ 0 & 0 & 0.5 & 0 & -0.6 & 0.5 & 0 & 0.6 & 0.6 & 0 \\ 1 & -1.5 & -1.5 & 0.5 & 0.4 & 0.5 & 0 & 0.6 & 0.6 & 0 \end{bmatrix}$$

Chung and Kwon [50], also extended the 2nd order ORW to develop 3rd order polynomial ORW3 closure using new flow-data from ODFC calculation. For the improved 3rd order orthotropic fitted closure ORW3 by Chung and Kwon [50], the constant coefficient matrix is given as

$$\left[\underline{\underline{\mathfrak{C}}}^{(3)}\right]^T = \begin{vmatrix} -0.1480648093 & -0.2106349673 & 0.4868019601 \\ 0.8084618453 & 0.9092350296 & 0.5776328438 \\ 0.7765597096 & 1.1104441966 & 0.4605743789 \\ 0.3722003446 & -1.2840654776 & -2.2462007509 \\ -1.7366749542 & -2.5375632310 & -4.8900459209 \\ -1.3431772379 & 0.1260059291 & -1.9088154281 \end{vmatrix}$$

$$
\begin{vmatrix}
-0.0324756095 & 0.5856304774 & 1.1817992322 \\
0.8895946393 & 1.9988098293 & 4.0544348937 \\
1.7367571741 & 1.4863151577 & 3.8542602127 \\
0.6631716575 & -0.0756740034 & 0.9512305286
\end{vmatrix}
$$

Eigenvalue-based optimal fitting (EBOF) closures were developed which are extensions to the ORF but with higher order polynomial fits for improved accuracy. These include the rational elliptical (RE) closure by Wetzel [45] and the 4th order polynomial fitted closure by Verweyst [51] (ORT).

- Rational Fitted Closure

The rational elliptical (RE) closure developed by Wetzel [41] is a higher order extension to the ORF using Carlson elliptical integrals. The rational ellipsoid fitted closure has two permutation vectors, the denominator being one order less than the numerator, i.e., $m = n - 1$ which is of cubic order, $n = 3$. The corresponding constant coefficient matrix for the Wetzel numerator $(n = 3)$ is [45]

$$
\left[\underline{\underline{\mathfrak{C}}}^{(3)}\right]^T =
\begin{vmatrix}
0.1433751825 & 0.1433751825 & 0.9685744898 \\
-0.6566650339 & -0.5209453949 & -2.5526857671 \\
-0.5106016916 & -0.6463213306 & -2.5756669706 \\
3.5295952199 & 0.6031924921 & 2.2044050704 \\
4.4349137241 & 2.3303190917 & 4.4520903005 \\
0.1229618909 & 5.1539592511 & 2.2485545147 \\
-2.9144388828 & -0.2256222796 & -0.6202937932 \\
-5.5556896198 & -1.6481269200 & -1.8811803355 \\
-2.8284365891 & -5.4494528976 & -1.9023485762 \\
0.2292109036 & -3.7461520908 & -0.6414620339
\end{vmatrix}
$$

And for the denominator $(m = n - 1 = 2)$

$$
\left[\underline{\underline{\mathfrak{C}}}^{(2)}\right]^T =
\begin{vmatrix}
1.0000000000 & 1.0000000000 & 1.0000000000 \\
0.7257989503 & 0.6916858207 & -1.2134964928 \\
3.0941511876 & 3.1282643172 & -1.2128608265 \\
-1.6239324646 & -1.5898193351 & 0.2393747647 \\
-4.7303686308 & -4.7303686308 & 0.6004510415 \\
-3.1742364608 & -3.2083495904 & 0.2162486576
\end{vmatrix}
$$

Mullens et al. [42] developed several high order polynomials fitted closures for short fiber composites and introduced the time derivative fitted closures. For the LAR32 closure by Mullens [36] the corresponding constant coefficient matrix for the numerator $(n = 3)$ is

$$\left[\underline{\underline{\mathfrak{C}}}^{(3)}\right]^T = \begin{vmatrix} 0.087602233 & 0.156805152 & 1.072423739 \\ 0.028205550 & -0.577818864 & -2.803554028 \\ -0.426784335 & -0.514280920 & -2.661576129 \\ 1.274677110 & 0.684250887 & 2.389379765 \\ 0.876469059 & 2.132305029 & 4.566728489 \\ 0.602031647 & 3.454835266 & 2.097523143 \\ -1.066583115 & -0.263237143 & -0.658248930 \\ -1.918931146 & -1.614122610 & -1.904704744 \\ -0.934291306 & -4.005261132 & -1.754978355 \\ -0.262854903 & -2.228133231 & -0.508282668 \end{vmatrix}$$

And for the denominator ($m = n - 1 = 2$)

$$\left[\underline{\underline{\mathfrak{C}}}^{(2)}\right]^T = \begin{vmatrix} 1.000000000 & 1.000000000 & 1.000000000 \\ -0.244001948 & 0.365652907 & -1.068512526 \\ -0.574150861 & 1.385725477 & -0.771356469 \\ -0.432097367 & -1.359687152 & 0.067386858 \\ -0.895226091 & -2.866357848 & 0.206908269 \\ -0.462709527 & -1.518996192 & -0.248999874 \end{vmatrix}$$

4. Quartic Orthotropic Fitted Closures ($n = 4$)

The constant coefficient matrix $\underline{\underline{\mathfrak{C}}}'$ based on regression analysis by Verweyst [51] developed from Carlson elliptic integrals.

$$\left[\underline{\underline{\mathfrak{C}}}^{(4)}\right]^T = \begin{vmatrix} 0.6363 & 0.6363 & 2.7405 \\ -1.8727 & -3.3153 & -9.1220 \\ -4.4797 & -3.0371 & -12.2571 \\ 11.9590 & 11.8273 & 34.3199 \\ 3.8446 & 6.8815 & 13.8295 \\ 11.3421 & 8.4368 & 25.8685 \\ -10.9583 & -15.9121 & -37.7029 \\ -20.7278 & -15.1516 & -50.2756 \\ -2.1162 & -6.4873 & -10.8802 \\ -12.3876 & -8.6389 & -26.9637 \\ 9.8160 & 9.3252 & 27.3347 \end{vmatrix}$$

$$
\begin{vmatrix}
3.4790 & 7.7468 & 15.2651 \\
11.7493 & 7.4815 & 26.1135 \\
0.5080 & 2.2848 & 3.4321 \\
4.8837 & 3.5977 & 10.6117
\end{vmatrix}
$$

The constant coefficient matrix $\underline{\underline{\mathfrak{C}}}'$ based on regression analysis for the FFLAR4 closure by Mullens [36]

$$
\left[\underline{\underline{\mathfrak{C}}}^{(4)}\right]^T =
\begin{vmatrix}
0.678225884 & 0.748226727 & 3.167356369 \\
-3.834359034 & -4.249612053 & -13.288266400 \\
-2.664862865 & -2.987266447 & -11.680179330 \\
9.746185193 & 8.641488072 & 23.788431340 \\
14.209962670 & 14.938209410 & 43.700607680 \\
2.700369681 & 5.974489008 & 17.383121430 \\
-8.013024236 & -7.521216405 & -19.959054610 \\
-22.447252700 & -21.757217160 & -58.354308000 \\
-13.078649640 & -15.798676320 & -49.513705640 \\
-0.125467689 & -3.616551654 & -11.755525930 \\
2.417857515 & 2.376441613 & 6.291273472 \\
10.563248410 & 10.222185780 & 25.844317920 \\
12.689484570 & 12.640352670 & 35.425354130 \\
2.487386515 & 4.788201652 & 18.226443930 \\
-0.328195677 & 1.056519961 & 2.925785795
\end{vmatrix}
$$

The constant coefficient matrix $\underline{\underline{\mathfrak{C}}}'$ based on regression analysis for the LAR4 closure by Mullens [36] is

$$
\left[\underline{\underline{\mathfrak{C}}}^{(4)}\right]^T =
\begin{vmatrix}
0.813175172 & 1.768619587 & 4.525066937 \\
-3.065410883 & -9.826017151 & -19.259137620 \\
-4.659333003 & -6.484058476 & -17.650178090 \\
6.329870878 & 19.986994700 & 33.901239610 \\
14.747639770 & 28.905936750 & 61.543979540 \\
9.739797775 & 10.759963010 & 28.467355970 \\
-4.216519964 & -17.715409270 & -27.768082700 \\
-15.922240910 & -40.492387100 & -76.738638810
\end{vmatrix}
$$

| | | |
|---|---|---|
| -20.818571900 | -27.442217500 | -68.977583290 |
| -8.993993112 | -7.230748101 | -22.399036130 |
| 1.138888034 | 5.785725498 | 8.600822308 |
| 5.834142985 | 18.709047480 | 32.480679940 |
| 11.470974520 | 19.729631240 | 43.875135630 |
| 9.874209286 | 8.882877701 | 26.928320210 |
| 3.100457733 | 2.224834058 | 7.101978254 |

**Invariant Based Optimal Fitting Closure (IBOF) Approximation**

The unknown independent coefficients of the binomial expansion for the six parameters in the IBOF closure representation based on regression fitting by Chung et al. [48] of actual flow data obtained from the distribution function closure considering different flow types like EBF closures are given in Table 18. below.

**Table 18.** 5th order binomial fitting coefficients for the IBOF closure approximation

| $k\backslash m$ | *1* | *2* | *3* |
|---|---|---|---|
| *0* | 2.49409081657860E+01 | -4.97217790110754E-01 | 2.34146291570999E+01 |
| *1* | -4.35101153160329E+02 | 2.34980797511405E+01 | -4.12048043372534E+02 |
| *2* | 7.03443657916476E+03 | 1.53965820593506E+02 | 5.73259594331015E+03 |
| *3* | 3.72389335663877E+03 | -3.91044251397838E+02 | 3.19553200392089E+03 |
| *4* | -1.33931929894245E+05 | -2.13755248785646E+03 | -6.05006113515592E+04 |
| *5* | 8.23995187366106E+05 | 1.52772950743819E+05 | -4.85212803064813E+04 |
| *6* | -1.59392396237307E+04 | 2.96004865275814E+03 | -1.10656935176569E+04 |
| *7* | 8.80683515327916E+05 | -4.00138947092812E+03 | -4.77173740017567E+04 |
| *8* | -9.91630690741981E+06 | -1.85949305922308E+06 | 5.99066486689836E+06 |
| *9* | 8.00970026849796E+06 | 2.47717810054366E+06 | -4.60543580680696E+07 |
| *10* | 3.22219416256417E+04 | -1.04092072189767E+04 | 1.28967058686204E+04 |
| *11* | -2.37010458689252E+06 | 1.01013983339062E+05 | 2.03042960322874E+06 |
| *12* | 3.79010599355267E+07 | 7.32341494213578E+06 | -5.56606156734835E+07 |
| *13* | -3.37010820273821E+07 | -1.47919027644202E+07 | 5.67424911007837E+08 |
| *14* | -2.57258805870567E+08 | -6.35149929624336E+07 | -1.52752854956514E+09 |
| *15* | -2.32153488525298E+04 | 1.38088690964946E+04 | 4.66767581292985E+03 |
| *16* | 2.14419090344474E+06 | -2.47435106210237E+05 | -4.99321746092534E+06 |
| *17* | -4.49275591851490E+07 | -9.02980378929272E+06 | 1.32124828143333E+08 |
| *18* | -2.13133920223355E+07 | 7.24969796807399E+06 | -1.62359994620983E+09 |
| *19* | 1.57076702372204E+09 | 4.87093452892595E+08 | 7.92526849882218E+09 |
| *20* | -3.95769398304473E+09 | -1.60162178614234E+09 | -1.28050778279459E+10 |

```
clear, clc
format shortg
[Rung, ~,ddA,dMdA]=BfuncC; v=1; a=3; re=1000;
% A=rand(3); A=1/2*(A+A'); Ad=diag(A,0); Adev=.1*(A-diag(Ad,0));
% Aext=diag(Ad+1/3*(1-trace(A)),0);A=Aext+Adev;Av=[A(1,1:3), A(2,2:3)];
Av=[0.0622,0.0765,0.0398,0.5521,0.0186];
G=1; E=1; dU=zeros(3); dU([1,5, 6, 9])=[-2*E, E, G, E];
% -------------------------------------------------------------------------
mdl_id=[1.0,3.1:.1:3.3, 3.5,3.6, 5.0, 6.3, 7.2] ;
mdl_nm={'FT','PT','iARD','pARD','WPT','Dz','NEM', 'pARD-RSC', 'iARD-RPR'};
cls_id=[1.1, 1.2, 2.1, 2.2, 2.3, 2.4, 2.5, 2.6, 3.0,...
       -1.1,-2.1,-2.3,-2.4,-3.1,-3.2,-3.3,-3.4,-4.1,-4.2,-4.3];
cls_nm=["HYB_1","HYB_2","ISO","LIN","QDR","SF2","HL1","HL2",...
        "IBOF","ORS","ORT","NAT_1","ORW","NAT_2",...
        "WTZ","LAR32","ORW3","VST","FFLAR4","LAR4"];
n1=length(mdl_id); n2=length(cls_id); err=zeros(n1,n2);
%
for j1=1:n1
    jmdl=mdl_id(j1);
    for j2=1:n2
        invars={v, a, @(t) dU ,@(t) [], jmdl,cls_id(j2),mdl(jmdl),{2,4}};
        [~, J1]=ddA([],Av, invars{:});
        [~ ,J2]=dMdA(@(x) ddA([], x, invars{:}), Av, 1.e-6, 6);
        J2=reshape(J2,5,5); err(j1,j2)=norm(J1-J2);
    end
end
arr={1:8,9:14, 15:n2};
for j=1:3
Tj=array2table(err(:,arr{j})); Tj=varfun(@(x) num2str(x, '%.4e'),Tj);
Tj.Properties.VariableNames=cls_nm(arr{j});Tj.Properties.RowNames=mdl_nm;
T{j}=Tj;
end
% -------------------------------------------------------------------------
function f=mdl(flg)
flg1=floor(flg); flg2=round(10*(flg-flg1));
switch flg1
    case 1 % FT - SRF
        f={1. .0311}; % ['kp' 'CI']
    case 2 % RSC
        f={.2 .01}; % ['kp' 'CI']
end
switch flg1
    case 3
        switch flg2
            case 0 % IRD
                f={.0311}; % 'CI'
            case 1 % ARD
                f={[1.924,58.390,400.,.1168,0.]*1e-4}; % 'bta 1-5'
            case {2 4} % iARD
                f={.0562 .9977};      % ['CI' 'CM'] ;
            case 3 % pARD
                c=.9868; D=[1,c,1-c]; D=diag(D,0);
```

```
                f={.0169 D};      % ['CI' 'D' ] ;
            case 5 % WPT
                f={.0504 .9950}; % ['CI', 'w']
            case 6 % Dz
                n=[0,0,1];
                f={.0258, .0051, n}; % ['CI', 'Dz', 'n']
        end
    case 5 % NEM
        f={.01; .063};  % ['CI, U0']
    case 6 % pARD-RSC
        CI=.027; c=.95; kp=.8; D=[1,c,1-c]; D=diag(D,0);f={CI D};f=[{kp} f];
    case 7 % iARD-RPR
        CI=.063; CM=.995; lfa=.2; bta=.01;
        f={lfa bta, CI, CM};
end
end
```

*Published with MATLAB® R2024b*

```
clear,  close, clc
[Rung, Newt]=Bfunc;
Av0N=[.3,.0001,.0001,.6,.1]'; Av0R=1/3*[1.,1.e-4,1.e-4,1.,1.e-4]';
Ar=1000; a=3; t0=0; ti=0.005; tn=500;
mkr={'o','+','*','.','x','_','|','s','d','^','v','>','<','p','h'};
% -------------------------------------------------------------------------
id2=[1.0,3.1:.1:3.3, 3.5,3.6, 4.3] ; nid2=length(id2);
lgn={'FT','PT','iARD','pARD','WPT','Dz','MRD'};
clr2=rand(nid2,3); cls=-4.1; dU=@(t) DU(t,1,1);
for j=1:nid2
    val=mdl(id2(j));
    var={Ar, a, dU, @(t) [],id2(j),cls,val,{2,4}};
           [~,Avn2(j,:)]=Newt([],Av0N,var{:});
    [tn(j),Avn1(j,:),t(:,j),Av(:,:,j)]=Rung(t0,ti,tn,Av0R,var{:});
end
%
for m=1:2
    f=figure(m);clf,f.Color='w'; hold on
    Axx=['\rm\it A_{' num2str(11*(3-m)) '}'];
    title([Axx '-Component']);
    for j=1:nid2
        mkrj=mkr(mod(j-1,15)+1);
        plot(t(:,j),Av(:,3*m-2,j),'Marker',mkrj,'Color',clr2(j,:),...
            'MarkerSize',.5,'LineWidth',.5);
    end
    xlabel('\it\.{$\gamma$}t','Interpreter','latex');
    ylabel(Axx); xlim([0 100]); legend(lgn{:});
    grid on
end
%
err=round(Avn1-Avn2,6)./round(Avn1,6)*100        ;err(isnan(err))=0;
T=array2table(abs(err(:,[1 4 5],1)));T=varfun(@(x) num2str(x, '%.4f'),T);
T.Properties.RowNames=lgn;T.Properties.VariableNames={'A_11','A_22','A_13'};
% -------------------------------------------------------------------------
cls{1}=[1.1,1.2,2.1,2.2,2.3,2.4,2.5,2.6];
dsn{1}=["HYB_1","HYB_2","ISO","LIN","QDR","SF2","HL1","HL2"];

cls{2}=[3.0,-1.1,-2.1,-2.3,-2.4,-3.1,-3.2,-3.3,-3.4,-4.1,-4.2,-4.3];
dsn{2}=["IBOF","ORS","ORT","NAT_1","ORW","NAT_2",...
        "WTZ","LAR32","ORW3","VST","FFLAR4","LAR4"];
G=1; dU=@(t) DU(t,G,1);ncls(1)=length(cls{1});ncls(2)=length(cls{2});
%
for i=1:2
nclsi=ncls(i);clr2=rand(nclsi,3);
for j=1:nclsi
    var={Ar, a, dU, @(t) [], 1.,cls{i}(j),{1., .01},{2,4}};
    [~      ,Avn2(j,:,i)                         ]=Newt([]      ,Av0N,var{:});
    [tn(j,i),Avn1(j,:,i),t(:,j,i),Av(:,:,j,i)]=Rung(t0,ti,tn,Av0R,var{:});
end
end
%
for i=1:2
```

```
    for j=1:2
        f=figure(j+2*(i-1)); clf,f.Color='w'; hold on
        Axx=['\rm\it A_{' num2str(10+j) '}'];
        title([Axx '-Component']);
        for k=1:ncls
            mkrk=mkr(mod(k-1,15)+1);
            plot(G*t(:,k,i),Av(:,j+3,k,i),'Marker',mkrk,'Color',clr2(k,:),...
                'MarkerSize',.5,'LineWidth',.5);
        end
        xlabel('\it\.{$\gamma$}t','Interpreter','latex');
        ylabel(Axx); xlim([0 500]);  ylim([0,1]-[.2,.7]*(j-1));
        grid on; set(gca,'Box','on'); set(gca,'Tickdir','both');
        legend(dsn{i}(:),'Location','south','NumColumns',5,...
            'Orientation','horizontal','FontSize',7);
    end
end

err=round(Avn1-Avn2,6)./round(Avn1,6)*100;
err(isnan(err))=0;err=err(:,[1 4 5],:);
for j=1:2
    Tj=array2table(abs(err(1:ncls(j),:,j)));
    T{j}=varfun(@(x) num2str(x, '%.4f'),Tj);
    T{j}.Properties.VariableNames={'A_11','A_22','A_12'};
    T{j}.Properties.RowNames=dsn{j};
end
```

**-------------------------------------------------------------------------**

**-----**

```
function f=mdl(flg)
flg1=floor(flg); flg2=round(10*(flg-flg1));
switch flg1
    case 1 % FT - SRF
        f={ 1 .0311};
    case 2 % RSC
        f={.2 .0311};
end
switch flg1
case {3 4 6 7}
switch flg2
    case 0 % IRD
        f={.0311}; % 'CI'
    case 1 % ARD
        f={[1.924,58.390,400.,.1168,0.]*1e-4}; % 'bta 1-5'
    case {2 4} % iARD
        f={.0562 .9977};     % ['CI' 'CM'] ;
    case 3 % pARD
        switch flg1
        case 4 % mARD
           %D=zeros(3); D([1 5 9])=[1.,.8,.15];f={.04796 D}; % ['CI' 'D' ] ;
            D=zeros(3); D([1 5 9])=[1.,.7946,.012];f={.0198 D}; % ['CI'
'D' ] ;
```

```
        otherwise % pARD
            c=.9868; D=[1,c,1-c]; D=diag(D,0);
            f={.0169 D};      % ['CI' 'D' ] ;
        end
    case 5 % WPT
        f={.0504 .9950}; % ['CI', 'w']
    case 6 % Dz
        n=[0,0,1];
        f={.0258, .0051, n}; % ['CI', 'Dz', 'n']
end
case 5 % NEM
    f={.0311; 1};  % ['CI, U0']
end
switch flg1
    case 6 % (p)ARD-RSC
        f=[{1/30} f]  ; %['kpa' ARD-vars]
    case 7 % X-RPR
        f=[{1-1/30, .01} f]; %['lfa' 'bta', X-vars]
end
end
```

**-----------------------------------------------------------------------------**

**-----**

```
function f=DU(t,flg,varargin)
f=zeros(3);
switch flg
    case 1 % Simple shear
        G=varargin{1};  f(6)=G;
    case 2 % Shearing/stretching
        [E,G]=varargin{1:2};
        f([1,5,6,9])=[2*E,-E,G,-E];
    case 3 % Uniaxial
        E=varargin{1};
        f([1,5,9])=[-E,-E,2*E];
    case 4 % Biaxial
        E=varargin{1};
        f([1,5,9])=[-2*E, E, E];
    case 5 % Shearing/planar stretching
        [E,G]=varargin{1:2};
        f([3,5,9])=[G,E,-E];

end
end
```

*Published with MATLAB® R2024b*

```
clear,  close, clc
T0=0; h=1.; Tn=5000; [Rung, Newt]=Bfunc;
Ar=1000; a=3;CI=.01;CLS=-3.1;
Av0R=1/3*[1.,1.e-4,1.e-4,1.1,.1]'; Av0N=[.35,.0001,.0001,.55,.1]';
E=.1;G=10*E; dU=zeros(3,3,2); dU(6)=G; dU(9+[1 6 9])=[E G -E];
id=[2 1 1 7.0] ;nid=length(id); kp={{.1}, {1.}, {.1},{.9, .0}};
for i=1:2
    for j=1:nid
        k=j+(i-1)*nid;
        var={Ar, a,@(t) dU(:,:,i),@(t) [],id(j),CLS,[kp{j}, {CI}],{2,3,1}};
        [~,Avn2(j,:,i)]=Newt([],Av0N,var{:});
        [tn(j,i),Avn1(j,:,i),t{i}(:,j),Av{i}
(:,:,j)]=Rung(T0,h,Tn,Av0R,var{:});k
    end
end
figure(1),clf
Lns={'k--', 'g-','r:','c-.'}; pk=[1 4 5];nms={'RSC', 'FT', 'SRF','RPR'};
for i=1:2
    p=subplot(1,2,i); hold on
    for j=1:nid
        for k=1:3
            Axx=[nms{j} '-A_{' num2str(5*(pk(k)-k)+10+k) '}'];
            plot(p,G*t{i}(:,j),Av{i}(:,pk(k),j),Lns{j},'DisplayName',Axx,...
                'LineWidth',.5);
        end
    end
    xlabel('\it\.{$\gamma$}t','Interpreter','latex');
    ylabel(p,'\it A_{xx}'); xlim(p,[0 400]);ylim(p,[0,1]);
     p.FontSize=6;legend(p,'Location','north','Orientation','horizontal',...
        'NumColumns',3,'FontSize',5); grid on;
end
err=round(Avn1-Avn2,6)./round(Avn1,6)*100; err(isnan(err))=0;
T=array2table(abs([err(:,[1 4 5],1) err(:,[1 4 5],2)]));
T=varfun(@(x) num2str(x, '%.4f'),T);
T.Properties.VariableNames=...
    {'A_11_L1','A_22_L1','A_13_L1','A_11_L2','A_22_L2','A_13_L2'};
T.Properties.RowNames={'RSC', 'FT', 'SRF','RPR'};
% -----------------------------------------------------------------------
Av0N=[.35,.0001,.0001,.55,.1]'; Av0R=1/3*[1.,1.e-4,1.e-4,1.,.01]';
Ar=1000; a=3;CLS=-4.1; G=1;
dU=zeros(3); dU(6)=G;
CI={.0165;.0630;.0060}; ck={.988; .965; .900}; CM={.999; 1.010; .900};
kpa={1/30; 1/30; 1/20}; ak=num2cell(1-[kpa{:}]'); bk={0.0;0.0;0.0};
bta={[3.842 -17.86  525.  .1168  -5. ]*1e-4;
     [37.28 -169.5  1750. -33.67 -100]*1e-4;
     [4.643 -6.169  190.    9.65  7  ]*1e-4};
for j=1:3, Dk{j,1}=diag([1, ck{j}, 1-ck{j}]); end
val={[ak, bk, CI, CM], [ak, bk, CI, Dk], [kpa, bta]};
id=[7.4, 7.3, 6.1] ;nid=length(id);
for i=1:3
    for j=1:nid
        k=j+nid*(i-1);
```

```
        var={Ar, a,@(t) dU,@(t) [],id(j), CLS,val{j}(i,:),{2,4}};
        [~, Avn2(k,:)]=Newt([],Av0N,var{:});[i,j]
        [tn(k),Avn1(k,:),t{i,j},Av{i,j}]=Rung(T0,h,Tn,Av0R,var{:});
    end
end
f=figure(1);f.Color='w'; clf, hold on
Lns={'x--', 'p:','+-.'}; pk=[1 4 5];
nms={'iARD-RPR','pARD-RPR', 'iARD-RSC'}; clr={'r','g','c'};
ttl={'40%wt.glass-fiber/PP', '31%wt.carbon-fiber/PP', '40%wt.glass-fiber/
nylon'};
for i=1:3
    p=subplot(2,2,i); hold on
    title(p,['\rm\it' ttl{i}]);
    for j=1:nid
        for k=1:3
            Axx=[nms{j} '-A_{' num2str(5*(pk(k)-k)+10+k) '}'];
            plot(p,G*t{i,j},Av{i,j}(:,pk(k)),Lns{k},'DisplayName',Axx,...
                'LineWidth',.5,'Color',clr{j},'MarkerSize',1);
        end
    end
    xlabel('\it\.{$\gamma$}t','Interpreter','latex');
    ylabel(p,'\it A_{xx}'); xlim(p,[0 G*Tn]);ylim(p,[0,1]);
     p.FontSize=6;legend(p,'Location','north','Orientation','horizontal',...
        'NumColumns',3,'FontSize',4); grid on;
end

err=round(Avn1-Avn2,6)./round(Avn1,6)*100; err(isnan(err))=0;
T=array2table(abs(err(:,[1 4 5])));T=varfun(@(x) num2str(x, '%.4f'),T);
T.Properties.VariableNames={'A_11','A_22','A_13'};
Tn=['iARD-RPR';'pARD-RPR'; 'iARD-RSC'];nTn=size(Tn,1);
Tn=[repmat(Tn,3,1),repelem(num2str((1:3)'),nTn,1)];nTn=size(Tn);
T.Properties.RowNames=mat2cell(Tn,ones(nTn(1),1),nTn(2));
%------------------------------------------------------------------------
```

*Published with MATLAB® R2024b*

```
clear,  close, clc
[Rung, Newt]=BfuncR;
a=3; cls=3.0; Av0=1/3*[1.,0.,0.,1.,0.]'; tol=1.e-2;
v=1.;  r0=1; rn=30; nz=31; z=linspace(-1.+tol,1.-tol,nz);
%
%
mdl_id=[6.1,7.4, 7.3]; nid=length(mdl_id);
mdl_pr={
        {1/30 [3.842 -17.86  525.  .1168  -5. ]*1e-4} ; % kpa bta(i)   - ARD-
RSC
        {1-1/30 0. .0165 .999                       } ; % lfa bta CI CM-iARD-
RPR
        {1-1/30 0. .0165 diag([1. .988 0.012],0)    }}; % lfa bta CI  D-pARD-
RPR
%
for i=1:nid
    for j=1:nz
        var={ v, a, @DU, z(j),mdl_id(i),cls,mdl_pr{i},{2,4}};
        [rRK(j,i), AvRK(j,:,i)]=Rung(r0,.001,rn,Av0,var{:});[i j]
    end
end
%
f=figure(1);clf, f.Color='w'; hold on
Lns={'r--', 'c:','k-.'};nms={'ARD-RSC','iARD-RPR','pARD-RPR'};
for i=1:4
    s=1/2*(i-1)*(i-2); k=s/3*(i-3);
    for j=1:3
        ydat=(1-k)*(s+(2-i)*AvRK(:,1,j)+(i*(1-s)-1)*AvRK(:,4,j))
+k*AvRK(:,5,j);
        plot(z,ydat,Lns{j},'DisplayName',['\it' nms{j}]);
    end
end
p=gca; p.XTick=linspace(-1.,1.,9);p.YTick=linspace(-.2, 1,9);
p.XAxisLocation="bottom";p.YAxisLocation="origin";p.Box='on';
grid(p,"on"); xlabel(p,'\it z/h');ylabel(p,'\it A_{xx}');
p.XAxis.TickDirection="out"; p.YAxis.TickDirection="both";
p.YAxis.TickLabelFormat='%.2f';p.FontSize=7;
legend(p,nms{:},'FontSize',8,'Location','northeast','Box','off');
xlim(p,[-1 1]);ylim(p,[-0.2, 1]);
%
aX=[.38 .325;.65 .575;.65 .6;.75 .7];
aY=[0.72 0.63;0.7 0.615; 0.4 0.28; 0.3 0.22];
for j=1:4
    s=1/6*(j-1)*(j-2)*(j-3); Ann=num2str(11*j*(1-s)+13*s);
    annotation('textarrow',aX(j,:),aY(j,:),'String',['\it A_{' Ann '}'],...
        'FontSize',7,'HeadStyle','deltoid','HeadWidth',5);
end
%---------------------------------------------------------------------------
function [f,df]=DU(r,z)
Q=100; b=1.5e-1; sz=sign(z);
vb=(3*Q)/(8*pi*b^2); vr=1/r*(1-z^2) ; vz=sz*z/r;
f =zeros(3);   f([5,1,6])=vb/b*[-vr/r,vr/r, -2*vz];
```

```
if nargout>1
df=zeros(3);  df([5,1,6])=(vb/b^2)*(-2/r)*[-vr/r,vr/r, -vz];
end
end
```

*Published with MATLAB® R2024b*

```
clear,  close, clc
[Rung, Newt]=Bfunc;
Av0R=1/3*[1.,1.e-4,1.e-4,1.,1.e-4]'; Ar=1000; a=3; t0=0; ti=0.01; tn=100;
mkr={'o','+','*','.','x','_','|','s','d','^','v','>','<','p','h'};
lgd={'SH','SH-ST-1','SH-ST-2','UA','BA','SH-PE-1','SH-PE-2',...
    'SH-BE-1','SH-BE-2','TA','SH-TE-1','SH-TE-2'};
% ------------------------------------------------------------------------
id2=1.0; kp=1.; CI=0.01; val2={kp CI}; cls=-2.3; % -2.3, -3.1, -3.2, 3.0
% ------------------------------------------------------------------------
val1={1;[0.1, 1;10 1]; 1; 1;[0.1, 1;10, 1];[0.2, 1;.5, 1];1;[0.2, 1;.5, 1]};
nf=length(val1); jj=0;
for j=1:nf
    valj=val1{j}; ns(j)=size(valj,1);
    for k=1:ns(j)
        jj=jj+1; [dU,Av0N]=DU(j,valj(k,:));
        var={Ar, a,@(t) dU , @(t) [],id2,cls,val2,{2,4}};
        [~     , Avn2(jj,:)                   ]=Newt([],Av0N,var{:}); [j,k]
        [tn(jj), Avn1(jj,:),t(:,jj),Av(:,:,jj)]=Rung(t0,ti,tn,Av0R,var{:});
    end
end
clr2=rand(jj,3);
%
for m=1:3
    f=figure(m);clf,f.Color='w'; hold on
    n=-(m^2-6*m+4); Axx=['\rm\it A_{' num2str(.5*(13*m^2-61*m+92)) '}'];
    title([Axx '-Component']);  jj=0;
    for j=1:nf
        nf2=ns(j);
        for k=1:nf2
            jj=jj+1; mkrj=mkr(mod(jj-1,15)+1);
            plot(t(:,jj),Av(:,n,jj),'Marker',mkrj,'Color',clr2(jj,:),...
                'MarkerSize',1,'LineWidth',.5,'DisplayName',...
                ['\it' lgd{j} '-' num2str(k)]);
        end
    end
    xlabel('\it\.{$\gamma$}t','Interpreter','latex');
    ylabel(Axx); xlim([0 30]);
    legend(lgd{:},'Location','best');
    set(gca,"Box","on");set(gca,"TickDir",'both'); grid on
end
%
err=round(Avn1-Avn2,6)./round(Avn1,6)*100        ;
err(isnan(err))=0; T=array2table(abs(err(:,[1 4 5]))) ;
T=varfun(@(x) num2str(x, '%.4f'),T);
T.Properties.VariableNames={'A_11','A_22','A_13'};
T.Properties.RowNames=lgd;
% ------------------------------------------------------------------------
nv=21; va=linspace(.01,8,nv)'; vb=logspace(-1.5,2,nv)';
v1=ones(nv,1); v2=1./va;  v3=1./vb;
val1={1;[v2 v1]; 1; 1;[v2 v1];[v3 v1];1;[v3 v1]};nf=length(val1);

jj=0;
```

```
for j=1:nf
    valj=val1{j}; ns(j)=size(valj,1);
    for k=1:ns(j)
        jj=jj+1; [dU,Av0N]=DU(j,valj(k,:));
        var={Ar, a,@(t) dU , @(t) [],id2,cls,val2,{2,4}};
        [~      , Avn2(jj,:)                    ]=Newt([],Av0N,var{:});
    end
end
f=figure(4); f.Color='w';
for j=1:size(Avn2,1)
    v=sort(eig(v2M(Avn2(j,:),1)),'descend');
    plot(v(1),v(2),'+'); hold on
end
X=[
    0    1    1    0
    0  0.5    1  0.5
    0  0.5    1    0
    0    1  0.5    0
    0  0.5    0  0.5];
for j=1:5
    line(X(j,1:2),X(j,3:4),'Color','k')
end
text(1/3,1/3,'(1/3,1/3)','FontSize',8);
text(1/2,1/2,'(1/2,1/2)','FontSize',8);
text(1  ,0  ,'(1,0)'    ,'FontSize',8);
text(0  ,1  ,'(0,1)'    ,'FontSize',8);
set(gca,'Box','on');
ax=gca; ax.XAxis.TickLabels='';ax.YAxis.TickLabels='';
xlabel('\it \lambda_{2}');ylabel('\it \lambda_{1}');
% -------------------------------------------------------------------------
function [f, Av0]=DU(flg,val)
f=zeros(3);
switch flg
    case 1 % Simple shear
        G=val(1);  f(6)=G;
        Av0=[.35,1.e-4,1.e-4,.55,1.e-4]';
    case 2 % Shearing/stretching
        E=val(1); G=val(2); f([1,5,6,9])=[2*E, E, G, -E];
        Av0=[.7,1.e-4,1.e-4,.2,1.e-4]';
    case 3 % Uniaxial
        E=val(1); f([1,5,9])=[-E,-E,2*E];
        Av0=[.1,1.e-4,1.e-4,.1,1.e-4]';
    case 4 % Biaxial
        E=val(1); f([1,5,9])=[E, E,-2*E];
        Av0=[.4,1.e-4,1.e-4,.4,1.e-4]';
    case 5 % Shearing/planar stretching
        E=val(1);G=val(2); f([1, 6, 9])=[E,  G, -E];
        Av0=[.7,1.e-4,1.e-4,.2,.01]';
    case 6 % Balanced shear/biaxial elongation flow
        E=val(1);G=val(2); f([1,5, 6, 9])=[-2*E, E, G, E];
        Av0=[.1,1.e-4,1.e-4,.8,1.e-4]';
    case 7 % Triaxial
        E=val(1); f([1,5,9])=[E, E, E];
        Av0=[.4,1.e-4,1.e-4,.4,1.e-4]';
```

```
    case 8 % Balanced shear/triaxial elongation flow
        E=val(1);G=val(2); f([1,5, 6, 9])=[E, E, G, E];
        Av0=[.4,1.e-4,1.e-4,.4,1.e-4]';
end
end

function A=v2M(Av,flg)
for m=1:2
    for n=m:3
        k=2*(m-1)+n;
        switch flg
            case 1
                A(m,n)=Av(k); A(n,m)=A(m,n);
            case 2
                A(k,1)=Av(m,n);
        end
    end
end
if flg==1, A(3,3)=1-A(1,1)-A(2,2); end
end
```

*Published with MATLAB® R2024b*

```
function varargout=Bfunc
varargout={@Rung, @Newt,@ddA,@dMdA};
end


ans =

  function_handle with value:

    @Rung


function [tn, Avn,t,Av]=Rung(t0,ti,tn,Av0,re,varargin)
v=(re^2-1)/(re^2+1); ztol=5;
% opts = odeset('Events',@dzero);
% [t, Av]=ode45(@(t,y) ddA(t,y,v,varargin{:}),t0:ti:tn,Av0,opts);
[t, Av]=ode45(@(t,y) ddA(t,y,v, varargin{:}),t0:ti:tn,Av0);
tn=t(end); Avn=Av(end,:);
    function [pos,ter,dir] = dzero(t,y)
        pos = round(norm(ddA(t,y, v,varargin{:}),2),ztol);
        ter = 1;
        dir = 0;
    end
end

function [t,Av]=Newt(t,Av,re,varargin)
v=(re^2-1)/(re^2+1);
err=1;
while err>1e-5
    [R, J]=ddA(t,Av, v, varargin{:});
%     [R,J2]=dMdA(@(x) ddA(t,x, v, varargin{:}), Av, 1.e-5, 2);
%     J2=reshape(J2,5,5); norm(J-J2);
    Av=Av-J\R; err=norm(R);
end
%
end
% upper case for matrix, lower case for scalar, v is lambda

function [R, J]=ddA(t,Av, varargin)
[v, a, DU, DUr, flg1,flg2,var1,var2]=varargin{:};
dU=DU(t); dUr=DUr(t);
%
A=v2M(Av,1); I=eye(3);
if nargout>1
    dA=DA(I); s3=size(dA,3); dI=zeros(3,3,s3);
end
%
if nargout>1
    [F,dF]=mdlR(A, I, dA, dI, v, a, dU, dUr,flg1,flg2,var1,var2);
    R=[F(1,:) F(2,2:3)]';J=permute([dF(1,:,:) dF(2,2:3,:)],[2,3,1]);
else
    F     =mdlR(A, I, [], [], v, a, dU,  [],flg1,flg2,var1,var2);
    R=[F(1,:) F(2,2:3)]';
```

```
end
if nargout>1
    % Rank Deficient - Normalization
    nc=size(J,2); nr=rank(J);
    if (nc-nr)>0
        if nargout>1
            Rn=norm(R);  Jn=1/Rn*(R'*J);  R =[Rn; R]; J= [Jn;J] ;
        end
    end
end
%
end

function [F,dF]=mdlR(A, I, dA, dI, v, a, dU, dUr, flg1,flg2,var1,var2)
%
flg1_1=floor(flg1); flg1_2=round(10*(flg1-flg1_1));
bol1=(nargout>1); bol2=isempty(dUr); bol=logical(abs(bol1-bol2));
%
switch flg1_1
    case {1 2 6} % SRF RSC
        kp=var1{1};
        switch flg1_1
            case {1 2} % IRD
                CI=var1{2}; var1=[];
            case 6 % ARD
                var1=var1(2:end);
        end
    case 5 % NEM
        [U0, CI]=var1{1:2}; var1=[];
end
switch flg1_1
    case 7 % RPR
        [ak, bk]=var1{1:2};
        switch flg1_2
            case 0
                CI=var1{3}; var1=[];
            otherwise
                var1=var1(3:end);
        end

end
%
switch floor(flg2)
    case {0 1}
        var2={a};
end
%
W =1/2*(dU'-dU); Y=1/2*(dU+dU');  y=sqrt(2*sum(Y.* Y' ,'all'));
if bol1
    s3 =size(dA,3);
    if ~bol2
    dWr=1/2*(dUr'-dUr); dYr=1/2*(dUr+dUr'); dyr=2/y*sum(Y.*dYr','all');
    end
end
```

```
%------------------------------------------------------------------------------
switch flg1_2
    case 4
        Yin=dU;
    otherwise
        Yin=Y;
end
if bol
    varin_1={A,I,[],[],Yin,[],flg1,flg2,var1,var2};
    switch flg1_1
        case {1 7}
            if flg1_2==0
                varin_1{7}=1;  YA4 =BA4(varin_1{:});
            end
    end
    switch flg1_1
        case 2
            [YA4, YLMA4]=BA4(varin_1{:});
        case {3 5 7}
            if flg1_2>0 || flg1_1==5
                [YA4, CA4,C]=BA4(varin_1{:});
            end
        case 4
            [YA4, C    ]=BA4(varin_1{:});
        case 6
            [YA4, YLMA4, CA4, CM4, CLMA4, C]=BA4(varin_1{:});
    end
end
if bol1
    if ~bol2
    switch flg1_2
        case 4
            dYin=dUr;
        otherwise
            dYin=dYr;
    end
    else
        dYin=[];
    end
    varin_2={A,I,dA,dI,Yin,dYin,flg1,flg2,var1,var2};
    switch flg1_1
        case {1 7}
            if flg1_2==0
                varin_2{7}=1;  [dYrA4, YdA4 ]=BA4(varin_2{:});
            end
    end
    switch flg1_1
        case 2
            [dYrA4, dYrLMA4,YdA4, YdLMA4]=BA4(varin_2{:});
        case {3 5 7}
            if flg1_2>0 || flg1_1==5
                [dYrA4, dCrA4, C, YdA4,dCA4, dC]=BA4(varin_2{:});
            end
        case 4
```

```
            [dYrA4, C, YdA4, dC]=BA4(varin_2{:});
        case 6
            [dYrA4, dYrLMA4, dCrA4, dCrM4, dCrLMA4,  C, ...
              YdA4,  YdLMA4,  dCA4,  dCM4,  dCLMA4, dC]=BA4(varin_2{:});
    end
    if bol2
        YA4=dYrA4;
        switch flg1_1
            case {2 6}
                YLMA4=dYrLMA4;
                switch flg1_1
                    case 6
                        CM4=dCrM4; CLMA4=dCrLMA4;
                end
        end
        switch flg1_1
            case {3 5 6 7}
                if flg1_2>0 || flg1_1==5
                    CA4=dCrA4;
                end
        end
    else
        switch flg1_1
            case {3 4 5 6 7}
                if flg1_2>0 || flg1_1==5
                    dCr=dC(:,:,1); dC=dC(:,:,2:end);
                end
        end
    end
end
% --------------------------------------------------------------------------
Ah=(W*A-A*W) + v*(Y*A+A*Y-2*YA4);
if bol1
for m=1:s3
    dAh(:,:,m)=(W*dA(:,:,m)-dA(:,:,m)*W)+...
                v*(Y*dA(:,:,m)+dA(:,:,m)*Y)-2*v*YdA4(:,:,m);
end
if ~bol2
    dAhr=(dWr*A-A*dWr)+ v*(dYr*A+A*dYr-2*dYrA4) ;
end
end
%
switch flg1_2
case 0
    Ad=2*y*CI*(I-a*A);
    if bol1
        for m=1:s3
            dAd(:,:,m)=-2*CI*y*a*dA(:,:,m);
        end
        if ~bol2
            dAdr=2*CI*dyr*(I-a*A);
        end
    end
end
```

```
% ------------------------------------------------------------------------
switch flg1_1
case 1 % SRF
    FT=Ah+Ad;  F=kp*FT;
    if bol1
        dFT=dAh+dAd;
        if ~bol2
            dFTr=dAhr+dAdr; dFT=cat(3,dFTr,dFT);
        end
        dF=kp*dFT;
    end
case 2 % RSC
    FT=Ah+Ad;  Arsc=2*v*YLMA4+Ad; F=FT-(1-kp)*Arsc;
    if bol1
        dFT =dAh +dAd ; dArsc =2*v*YdLMA4 +dAd ; dF =dFT -(1-kp)*dArsc ;
        if ~bol2
            dFTr=dAhr+dAdr; dArscr=2*v*dYrLMA4+dAdr; dFTr=dFTr-(1-kp)*dArscr;
            dF=cat(3,dFTr,dF);
        end
    end
case {3 7} % ARD , X-RPR
    if flg1_2>0
    Aard=y*(2*C-2*trace(C)*A-5*(C*A+A*C)+10*CA4);
    if bol1
        if ~bol2
            dAardr=dyr/y*Aard+y*(2*dCr-2*trace(dCr)*A-5*(dCr*A+A*dCr)
+10*dCrA4);
        end
        %
        for m=1:s3

dAard(:,:,m)=y*(2*dC(:,:,m)-2*(trace(dC(:,:,m))*A+trace(C)*dA(:,:,m))+...
                -5*(dC(:,:,m)*A+C*dA(:,:,m)+dA(:,:,m)*C+A*dC(:,:,m))
+10*dCA4(:,:,m));
        end
    end
    end
    switch flg1_1
        case 3 % ARD
            F=Ah+Aard;
            if bol1
                dF=dAh+dAard;
                if ~bol2
                    dFr=dAhr+dAardr;  dF=cat(3,dFr,dF);
                end
            end
        case 7 % X-RPR
            nk=[1 2 3];
            switch flg1_2
                case 0    % IRD
                    F=Ah+Ad;
                    if bol1
                        dF=dAh+dAd; if ~bol2, dFr=dAhr+dAdr; end
                    end
```

```
                otherwise % ARD
                    F=Ah+Aard;
                    if bol1
                        dF=dAh+dAard; if ~bol2, dFr=dAhr+dAardr; end
                    end
            end
            % RPR Correction
            if bol1 %
                [~, Q, ~, dQ]=deigvdA(A,I,dA,dI,var2{:});
                eigv=diag(Q\F*Q,0);
                for m=1:s3
                    dv(:,m)=diag(-(Q\dQ(:,:,m)/Q)*F*Q+...
                        Q\dF(:,:,m)*Q+Q\F*dQ(:,:,m),0);
                end
            else
                [~, Q]=deigvdA(A,I,[],[],var2{:}); eigv=diag(Q\F*Q,0);
            end
            eigvt=transpose(eigv);
            bol=1-I; V_iok=ak*(eigv-bk*(eigv.^2+2*prod(eigvt.*bol+I,2)));
            V_iok=diag(V_iok,0); A_iok=-Q*V_iok/Q;
            F=F+A_iok;
            if bol1
                for m=1:3
                    jk=nk(nk~=m); kj=6-m-jk;
                    dV_iok(m,:)=ak*(dv(m,:)-2*bk*(eigv(m).*dv(m,:)+...
                        sum(eigv(jk).*dv(kj,:),1)));
                end
                dV_iok=permute(dV_iok,[1 3 2]).*repmat(I,1,1,s3);
                for m=1:s3
                    dA_iok(:,:,m)=-(dQ(:,:,m)*V_iok/Q+Q*dV_iok(:,:,m)/Q-...
                        Q*V_iok*(Q\dQ(:,:,m)/Q));
                end
                dF=dF+dA_iok; if ~bol2, dF=cat(3,dFr,dF); end
            end
    end
case 4 % MRD
    Amrd=2*y*(C-trace(C)*A); F=Ah+Amrd;
    if bol1
        for m=1:s3
            dAmrd(:,:,m)=2*y*(dC(:,:,m)-trace(dC(:,:,m))*A-
trace(C)*dA(:,:,m));
        end
        dF=dAh+dAmrd;
        if ~bol2
            dAmrdr=2*(dyr*(C-trace(C)*A)+y*(dCr-trace(dCr)*A));
            dFr=dAhr+dAmrdr; dF=cat(3,dFr,dF);
        end
    end
case 5 % NEM
    AA4=CA4; Anem=y*U0*(A*A-AA4); F=Ah+.5*Ad+Anem;
    if bol1
        dAA4=dCA4;
        for m=1:s3
            dAnem(:,:,m)=y*U0*(dA(:,:,m)*A+A*dA(:,:,m)-dAA4(:,:,m));
```

```
        end
        dF=dAh+.5*dAd+dAnem;
        if ~bol2
            dAnemr=dyr/y*Anem; dFr=dAhr+.5*dAdr+dAnemr; dF=cat(3,dFr,dF);
        end
    end
case 6 % ARD-RSC
    Arsc=-2*v*YLMA4;
    Aard_rsc=y*(2*(C-(1-kp)*CM4)-2*kp*trace(C)*A-5*(C*A+A*C)+...
         10*(CA4+(1-kp)*CLMA4));
    F=Ah+(1-kp)*Arsc+Aard_rsc;
    if bol1
        for m=1:s3
            dArsc(:,:,m)=-2*v*YdLMA4(:,:,m);
            dAard_rsc(:,:,m)=y*(...
                2*(dC(:,:,m)-(1-kp)*dCM4(:,:,m))-...
                2*kp*(trace(dC(:,:,m))*A+trace(C)*dA(:,:,m))-...
                5*(dC(:,:,m)*A+C*dA(:,:,m)+dA(:,:,m)*C+A*dC(:,:,m))+...
                10*(dCA4(:,:,m)+(1-kp)*dCLMA4(:,:,m)));
        end
        dF=dAh+(1-kp)*dArsc+dAard_rsc;
        if ~bol2
            dArscr=-2*v*dYrLMA4;
            dAard_rscr=dyr/y*Aard_rsc+y*(2*(dCr-(1-kp)*dCrM4)-...
                2*kp*trace(dCr)*A-5*(dCr*A+A*dCr)+10*(dCrA4+(1-kp)*dCrLMA4));
            dFr=dAhr+(1-kp)*dArscr+dAard_rscr; dF=cat(3,dFr,dF);
        end
    end
end
%
end

function varargout=BA4(A,I,dA ,dI ,Y, dY,flg1,flg2,var1,var2)
%
flg1_1=floor(flg1); flg1_2=round(10*(flg1-flg1_1));
outflg=[1 2 3 2 3 6 3]; nout=outflg(flg1_1);
%
YA4=zeros(3);
if nargout>nout
    s3=size(dA,3); YdA4=zeros(3,3,s3);
    [A4, dA4 ]=aijkl(A,I,dA,dI,flg2,var2{:});
else
     A4       =aijkl(A,I,[],[],flg2,var2{:});
end
% Initialization
switch flg1_1
    case {2 6}
        if nargout>nout
             YdLMA4=zeros(3,3,s3);
            [L4,M4, dL4, dM4]=LM(A,I,dA,dI,var2{:});
            [MA4,dMA4]=fun_MA(M4,A4, dM4,dA4);
        else
             YLMA4=zeros(3,3);
            [L4,M4]=LM(A,I,[],[],var2{:}); MA4=fun_MA(M4,A4);
```

```
        end
        switch flg1_1
            case 6
                CM4=zeros(3,3); CLMA4=zeros(3,3);
                if nargout>nout
                    dCM4=zeros(3,3,s3); dCLMA4=zeros(3,3,s3);
                end
        end
end
switch flg1_1
    case {3 5 6 7}
        CA4=zeros(3,3); if nargout>nout, dCA4=zeros(3,3,s3); end
end
%
switch flg1_1
    case {3 4 6 7} % Spatial Diffusion Tensor
        if nargout>nout
            [C, dC]=Cij(A,I,dA,dI,Y,dY,flg1_2,var1,var2);
            if ~isempty(dY)
                switch flg1_2
                    case {1 2 4}
                        dCr=dC(:,:,1); dC=dC(:,:,2:end);
                end
            end
        else
            C          =Cij(A,I,[],[],Y,[],flg1_2,var1,var2);
        end
        switch flg1_2
            case 4
                L=Y; Y=1/2*(L+L');
                if nargout>nout
                    if ~isempty(dY)
                        dL=dY; dY=1/2*(dL+dL');
                    end
                end
        end
    case 5
        C=A; dC=dA; if ~isempty(dY), dCr=zeros(3); end
end
% -------------------------------------------------------------------------
if ~isempty(dY), Yr=dY; else,  Yr=Y ; end
%
switch flg1_1
    case {3 5 6 7}
        if ~isempty(dY), Cr=dCr; else,Cr=C  ;end
end
%
for i=1:3
    for j=1:3
        YA4(i,j)=0  ; if nargout>nout, YdA4(i,j,1:s3)=0;   end
        switch flg1_1
            case {2 6}
                YLMA4(i,j)=0 ;
                if nargout>nout
```

```
                YdLMA4(i,j,1:s3)=0;
            end
            switch flg1_1
                case 6
                    CLMA4(i,j)=0; CM4(i,j)=0      ;
                    if nargout>nout
                        dCLMA4(i,j,1:s3)=0; dCM4(i,j,1:s3)=0  ;
                    end
            end
    end
    %
    switch flg1_1
        case {3 5 6 7}
            CA4(i,j)=0 ;  if nargout>nout, dCA4(i,j,1:s3)=0  ; end
    end
    %
    for k=1:3
        for l=1:3
            %
            YA4(i,j) = YA4(i,j)+ Yr(k,l)*A4(i,j,k,l);
            if nargout>nout
                for n=1:s3
                    YdA4(i,j,n) = YdA4(i,j,n)+Y(k,l)*dA4(i,j,k,l,n);
                end
            end
            %
            switch flg1_1
                case {2 6}
                    YLMA4(i,j) = YLMA4(i,j)+...
                        Yr(k,l)*(L4(i,j,k,l)-MA4(i,j,k,l));
                    if nargout>nout
                        for n=1:s3
                            YdLMA4(i,j,n) = YdLMA4(i,j,n)+...
                                Y(k,l)*(dL4(i,j,k,l,n)-dMA4(i,j,k,l,n));
                        end
                    end
                    %
                    switch flg1_1
                        case 6
                            CLMA4(i,j) = CLMA4(i,j)+...
                                Cr(k,l)*(L4(i,j,k,l)-MA4(i,j,k,l));
                            CM4(i,j) = CM4(i,j)+Cr(k,l)*M4(i,j,k,l);
                            if nargout>nout
                                for n=1:s3
                                    dCLMA4(i,j,n) = dCLMA4(i,j,n)+...
                                        C(k,l  )*(dL4(i,j,k,l,n)-
dMA4(i,j,k,l,n))+...

                                        dC(k,l,n)*( L4(i,j,k,l  )-
MA4(i,j,k,l  ));
                                    dCM4(i,j,n) = dCM4(i,j,n)+...
                                        C(k,l)*dM4(i,j,k,l,n)

+dC(k,l,n)*M4(i,j,k,l);
                                end
                            end
```

```
                        end
                end
                %
                switch flg1_1
                    case {3 5 6 7}
                        CA4(i,j) = CA4(i,j)+Cr(k,l)*A4(i,j,k,l);
                        if nargout>nout
                            for n=1:s3
                                dCA4(i,j,n) = dCA4(i,j,n)+...
                                    C(k,l)*dA4(i,j,k,l,n)
+dC(k,l,n)*A4(i,j,k,l);
                            end
                        end
                end
            end
        end
    end
end
%
switch flg1_1
    case 1
        varargout={YA4};
        if nargout>nout
            varargout=[varargout {YdA4}];
        end
    case 2
        varargout={YA4, YLMA4};
        if nargout>nout
            varargout=[varargout {YdA4 YdLMA4}];
        end
    case {3 5 7}
        % AA4=CA4 - 5
        varargout={YA4, CA4, C};
        if nargout>nout
            if ~isempty(dY), dC=cat(3,Cr,dC); end
            varargout=[varargout {YdA4, dCA4, dC}];
        end
    case 4
        varargout={YA4, C};
        if nargout>nout
            if ~isempty(dY), dC=cat(3,Cr,dC); end
            varargout=[varargout {YdA4, dC}];
        end
    case 6
        varargout={YA4, YLMA4, CA4, CM4, CLMA4, C};
        if nargout>nout
            if ~isempty(dY), dC=cat(3,Cr,dC); end
            varargout=[varargout ...
                {YdA4 YdLMA4, dCA4, dCM4, dCLMA4, dC}];
        end
end
end
```

# Spatial Diffusion Tensor

```
function [C, dC]=Cij(A,I,dA,dI,Y,dY,flg,var1,var2)
if nargout>1
    s3=size(dA,3); dC=zeros(3,3,s3);
end
switch flg
    case {1 2}
        y=sqrt(2*sum(Y.*Y','all'));
        if nargout>1, if ~isempty(dY), dy=2/y*sum(Y.*dY','all');end; end
        switch flg
            case 1 % ARD
                b=var1{1}; C=b(1)*I+b(2)*A+b(3)*(A*A')+b(4)*Y/y+b(5)*(Y*Y')/
y^2;
                if nargout>1
                    dC =b(2)*dA;
                    for m=1:s3
                        dC(:,:,m)=dC(:,:,m)+b(3)*(A*dA(:,:,m)'+dA(:,:,m)*A');
                    end
                    if ~isempty(dY)
                    dCr=b(4)*(y*dY-Y*dy)/y^2+...
                        2*b(5)*(y*(Y*dY'+dY*Y')-2*(Y*Y')*dy)/y^3;
                    end
                end
            case 2 % iARD
                [CI, CM]=var1{1:2};   C=CI*(I-4*CM*(Y*Y')/y^2);
                if nargout>1
                    if ~isempty(dY)
                    dCr=-4*CI*CM*(y*(Y*dY'+dY*Y')-2*(Y*Y')*dy)/y^3;
                    end
                end
        end
    case 3 % pARD
        [CI,D]=var1{1:2};
        if nargout>1
            [~, Q,~, dQ]=deigvdA(A,I,dA,dI,var2{:});
        else
            [~, Q]=deigvdA(A,I,[],[],var2{:});
        end
        C=CI*Q*D/Q;
        if nargout>1
            for m=1:s3
                dC(:,:,m)=CI*(dQ(:,:,m)*D/Q-Q*D*(Q\dQ(:,:,m)/Q));
            end
            if ~isempty(dY), dCr=zeros(3); end
        end
    case 4 % iARD
        [CI,CM]=var1{1:2};
        L=Y; lb=sum(L.*L,'all'); Lb=(L*L')/lb; C=CI*(I-CM*Lb);
        if nargout>1
            if ~isempty(dY)
            dL=dY;    dlbr=sum(2*L.*dL,'all');
            dLbr=(lb*(dL*L'+L*dL')-(L*L')*dlbr)/lb^2;
```

```
            dCr=-CI*CM*dLbr;
            end
        end
    case 5 % WPT
        [CI, w]=var1{1:2}; C=CI*((1-w)*I+w*(A*A'));
        if nargout>1
            for m=1:s3
                dC(:,:,m)=CI*w*(A*dA(:,:,m)'+dA(:,:,m)*A');
            end
            if ~isempty(dY), dCr=zeros(3); end
        end
    case 6 % Dz
        [CI, Dz, n]=var1{1:3}; C=CI*(I-(1-Dz)*(n'*n));
        if nargout>1
            if ~isempty(dY), dCr=zeros(3); end
        end
end
if nargout>1, if ~isempty(dY), dC=cat(3,dCr,dC); end, end
end

function [MA,dMA]=fun_MA(M,A,dM,dA)
MA=zeros(3,3,3,3);
if nargout>1, sz=size(dA); s3=sz(end); dMA=zeros(3,3,3,3,s3); end
for i=1:3
    for j=1:3
        for k=1:3
            for l=1:3
                MA(i,j,k,l)=0;
                if nargout>1, dMA(i,j,k,l,:)=0;end
                for m=1:3
                    for n=1:3
                        MA(i,j,k,l)=MA(i,j,k,l)+M(i,j,m,n)*A(m,n,k,l);
                        if nargout>1
                        dMA(i,j,k,l,:)=dMA(i,j,k,l,:)+...
                            dM(i,j,m,n,:).*A(m,n,k,l)
+M(i,j,m,n).*dA(m,n,k,l,:);
                        end
                    end
                end
            end
        end
    end
end
end

function [L,M, dL, dM]=LM(A,I,dA,dI,varargin)
%
M=zeros(3,3,3,3)  ; L=zeros(3,3,3,3);
if nargout>2
    sz=size(dA); s3=sz(end); dM=zeros(3,3,3,3,s3);dL=zeros(3,3,3,3,s3);
    [egv,Q , degv, dQ]=deigvdA(A,I,dA,dI,varargin{:});
else
    [egv,Q]=deigvdA(A,I,[],[],varargin{:});
end
```

```
for m=1:3
    Mm =zeros(3,3,3,3)  ; Lm=zeros(3,3,3,3);
    if nargout>2
        dMm=zeros(3,3,3,3,s3);dLm=zeros(3,3,3,3,s3);
    end
    for i=1:3
        for j=1:3
            for k=1:3
                for l=1:3
                    f=[i j k l];
                    Q4=1;
                    if nargout>2, dQ4=zeros(1,1,s3); end
                    for r =1:4
                        Q4=Q4*Q(f(r),m);
                        if nargout>2
                            dQr=dQ(f(r),m,:);
                            for s=1:4
                                if r~=s
                                    dQr=dQr.*Q(f(s),m);
                                end
                            end
                            dQ4=dQ4+dQr;
                        end
                    end
                    Mm(i,j,k,l)= Q4; Lm(i,j,k,l)=egv(m)*Q4;
                    if nargout>2
                        dMm(i,j,k,l,:)=dQ4;
                        dLm(i,j,k,l,:)=egv(m)*dQ4+permute(degv(m,:),[1 3
2]).*Q4;
                    end
                end
            end
        end
    end
    M=M + Mm ; L=L + Lm;
    if nargout>2, dM=dM + dMm ; dL=dL + dLm; end
end
end

function [F2, dF2]=aijkl(A,I,dA,dI,flg,varargin)
flgs=split(num2str(flg),'.');flg1=str2double(flgs{1});
if length(flgs)>1, flg2=str2double(flgs{2});end
%
if nargout>1
    s3=size(dA,3);
end
%
switch sign(flg)
    case {0 1}
        a=varargin{1};
    case -1
        if nargout>1
            [eigv,Q,deigv, dQ]=deigvdA(A,I,dA,dI,varargin{:});
        else
```

```
            [eigv, Q]=deigvdA(A,I,[],[],varargin{:});
        end

end
%
switch sign(flg)
    case  {0, 1}
        switch flg1
            case {0 1}
                if nargout>1
                    [F2_1, dF2_1]=gencls(A,dA,3); % quad
                    [F2_2, dF2_2]=gencls(A,dA,2); % lin
                else
                    F2_1=gencls(A,[],3); % quad
                    F2_2=gencls(A,[],2); % lin
                end
                switch flg1
                    case 0
                        f=flg; F2=f*F2_1+(1-f)*F2_2;
                        if nargout>1, dF2=f*dF2_1+(1-f)*dF2_2; end
                    case 1
                        switch flg2
                            case 1 % 1.1
                                a1=a/(a-1); b1=1/(a-1);
                                ysq=sum(A.*A','all'); f=a1*ysq-b1;
                                if nargout>1
                                    df=a1*sum(A.*permute(dA,[2 1 3])+dA.*A',
[1 2]);

                                end
                            case 2 % 1.2
                                detj=f_detj(A); f=1-a^a*detj;
                                if nargout>1
                                    [~,ddetj]=f_detj(A,dA); df=-a^a*ddetj;
                                end
                        end
                        F2=f*F2_1+(1-f)*F2_2;
                    if nargout>1
                        dF2=f*dF2_1+(1-f)*dF2_2;
                        for m=1:s3
                            dF2(:,:,:,:,m)=dF2(:,:,:,:,m)+df(m)*(F2_1-F2_2);
                        end
                    end
                end
                %
            case 2 % General permutation closure
                if nargout>1
                    [F2, dF2]=gencls(A,dA,flg2);
                else
                    F2=gencls(A,[],flg2);
                end
            case 3 % IBOF
                if nargout>1
                    [F2, dF2]=A4_IBOF(A,I,dA,dI,varargin{:});
                else
```

```
                    F2=A4_IBOF(A,I,[],[],varargin{:});
                end
            otherwise %

        end
        %
    case -1
        if nargout>1
            [av,dav]=binom(eigv, -flg1); av=av'; dav=dav';
            dav=dav(:,1)*deigv(1,:)+dav(:,2)*deigv(2,:);
        else
            av=binom(eigv, -flg1); av=av';
        end
        %
        switch flg1
            case -3 % Cubic
                switch flg2
                    case {2 3} % Rational ellipsoid
                        if nargout>1
                            [av2,dav2]=binom(eigv, -flg1+1); av2=av2';
dav2=dav2';
                            dav2=dav2(:,1)*deigv(1,:)+dav2(:,2)*deigv(2,:);
                            av=[av;av2]; dav=[dav;dav2];
                        else
                            av2=binom(eigv, -flg1+1); av2=av2'; av=[av;av2];
                        end
                end
        end
%
C=C_EBOF(flg1,flg2); D=[0 1 1;1 0 1;1 1 0];
%
if any(flg==[-3.2,-3.3])
    ut=1:10; vt=11:16;
    u_ex=C(:,ut)*av(ut); v_ex=C(:,vt)*av(vt); A_ex=u_ex./v_ex;
else
    A_ex=C*av;
end
A_dv=D\(eigv-A_ex); A=diag([A_ex;A_dv],0);
for m=1:3, for n=1:3, if m~=n, k=9-m-n; A(m,n)=A(k,k); end , end, end
%
if nargout>1
    if any(flg==[-3.2,-3.3])
        du_ex=C(:,ut)*dav(ut,:); dv_ex=C(:,vt)*dav(vt,:);
        dA_ex=(v_ex.*du_ex-u_ex.*dv_ex)./v_ex.^2;
    else
        dA_ex=C*dav;
    end
    dA_dv=D\(deigv-dA_ex); dA_d=[dA_ex;dA_dv];
    for k=1:s3
        for i=1:6
            for j=1:6
                if i==j
                    dA6(i,j,k)=dA_d(i,k);
                else
```

```
                    if i<=3 && j<=3
                        dA6(i,j,k)=dA_d(9-i-j,k);
                    end
                end
            end
        end
    end
end
%
fnc=@(m,n) (m==n)*m+(m~=n)*(9-m-n);
for i=1:3
    for j=1:3
        for k=1:3
            for l=1:3
                p1=fnc(i,j); p2=fnc(k,l);
                F(i,j,k,l)=A(p1,p2);
                if nargout>1, dF(i,j,k,l,:)=dA6(p1,p2,:); end
            end
        end
    end
end
%
for i=1:3
    for j=1:3
        for k=1:3
            for l=1:3
                ijkl=[i j k l]; val=0;
                if nargout>1, dval=zeros(1,s3); end
                for r=1:3
                    for s=1:3
                        for m=1:3
                            for n=1:3
                                rsmn=[r s m n]; Q4=1;
                                if nargout>1, dQ4=zeros(1,1,s3); end
                                for p=1:4
                                    Q4=Q4*Q(ijkl(p),rsmn(p));
                                    if nargout>1
                                        dQ4p=dQ(ijkl(p),rsmn(p),:);
                                        for q=1:4
                                            if p~=q
                                                dQ4p=dQ4p*Q(ijkl(q),rsmn(q));
                                            end
                                        end
                                        dQ4=dQ4+dQ4p;
                                    end
                                end
                                val=val+Q4*F(r,s,m,n);
                                if nargout>1
                                    for h=1:s3
                                        dval(h)=dval(h)+Q4*dF(r,s,m,n,h)
+dQ4(h)*F(r,s,m,n);
                                    end
                                end
                            end
```

```
                        end
                    end
                end
                F2(i,j,k,l)=val;  if nargout>1, dF2(i,j,k,l,:)=dval;end
            end
        end
    end
end
end
end

function [A4, dA4]=A4_IBOF(A,I,dA,dI,varargin)
A2v=v2M(A,2); d=I; s3=size(dA,3);
%
nA4=[1 1 1 1;1 1 2 2;1 1 2 3 ;...
     1 1 1 3;1 1 1 2;2 2 2 2 ;...
     2 2 2 3;2 2 1 3;2 2 1 2];

n4=size(nA4); A4v=zeros(n4(1),1); nn1=factorial(n4(2));
if nargout>1
    [bta, dbta]=bta_fnc(A,I,dA,dI,varargin{:}); dA4v=zeros(n4(1),s3);
else
    bta=bta_fnc(A,I,[],[],varargin{:});
end

for i=1:n4(1)
    ijkl=nA4(i,:);    A4v(i)=0;
    for j=1:3
        np=[perms(1:n4(2)) n4(2)+(1:j-1).*ones(nn1,1)];
        %
        nn=size(np); nm=nn(2)-j+1;
        np2=zeros(nn(1),2*(j+1)); np2(:,[1 end-nm+2:end])=np(:,1:nm);
        for jj=1:j-1
            jq1=2*j-4*jj+4;  jq2=nm-2*jj; j3=1/2*(j-1)*(j-2)*(jj-1);
            np2(:,jq1  :jq1+1)=repmat(np(:,nm+j-jj),1,2);
            np2(:,jq1+j3-2:jq1-1)=np(:,jq2+j3:jq2+1);
        end
        %
        ijklx=catj(ijkl ,j-1); nq=size(ijklx,1);
        %
        f0=0; df0=0; M2=repmat(A,1,1,j+1);
        for k=1:4-j
            f1=0; df1=0; kj=4*j-j*(j+1)/2+1-k;
            for m=1:nn(1)
                mn4=ijklx(:,np2(m,:));
                for r=1:nq
                    f2=1;
                    for n=1:j+1
                        p=2*n-1;
                        f2=f2*M2(mn4(r,p),mn4(r,p+1),n);
                    end
                    f1=f1+f2;
                    %
                    if nargout>1
```

```
                        df2=0;
                        for n=1:j-k+2
                            p=2*n-1;
                            df3(1,:)=dA(mn4(r,2*j-p+2),mn4(r,2*j-p+3),:);
                            for s=1:j+1
                                if s~=n
                                    q=2*s-1;
                                    df3=df3*M2(mn4(r,2*j-q+2),mn4(r,2*j-
q+3),j+2-s);
                                end
                            end
                            df2=df2+df3;
                        end
                        df1=df1+df2;
                    end
                end
            end
            f1= f1/nn(1) ;  f0=f0+bta(kj)*f1;  M2(:,:,k)=d;
            if nargout>1
                df1=df1/nn(1); df0=df0+bta(kj)*df1+dbta(kj,:)*f1;
            end
        end
        A4v(i)=A4v(i)+f0; if nargout>1, dA4v(i,:)=dA4v(i,:)+df0; end
    end
end
if nargout>1
    [A4,dA4]=A4v2M(A2v,A4v,dA4v,dA);
else
    A4=A4v2M(A2v,A4v);
end
end

function [b,db]=bta_fnc(A,M,dA,dM,varargin)
s3=size(dA,3);
if nargout>1
    [I,dI]=INV(A,M,dA,dM,varargin{:});
else
    I=INV(A,M,[],[],varargin{:});
end
% I=[trace(A2);.5*(trace(A2)^2-trace(A2^2));det(A2)];
pv=[1,2,4,3,6,5,9,7,8,10,14,11,12,13,15,20,16,17,18,19,21];
b0=zeros(3,1); db0=zeros(3,s3);
for k=1:3
    b0(k)=C_IBOF(1,k);
    for i=1:5
        for j=0:i
            ij=j+.5*(i+1)*i+1; fac=C_IBOF(pv(ij),k);
            b0(k)=b0(k)+fac*I(2)^(i-j)*I(3)^j;
            if nargout>1
            db0(k,:)=db0(k,:)+fac*(i-j)*I(2)^(i-j-1)*I(3)^j*dI(2,:)+...
                                  fac*j*I(2)^(i-j)*I(3)^(j-1)*dI(3,:);
            end
        end
    end
```

```
end

%
b=zeros(6,1); b([3 4 6])=b0;
b(1)=3/5*(-1/7+1/5*b(3)*(1/7+4/7*I(2)+8/3*I(3))-
b(4)*(1/5-8/15*I(2)-14/15*I(3))+...

-b(6)*(1/35-24/105*I(3)-4/35*I(2)+16/15*I(2)*I(3)+8/35*I(2)^2));
b(2)=6/7*(1-1/5*b(3)*(1+4*I(2))+7/5*b(4)*(1/6-I(2))-
b(6)*(-1/5+2/3*I(3)+4/5*I(2)-8/5*I(2)^2));
b(5)=-4/5*b(3)-7/5*b(4)-6/5*b(6)*(1-4/3*I(2));
%
if nargout>1
db=zeros(6,s3); db([3 4 6],:)=db0;
db(1,:)=3/5*(1/5*db(3,:)*(1/7+4/7*I(2)+8/3*I(3))-
db(4,:)*(1/5-8/15*I(2)-14/15*I(3))+...

-db(6,:)*(1/35-24/105*I(3)-4/35*I(2)+16/15*I(2)*I(3)+8/35*I(2)^2)+...
             1/5*b(3)*(4/7*dI(2,:)+8/3*dI(3,:))-
b(4)*(-8/15*dI(2,:)-14/15*dI(3,:))+...

-b(6)*(-24/105*dI(3,:)-4/35*dI(2,:)+16/15*(dI(2,:)*I(3)+I(2)*dI(3,:))+...
                 16/35*I(2)*dI(2,:)));
db(2,:)=6/7*(-1/5*db(3,:)*(1+4*I(2))+7/5*db(4,:)*(1/6-I(2))-
db(6,:)*(-1/5+2/3*I(3)+4/5*I(2)-8/5*I(2)^2)+...
             -4/5*b(3)*dI(2,:)-7/5*b(4)*dI(2,:)-
b(6)*(2/3*dI(3,:)+4/5*dI(2,:)-16/5*I(2)*dI(2,:)));
db(5,:)=-4/5*db(3,:)-7/5*db(4,:)-6/5*db(6,:)*(1-4/3*I(2))+8/5*b(6)*dI(2,:);
end
end

function [f,df]=INV(A,I,dA,dI,varargin)
if nargout>1
    [v, ~, dv]=deigvdA(A,I,dA,dI,varargin{:});
else
    v=deigvdA(A,I,[],[],varargin{:});
end
f=zeros(3,1);
for i=1:3
    ij=nchoosek(1:3,i); nij=size(ij);
    f0=0; df0=0;
    for j=1:nij(1)
        f1=1; df1=0;
        for k=1:nij(2)
            f1=f1*v(ij(j,k));
            %
            if nargout>1
                df2=dv(ij(j,k),:);
                for m=1:nij(2)
                    if m~=k
                        df2=df2*v(ij(j,m));
                    end
                end
                df1=df1+df2;
```

```
            end
            %
        end
        f0=f0+f1; if nargout>1,  df0=df0+df1; end
    end
    f(i)=f0; if nargout>1, df(i,:)=df0; end
end
end

function [f,df]=A4v2M(A2v,A4v,varargin)
% A2v=['A11','A12','A13','A22','A23']
%
A4v=['A1111','A1122','A1123','A1113','A1112','A2222','A2223','A2213','A2212']
d=eye(3); f1=@(i,j) (9-i-j).*(1-d(i+3*(j-1)))+d(i+3*(j-1)).*i;
f2=@(m) [m, m]+(m>3)*([-m, 9-2*m]+[1 -1].*(m^2-11*m+32)/2);
Mij=@(M,i,j) logical(prod((M==i)+(M==j),2));
%
jj=0; Axx=[];B2=[];
for i =1:3
    for j=i:3
        for k=j:3
            for l=k:3
                ijkl=[i j k l]; jj=jj+1;  Axx=[Axx; f1(i,j) f1(k,l)];
                B1=unique(perms(ijkl),'Rows'); njj=jj*ones(size(B1,1),1);
                B2=[B2;njj sum((3.^(0:3)).*(B1-1),2)+1];
            end
        end
    end
end
%
k=0; for i=1:3,for j=i:6, k=k+1; D1(i,j)=k; end, end
%
Axx(Mij(Axx,4,6),:)=[2 5]; Axx=sort(Axx,2);
[~,ndk]=sortrows(Axx,[1 2]); D2=(1:15)'; D3=D2(ndk,1);
%
jj=0;
for i=1:2
    for j=i:6
        if j~=3
            jj=jj+1;
            n_inp(jj,1)=D3(D1(i, j),1);
        end
    end
end
D2(n_inp,2)=A4v; A2=v2M(A2v,1);
if nargin>2
    [dA4v,dA]=varargin{1:2}; dD2(n_inp,:)=dA4v; s3=size(dA,3);
end
%
for k=1:2
    for i=1:3
        m=i+3*(k-1); r=f2(m); s=[(2-k) (k-1)]*[m 3;3 m];
        val=A2(r(1),r(2));
        if nargin>2, dval(1,:,1)=dA(r(1),r(2),:); end
```

```
        for j=1:2
            p=sort([j m],2);
            val=val-D2(D3(D1(p(1),p(2)),1),2);
            if nargin>2
                dval=dval-dD2(D3(D1(p(1),p(2)),1),:);
            end
        end
        D2(D3(D1(s(1),s(2)),1),2)=val;
        if nargin>2
            dD2(D3(D1(s(1),s(2)),1),:)=dval;
        end
    end
end
f(B2(:,2))= D2(B2(:,1),2); f=reshape(f,3,3,3,3);
%
if nargout>1
    df(B2(:,2),:)= dD2(B2(:,1),:); df=reshape(df,3,3,3,3,s3);
end
%
end

function [v,Q, dv,dQ]=deigvdA(K,M,dK,dM,k,flg,varargin)
switch flg
    case {1 2 3}
        % syms v; v=vpasolve(det(K-v*M)==0,v); v=double(sort(v,'descend'));
        % det(K)-sum(M.*adjoint(K)','all')*v+sum(K.*adjoint(M)','all')*v^2-
det(M)*v^3=0;
        v=roots([-det(M),sum(K.*adjoint(M)','all'),-
sum(M.*adjoint(K)','all'),det(K)]);
        v=sort(v,'descend');
        %
        [s1,s2]=size(K); if nargout>2,s3=size(dK,3); end
        for i=1:s2
            S=K-v(i)*M;
            Qi=rand(s1,1); err=1;
            while err>1e-8
                Ci=Qi'*M*Qi;
                switch flg
                    case 1
                        Rk=Ci-1; dSk=2*Qi'*M;
                    case 2
                        [p,a]=varargin{1:2};Z=zeros(1,s2); Z(p)=1;
                        Rk=Z*Qi-a; dSk=Z;
                    case 3
                        b=varargin{1};
                        Rk=sqrt(Qi'*Qi)-b; dSk=Qi';
                end
                %
                R=S*Qi; R(k)=Rk; J=S; J(k,:)=dSk;
                Qn=Qi-2.*J\R; err=norm(Qn-Qi); Qi=Qn;
            end
            Q(:,i)=Qi;
            %
            if nargout>2
```

```
                for j=1:s3
                    dv(i,j)=1/Ci*(Qi'*(dK(:,:,j)-v(i)*dM(:,:,j))*Qi);
                    F=-(dK(:,:,j)-dv(i,j)*M-v(i)*dM(:,:,j))*Qi;
                    switch flg
                        case 1
                            Fk=-Qi'*dM(:,:,j)*Qi;
                        case {2 3}
                            Fk=0;
                    end
                    dF=F; dF(k)=Fk; P=S; P(k,:)=dSk; dQ(:,i,j)=P\dF;
                end
            end
        end
    case 4
        if nargout>2
            [v, Q, dv, dQ]=Nlsn(K,M,dK,dM,k);
        else
            [v, Q]=Nlsn(K,M,dK,dM,k);
        end
end
end

function [v, Q, dv, dQ]=Nlsn(K,M,dK,dM,k)
[Q, v]=eig(K);   [v, ind]=sort(diag(v,0),'descend');Q=Q(:,ind);
[s1,s2]=size(K);
if nargout>2
    s3=size(dK,3); dv=zeros(s2,s3); dQ=zeros(s1,s2,s3);
    for j=1:s2
        L=v(j); Qj=Q(:,j);
        Cj=Qj'*M*Qj; Sj=K-L*M; Sj(k,:)=0; Sj(:,k)=0; Sj(k,k)=1;
        for i=1:s3
            dKi=dK(:,:,i); dMi=dM(:,:,i);
            dLi=1/Cj*(Qj'*(dKi-L*dMi)*Qj); dv(j,i)=dLi;
            if nargout>3
                Fj=-(dKi-dLi*M-L*dMi)*Qj; Fj(k)=0;
                Vj=Sj\Fj; cj=-Qj'*M*Vj-1/2*Qj'*dMi*Qj;
                dQ(:,j,i)=Vj+cj.*Qj;
            end
        end
    end
end
end

function [detj,ddetj]=f_detj(A,dA)
detj=0; s1=size(A,1);
if nargout> 1, s3=size(dA,3); ddetj=zeros(1,1,s3); end
ijkn=perms(1:s1); [s1, s2]=size(ijkn);
for i=1:s1
    f=1; ijk=ijkn(i,:); eijk=par(ijk);
    for m=1:s2
        f=f*A(ijk(m),m);
        if nargout>1
            df=dA(ijk(m),m,:);
            for n=1:s2
```

```
                if m~=n
                    df=df*A(ijk(n),n);
                end
            end
            ddetj=ddetj+eijk*df;
        end
    end
    detj=detj+eijk*f;
end
end

function [R,J]=dMdA(func, A, del, flg, varargin)
A=A(:); n2=numel(A); nA=size(A); nA=num2cell(nA);
D=zeros(n2,1); tol=eps;
f=@(A) func(A,varargin{:}); R=f(A); n1=numel(R);
%
J=zeros(n1*n2,1);
Fac={
    1,          [-1  1 ; -1  0] ;
    1,          [-1  1 ;  0  1] ;
    2,          [-1  1 ; -1  1] ;
    2,      [1 -4  3 ; -2 -1 0] ;
    2,      [-3 4 -1 ;  0  1 2] ;
   12, [1 -8  8 -1 ; -2 -1 1 2]
                                    };
%
M1=Fac{flg,1}; M2=Fac{flg,2}; m=length(M2);
for i=1:n2
    na=(i-1)*n1+1; nb=na+n1-1;
    D(i,1)=1; J(na:nb,1)=zeros(n1,1); dA=del*D; h=del;
    if abs(A(i))>tol,  dA=dA.*A(:); h=h*A(i);  end
    for j=1:m
        Jj=M2(1,j)*f(A+M2(2,j)*reshape(dA,nA{:}));
        J(na:nb,1)=J(na:nb,1)+Jj(:);
    end
    J(na:nb,1)=J(na:nb,1)/(M1*h);
    D(i,1)=0;
end
nR=size(R); nR=num2cell(nR);J=reshape(J,nR{:},n2);
end

function [f,df]=gencls(A,dA,flg)
d=eye(3); f=zeros(3,3,3,3);
%
if nargout>1
    s3=size(dA,3); [b,db]=fbta(A,dA,flg); df=zeros(3,3,3,3,s3);
else
    b=fbta(A,[],flg);
end
%
for i=1:3
    for j=1:3
        for k=1:3
            for l=1:3
```

```
                f(i,j,k,l)=b(1)*(d(i,j)*d(k,l))+b(2)*(d(i,k)*d(j,l)
+d(i,l)*d(j,k))+...
                    b(3)*(d(i,j)*A(k,l)+A(i,j)*d(k,l))+...
                    b(4)*(A(i,k)*d(j,l)+A(j,l)*d(i,k)+A(i,l)*d(j,k)
+A(j,k)*d(i,l))+...
                    b(5)*(A(i,j)*A(k,l))+b(6)*(A(i,k)*A(j,l)+A(i,l)*A(j,k));
                if nargout>1
                    df(i,j,k,l,:)=db(1,1,:)*(d(i,j)*d(k,l))
+db(1,2,:)*(d(i,k)*d(j,l)+d(i,l)*d(j,k))+...
                        db(1,3,:)*(d(i,j)*A(k,l)+A(i,j)*d(k,l))+...
                        db(1,4,:)*(A(i,k)*d(j,l)+A(j,l)*d(i,k)+A(i,l)*d(j,k)
+A(j,k)*d(i,l))+...
                        db(1,5,:)*(A(i,j)*A(k,l))+db(1,6,:)*(A(i,k)*A(j,l)
+A(i,l)*A(j,k));
                    for q=1:s3
                        df(i,j,k,l,q)=df(i,j,k,l,q)+...
                            b(3)*(d(i,j)*dA(k,l,q)+dA(i,j,q)*d(k,l))+...
                            b(4)*(dA(i,k,q)*d(j,l)+dA(j,l,q)*d(i,k)
+dA(i,l,q)*d(j,k)+dA(j,k,q)*d(i,l))+...
                            b(5)*(dA(i,j,q)*A(k,l)+A(i,j)*dA(k,l,q))+...
                            b(6)*(dA(i,k,q)*A(j,l)+A(i,k)*dA(j,l,q)
+dA(i,l,q)*A(j,k)+A(i,l)*dA(j,k,q));
                    end
                end
                %
                for m=1:3
                    f(i,j,k,l)=f(i,j,k,l)+b(7)*(d(i,j)*A(k,m)*A(m,l)
+A(i,m)*A(m,j)*d(k,l));
                    if nargout>1
                        for q=1:s3
                            df(i,j,k,l,q)=df(i,j,k,l,q)
+db(1,7,q)*(d(i,j)*A(k,m)*A(m,l)+A(i,m)*A(m,j)*d(k,l))+...
                                b(7)*(d(i,j)*dA(k,m,q)*A(m,l)
+d(i,j)*A(k,m)*dA(m,l,q)+...
                                dA(i,m,q)*A(m,j)*d(k,l)
+A(i,m)*dA(m,j,q)*d(k,l));
                        end
                    end
                    for n=1:3
                        f(i,j,k,l)=f(i,j,k,l)
+b(8)*(A(i,m)*A(m,j)*A(k,n)*A(n,l));
                        if nargout>1
                            for q=1:s3
                                df(i,j,k,l,q)=df(i,j,k,l,q)
+db(1,8,q)*(A(i,m)*A(m,j)*A(k,n)*A(n,l))+...
                                    b(8)*(dA(i,m,q)*A(m,j)*A(k,n)*A(n,l)
+A(i,m)*dA(m,j,q)*A(k,n)*A(n,l)+...
                                    A(i,m)*A(m,j)*dA(k,n,q)*A(n,l)
+A(i,m)*A(m,j)*A(k,n)*dA(n,l,q));
                            end
                        end
                    end
                end
                %
```

```
            end
        end
    end
end
end

function [b,db]=fbta(A,dA,flg)
%
ysq=sum(A.*A','all'); a=exp(2*(1-3*ysq)/(1-ysq));
%
bta=[
    1/15        1/15        0           0       0       0       0       0;
   -1/35       -1/35      1/7         1/7       0       0       0       0;
       0           0        0           0       1       0       0       0;
       0           0        0           0       1       1       0  -2/ysq;
       0           0      2/5           0    -1/5     3/5    -2/5       0;
26*a/315    26*a/315  16*a/63     -4*a/21       1       1       0  -2/ysq];
%
b=bta(flg,:);
%
if nargout>1
    s3=size(dA,3); dysq=sum(2*dA.*A',[1 2]); da=(-4*a/(1-ysq)^2)*dysq;
    dbta=zeros(6,8,s3);  dbta([4 6],8,:)=2/ysq^2*dysq.*[1;1];
    dbta(6,1:4,:)=[26/315    26/315   16/63    -4/21].*da; db=dbta(flg,:,:);
end
end

function f=DA(d)
for i=1:3
    for j=1:3
        for m=1:2
            for n=m:3
                k=2*(m-1)+n; p=-.5*(i^2-3*i+2);
                f(i,j,k)=d(i,m)*d(j,n)+(1-d(i,j))*d(i,n)*d(j,m)
+p*d(i,j)*d(m,n);
            end
        end
    end
end
end

function A=v2M(Av,flg)
for m=1:2
    for n=m:3
        k=2*(m-1)+n;
        switch flg
            case 1
                A(m,n)=Av(k); A(n,m)=A(m,n);
            case 2
                A(k,1)=Av(m,n);
        end
    end
end
if flg==1, A(3,3)=1-A(1,1)-A(2,2); end
end
```

```
function A=catj(A,p)
pj=1;
while pj<=p
    n=size(A); B=zeros(3*n(1),n(2)+1);
    for k=1:3
        rk=(k-1)*n(1)+1:k*n(1);
        for j=1:n(1)
            B(rk,1:n(2))=A; B(rk,n(2)+1)=k;
        end
    end
    A=B; pj=pj+1;
end
end

function [t,X]=rk4(fnc,t, X0, varargin)
n=length(t);  X(1,:)=X0';
for j=1:n-1
    dt=t(j+1)-t(j);
    K1=fnc(t(j)     , X(j,:)'         ,varargin{:});
    K2=fnc(t(j)+dt/2, X(j,:)' +dt/2*K1,varargin{:});
    K3=fnc(t(j)+dt/2, X(j,:)' +dt/2*K2,varargin{:});
    K4=fnc(t(j)+dt  , X(j,:)' +dt*K3  ,varargin{:});
    dXj=dt/6*(K1+2*K2+2*K3+K4); X(j+1,:) = X(j,:) + dXj';
end
end

function [av,dav]=binom(v, n)
av=1; if nargout>1, dav=[0;0]; end
for i=1:n
    for j=0:i
        av=[av v(1)^(i-j)*v(2)^j];
        if nargout>1
            dav=[dav,...
                [(i-j)*v(1)^(i-j-1)*v(2)^j      ;
                j*v(1)^(i-j  )*v(2)^(j-1)]];
        end
    end
end
end

function [f,a]=par(a)
n=length(a); k=0; j=1; m=0;
while 1
    aj=a(j);
    if j~=aj
        c=a(aj); a(aj)=aj;  a(j)=c; k=k+1; m=m+abs(a(j)-j);
    end
    m=m*(j<n); j=j*(j<n)+1;
    if (j==n)*(m==0)==1, break, end
end
f=(-1)^k;
end

function f=C_IBOF(i,j)
C=[
```

```
 0.249409081657860E+02  -0.497217790110754E+00   0.234146291570999E+02
-0.435101153160329E+03   0.234980797511405E+02  -0.412048043372534E+03
 0.372389335663877E+04  -0.391044251397838E+03   0.319553200392089E+04
 0.703443657916476E+04   0.153965820593506E+03   0.573259594331015E+04
 0.823995187366106E+06   0.152772950743819E+06  -0.485212803064813E+05
-0.133931929894245E+06  -0.213755248785646E+04  -0.605006113515592E+05
 0.880683515327916E+06  -0.400138947092812E+04  -0.477173740017567E+05
-0.991630690741981E+07  -0.185949305922308E+07   0.599066486689836E+07
-0.159392396237307E+05   0.296004865275814E+04  -0.110656935176569E+05
 0.800970026849796E+07   0.247717810054366E+07  -0.460543580680696E+08
-0.237010458689252E+07   0.101013983339062E+06   0.203042960322874E+07
 0.379010599355267E+08   0.732341494213578E+07  -0.556606156734835E+08
-0.337010820273821E+08  -0.147919027644202E+08   0.567424911007837E+09
 0.322219416256417E+05  -0.104092072189767E+05   0.128967058686204E+05
-0.257258805870567E+09  -0.635149929624336E+08  -0.152752854956514E+10
 0.214419090344474E+07  -0.247435106210237E+06  -0.499321746092534E+07
-0.449275591851490E+08  -0.902980378929272E+07   0.132124828143333E+09
-0.213133920223355E+08   0.724969796807399E+07  -0.162359994620983E+10
 0.157076702372204E+10   0.487093452892595E+09   0.792526849882218E+10
-0.232153488525298E+05   0.138088690964946E+05   0.466767581292985E+04
-0.395769398304473E+10  -0.160162178614234E+10  -0.128050778279459E+11];
f=C(i,j);
end

function C=C_EBOF(flg1,flg2)
switch flg1
    case -1 % linear
        switch flg2
            case 1 % Smooth ortho - Cintra and Tucker
                C =[-0.15  1.15 -0.10 ;
                    -0.15  0.15  0.90 ;
                     0.60 -0.60 -0.60];
            case 2 % General ortho - Kuzmin
                C=1/7*[ -3/5  6   0;
                        -3/5  0   6;
                        27/5 -6  -6];
        end
    case -2 % quadratic
        %
        switch flg2
            case 1 % fitted ortho - Cintra & Tucker
                C=[
                    0.060964  0.371243 -0.369160 0.555301 0.371218 0.318266 ;
                    0.124711 -0.389402  0.086169 0.258844 0.544992 0.796080 ;
                    1.228982 -2.054116 -2.260574 0.821548 1.819756 1.053907];
            case 2 % general ortho - Cintra and Tucker
                C=[0   0   0 1   0   0 ;
                   0   0   0 0   0   1 ;
                   1  -2  -2 1   2   1];
            case 3 % natural ortho - exact midpoint fit -Kuzmin
                C=[0.0708  0.3236  -0.3776 0.6056  0.4124  0.3068;
                   0.0708 -0.2792   0.2252 0.2084  0.4124  0.7040;
                   1.1880 -2.0136  -2.1264 0.8256  1.7640  0.9384];
            case 4 % Improved ortho fitted ORW Chung et al 2002
```

```
            C=[0.070055  0.339376 -0.396796 0.590331 0.411944 0.333693 ;
               0.115177 -0.368267  0.094820 0.252880 0.535224 0.800181 ;
               1.249811 -2.148297 -2.290157 0.898521 1.934914 1.044147];
        end
    case -3 % Cubic
        switch flg2
            case 1 % Natural ortho - extended quad. fit
                C=[0  1/2     0 1/2  -3/5   0   0  3/5 3/5 0;
                   0    0   1/2   0  -3/5 1/2   0  3/5 3/5 0;
                   1 -3/2  -3/2 1/2   2/5 1/2   0  3/5 3/5 0];
            case {2 3} % Rational ellipsoid
                switch flg2
                    case 2 % Wertzel
                        C=[
                             0.1433751825  0.1433751825  0.9685744898
                            -0.6566650339 -0.5209453949 -2.5526857671
                            -0.5106016916 -0.6463213306 -2.5756669706
                             4.4349137241  2.3303190917  4.4520903005
                             3.5295952199  0.6031924921  2.2044050704
                             0.1229618909  5.1539592511  2.2485545147
                            -5.5556896198 -1.6481269200 -1.8811803355
                            -2.8284365891 -5.4494528976 -1.9023485762
                            -2.9144388828 -0.2256222796 -0.6202937932
                             0.2292109036 -3.7461520908 -0.6414620339
                             1.0000000000  1.0000000000  1.0000000000
                             0.7257989503  0.6916858207 -1.2134964928
                             3.0941511876  3.1282643172 -1.2128608265
                            -4.7303686308 -4.7303686308  0.6004510415
                            -1.6239324646 -1.5898193351  0.2393747647
                            -3.1742364608 -3.2083495904  0.2162486576]';
                    case 3 %LAR32 Matthew Mullens
                        C=[
                             0.0876022     0.1568052     1.0724237
                             0.0282056    -0.5778189    -2.8035540
                            -0.4267843    -0.5142809    -2.6615761
                             0.8764691     2.1323050     4.5667285
                             1.2746771     0.6842509     2.3893798
                             0.6020316     3.4548353     2.0975231
                            -1.9189311    -1.6141226    -1.9047047
                            -0.9342913    -4.0052611    -1.7549784
                            -1.0665831    -0.2632371    -0.6582489
                            -0.2628549    -2.2281332    -0.5082827
                             1.0000000     1.0000000     1.0000000
                            -0.2440019     0.3656529    -1.0685125
                            -0.5741509     1.3857255    -0.7713565
                            -0.8952261    -2.8663578     0.2069083
                            -0.4320974    -1.3596872     0.0673869
                            -0.4627095    -1.5189962    -0.2489999]';

                end
                pt=[1 2 3 5 4 6 9 7 8 10 11 12 13 15 14 16];
                C=C(:,pt);
            case 4 % Improved Orthotropic fitted - ORW3 - Chung & kwon
                C=[ -0.1480648093 -0.2106349673  0.4868019601
```

```
                0.8084618453  0.9092350296  0.5776328438
                0.3722003446 -1.2840654776 -2.2462007509
                0.7765597096  1.1104441966  0.4605743789
               -1.3431772379  0.1260059291    -1.9088154281
               -1.7366749542 -2.5375632310    -4.8900459209
                0.8895946393  1.9988098293     4.0544348937
                1.7367571741  1.4863151577     3.8542602127
               -0.0324756095  0.5856304774     1.1817992322
                0.6631716575 -0.0756740034     0.9512305286]';
           pt=[1 2 4 3 6 5 9 7 8 10]; C=C(:,pt);
    end
case -4 % quartic
    switch flg2
        case 1 % Verweyst
            C=[
                0.6363     0.6363     2.7405;
               -1.8727    -3.3153    -9.1220;
               -4.4797    -3.0371   -12.2571;
               11.9590    11.8273    34.3199;
                3.8446     6.8815    13.8295;
               11.3421     8.4368    25.8685;
              -10.9583   -15.9121   -37.7029;
              -20.7278   -15.1516   -50.2756;
               -2.1162    -6.4873   -10.8802;
              -12.3876    -8.6389   -26.9637;
                9.8160     9.3252    27.3347;
                3.4790     7.7468    15.2651;
               11.7493     7.4815    26.1135;
                0.5080     2.2848     3.4321;
                4.8837     3.5977    10.6117]';
        case 2 % FFLAR4 - Matthew Mullens
            C=[
                  0.678225884   0.748226727   3.167356369
                 -3.834359034  -4.249612053 -13.288266400
                 -2.664862865  -2.987266447 -11.680179330
                 14.209962670  14.938209410  43.700607680
                  9.746185193   8.641488072  23.788431340
                  2.700369681   5.974489008  17.383121430
                -22.447252700 -21.757217160 -58.354308000
                -13.078649640 -15.798676320 -49.513705640
                 -8.013024236  -7.521216405 -19.959054610
                 -0.125467689  -3.616551654 -11.755525930
                 12.689484570  12.640352670  35.425354130
                 10.563248410  10.222185780  25.844317920
                  2.487386515   4.788201652  18.226443930
                  2.417857515   2.376441613   6.291273472
                 -0.328195677   1.056519961   2.925785795]';
        case 3  %LAR4 Ref Matthew Mullens
            C=[
                  0.813175172   1.768619587   4.525066937
                 -3.065410883  -9.826017151 -19.259137620
                 -4.659333003  -6.484058476 -17.650178090
                 14.747639770  28.905936750  61.543979540
                  6.329870878  19.986994700  33.901239610
```

```
                     9.739797775  10.759963010  28.467355970
                   -15.922240910 -40.492387100 -76.738638810
                   -20.818571900 -27.442217500 -68.977583290
                    -4.216519964 -17.715409270 -27.768082700
                    -8.993993112  -7.230748101 -22.399036130
                    11.470974520  19.729631240  43.875135630
                     5.834142985  18.709047480  32.480679940
                     9.874209286   8.882877701  26.928320210
                     1.138888034   5.785725498      8.600822308
                     3.100457733   2.224834058      7.101978254]';

        end
        pt=[1 2 3 5 4 6 9 7 8 10 14 12 11 13 15];
        C=C(:,pt);
end
end
%
```

*Published with MATLAB® R2024b*